\documentclass[12pt]{article}
\pdfoutput=1
\usepackage{amssymb}
\usepackage{amsmath}
\usepackage{amsthm}
\usepackage{color}
\usepackage{comment}
\usepackage{geometry}		
\usepackage{graphicx}

\let\originalleft\left
\let\originalright\right
\renewcommand{\left}{\mathopen{}\mathclose\bgroup\originalleft}
\renewcommand{\right}{\aftergroup\egroup\originalright}

\usepackage{caption}
\usepackage{subcaption}

\usepackage{hyperref}
\hypersetup{
    colorlinks=true,
    linkcolor=black,
    urlcolor=blue,
    citecolor=black
}

\begin{document}

\newcommand\cO{\mathcal{O}}
\newcommand{\ee}{\varepsilon}

\newcommand{\removableFootnote}[1]{\footnote{#1}}

\newtheorem{theorem}{Theorem}[section]
\newtheorem{corollary}[theorem]{Corollary}
\newtheorem{lemma}[theorem]{Lemma}
\newtheorem{proposition}[theorem]{Proposition}

\theoremstyle{definition}
\newtheorem{definition}{Definition}[section]
\newtheorem{example}[definition]{Example}

\theoremstyle{remark}
\newtheorem{remark}{Remark}[section]



\newcommand{\orcid}[1]{\href{https://orcid.org/#1}{\includegraphics[width=8pt]{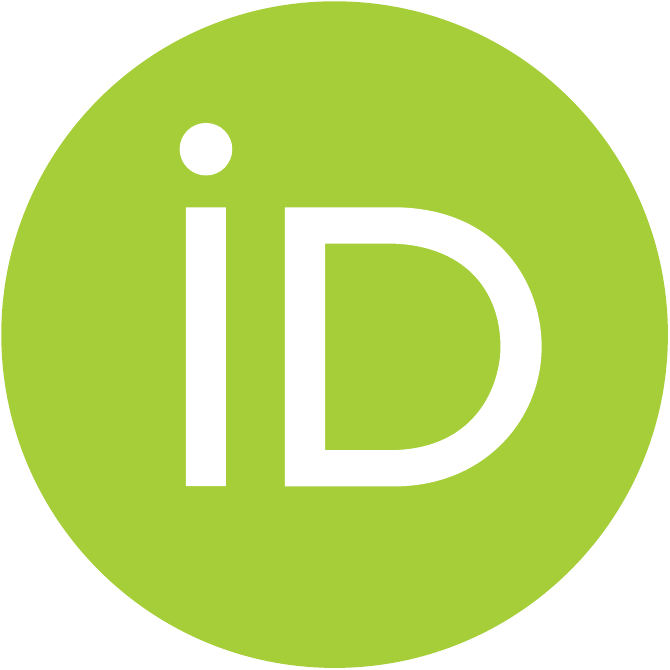}}}

\title{
Pattern Formation in a Spatially-Extended Model of Pacemaker Dynamics in Smooth Muscle Cells
}
\author{
H.O.~Fatoyinbo\orcid{0000-0002-6036-2957}, R.G.~Brown,
D.J.W.~Simpson\orcid{0000-0002-0284-6283}, B.~van Brunt\\\\
School of Fundamental Sciences\\
Massey University\\
Palmerston North, 4410\\
New Zealand\\
 }
 
 
\maketitle


\begin{abstract}
Spatiotemporal patterns are common in biological systems. For electrically-coupled cells previous studies of pattern formation have mainly used external forcing as the main bifurcation parameter. The purpose of this paper is to show that spatiotemporal patterns in electrically-coupled smooth muscle cells occur even in the absence of forcing. We study a reaction-diffusion system with the Morris-Lecar equations and observe a wide range of spatiotemporal patterns for different values of the model parameters. Some aspects of these patterns are explained via a bifurcation analysis of the system without coupling --- in particular Type I and Type II excitability both occur. We show the patterns are not due to a Turing instability and use travelling wave coordinates to analyse travelling waves.
\end{abstract}

\section{Introduction}
\label{sec:intro}
\setcounter{equation}{0}

Smooth muscle cells (SMCs) are widely spread across organs and tubes where they provide a variety of functions in the body. The contraction and relaxation of SMCs regulates organ function, such as the flow rate of blood vessels \cite{Lamboley2003,shaikh}. SMCs help with digestion and nutrient collection in the gastrointestinal tract \cite{Bitar2003, Harnett2005}, and regulate bronchiole diameter in the respiratory system \cite{Chung2000}. In the urinary system, they play a role in removing toxins and electrolyte balance \cite{Alexander73,Andersson2004}. Like other excitable cells (e.g.~neuron, endocrine, and skeletal cells), when stimulated SMCs can generate a large electrical signal (action potential), and contract in response. This process is known as electro-mechanical coupling.

Electro-mechanical coupling in the cell membrane of a SMC is mediated by the influx of extracellular $\textrm{Ca}^{2+}$ through voltage-gated $\textrm{Ca}^{2+}$ channels and $\textrm{Ca}^{2+}$ release from the cell's internal $\textrm{Ca}^{2+}$ store, the sacroplasmic reticulum. The elevation of the intracellular $\textrm{Ca}^{2+}$ concentration causes the membrane potential to increase rapidly, hence the cell membrane is depolarised, and this results in the opening of the $\textrm{K}^{+}$ channels. The efflux of $\textrm{K}^{+}$ then leads to the repolarisation of the cell. The repetition of this activity results in periodic oscillations that elicit vasomotion, that is, the contraction and relaxation of the vessel's cell wall. 

Oscillations are driven by applied current \cite{Hodgkin1952}, agonists \cite{Sneyd1995, koenigsberger2005}, temperature \cite{Anatoly2013,Fillafer2013}, and pressure \cite{Kubanek2018}. Several experimental studies have investigated the electrical activity induced by external stimuli in excitable cells \cite{FRIEL1995,Latchoumane2018,Liang2019}; see also \cite{Roth,FARHY2004,koenigsberger2005,Izhikevich2007DynamicalBursting} for computational studies. 

\begin{figure}[htb]
   \centering
   \includegraphics[scale=0.35]{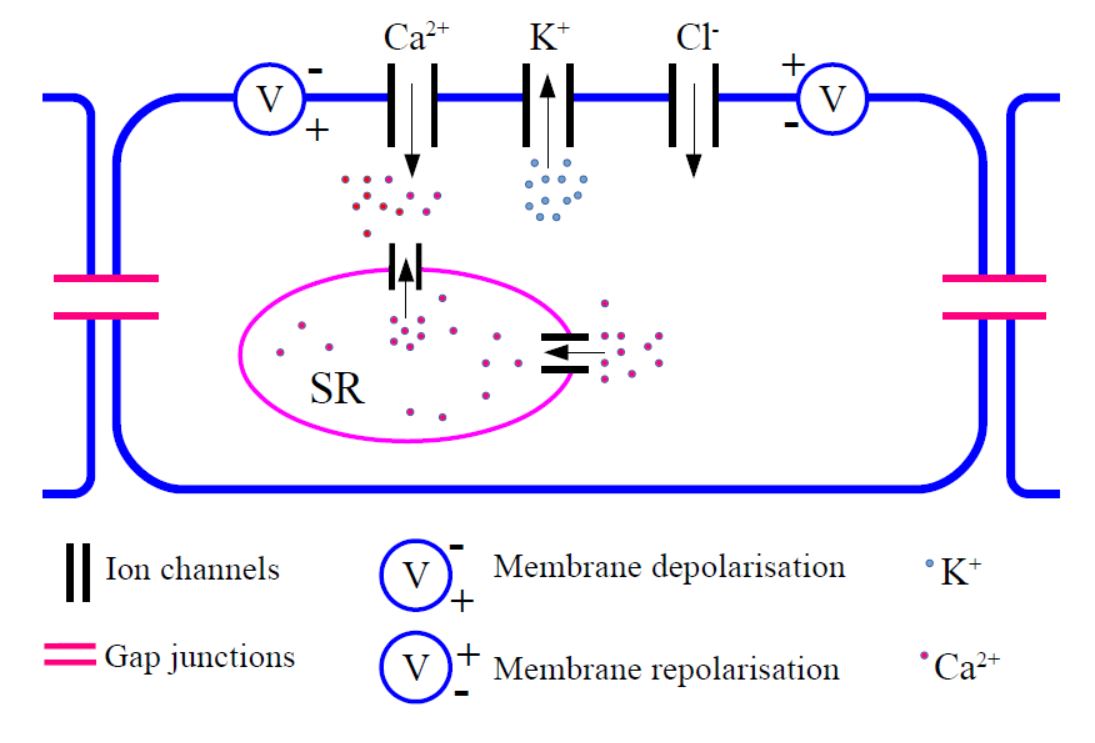}
   \caption{A schematic representation of coupled smooth muscle cells.}
    \label{fig:SMC_Schematic}
\end{figure}

Communication between cells, primarily excitable and non-excitable cells, helps regulate a wide range of cellular activities, for example, receiving and transmitting of signals in the central nervous system \cite{Duan2008Two-parameterModel, JamesKeener2009}, the release of hormones into extracellular fluid in endocrine cells \cite{Schwartz2000, Nakanishi2006, Yves2020}, and contractile activity in muscles \cite{Mege1994,Matchkov2010, Tirziu2010, Bian2015}. Cells are connected to their immediate neighbors through different mechanisms \cite{Jongsma2000,Giepmans2004,Shimizu2013}. SMCs are coupled through gap junctions which can be one of three types: $\textrm{Ca}^{2+}$, inositol triphospate ($\textrm{IP}_{3}$), or membrane potential (electrical) \cite{koenigsberger2004,Haddock05,koenigsberger2005}. A schematic representation of electrically coupled SMCs is shown in Fig.~\ref{fig:SMC_Schematic}. Gap junctional communications have been observed in other cell types, including germ cells in testis \cite{Decrouy2004}, fibroblasts \cite{Azzam2001}, and astrocytes \cite{Giaume1996}.

The dynamics across a large number of coupled cells can form of simple travelling waves, or complex spatiotemporal patterns. For example, as revealed in experiments, spiral waves during heart contractions can cause cardiac arrhythmia \cite{Hwang2005,pandit2013}. Epileptic seizures in the cortex and hallucinations in the retina or visual cortex can be induced by travelling waves \cite{Traub93,Huang2004,Pinto2005,Pearce2015}. 
 
The dynamical behaviour of a single cell is often modelled by a set of ordinary differential equations. Many models have been used and many are strongly related to that given by \cite{Hodgkin1952}. The communication between a large number of cells can be modelled by incorporating spatial-dependence and using a diffusion term. The reaction-diffusion framework for spatiotemporal pattern formation originates with Turing's seminal paper \cite{Turing1952a}. Turing used a reaction-diffusion system to explain how diffusion-driven instability in chemical reactions can result in what are now known as {\em Turing patterns}. Pattern formation has subsequently been studied extensively in diverse applications. In ecology, the Lotka-Volterra model for two interacting species exhibits both Turing and non-Turing patterns when a diffusion term is added \cite{Banerjee2012,shi,Liu2020}. In epidemiology, spatial patterns have been observed in diffusive epidemic models designed to investigate the spread and control of infectious diseases \cite{JiaYF,Chang2020}. Various spatiotemporal patterns have been observed in cellular dynamics due to electrophysiological processes in cells and tissues \cite{Izhikevich2007DynamicalBursting, ramos,kaper,Vo}. Also, spatial patterns have been generated in a variety of physical and mechanical systems \cite{Paul,Perez-Londono2010Zero-HopfStability, Hens}.

Many studies have been published on spatiotemporal patterns in networks of excitable cells \cite{Fujii,Hartle17,Keplinger2014,Lafranceschina2014,Mondal2018DynamicsSolutions,Calim}. Tsyganov {\em et.~al.}~(\cite{Tsyganov}) examined the Fitzhugh-Nagumo model with piecewise-linear reaction terms to study spatiotemporal behaviour in neurons. They found complex dynamical behaviour including the collision and reflection of excitation waves. Meier {\em et.~al.}~(\cite{Meier}) studied a one-dimensional (1D) Morris-Lecar reaction-diffusion system to investigate complex spatiotemporal formation in a network of neurons. Propagation of excitable waves and the spatiotemporal dynamics of excitable neuronal populations in 1D and 2D using Morris-Lecar model were explored by Mondal et al. (\cite{Mondal2019}). Zhu and Liu (\cite{Zhu-Liu}) studied a model with time delay between connected neurons and observed that the spatiotemporal dynamics depends critically on the bifurcation structure of individual neurons. Also, Ali {\em et.~al.}~(\cite{Rehman}) investigated pattern formation in a spatially-extended Wilson-Cowan system. However, most of these studies focused on patterns and waves that are driven by an external stimulus. The purpose of this paper is to stress that an external stimulus is not necessary for spatiotemporal patterns to occur. {\em Pacemaker dynamics} refers to spontaneous oscillations in cells. Experimental studies have shown that pacemaker dynamics occurs in many types of SMCs, such as the gastrointestinal tract, urinary tract, lymphatic vessels, arteries, and veins \cite{Tomita,Hashitani,Fukuta2002,VanHelden,McHale2006OriginMuscle}. Also, there have been several computational studies of pacemaker dynamics in muscle cells \cite{Youm2006,Rihana2009,Cho2012,Ho2016}.

In this paper we focus on pacemaker electro-mechanical coupling
activity in arterial SMCs due to changes in the vessel's transmural pressure, that is, the pressure gradient across the vessel wall. We study a spatially-extended two-variable nondimensionalised Morris-Lecar model with no applied current. As shown in our previous work \cite{hammed}, without diffusion the model is a reduced form of the three-dimensional ODE model of Gonzalez-Fernandez and Ermentrout (\cite{Gonzalez-Fernandez1994}) for the dynamics of pacemaker vasomotion in SMCs of small arteries. 
 
In Section 2 we state the model equations. In Section 3 we summarise the dynamics of the model without diffusion using the voltage associated with the ${\rm K}^{+}$ and ${\rm Ca}^{2+}$ channels as bifurcation parameters. Oscillations can arise via both Type I and Type II excitability. This distinction of two types of excitability was first described by \cite{Hodgkin1948TheAxon}.  For Type I excitability oscillations arise via a saddle-node on invariant circle (SNIC) bifurcation, whereas for Type II excitability oscillations arise via a Hopf bifurcation \cite{Rinzel1999AnalysisNetwork}. In Section \ref{sec:4} we show that the spatiotemporal patterns that emerge are non-Turing patterns due to violation of Turing's instability criteria. Numerical simulations of the reaction-diffusion model is carried out in Section \ref{sec:5}. Various spatiotemporal patterns including travelling pulses and fronts are explored. The existence of the travelling waves is analysed in Section \ref{sec:6}. In the final section some conclusions are presented.


\section{A nondimensionalised Morris-Lecar system with diffusion}\label{sec:2}
\setcounter{equation}{0}

We consider a nondimensionalised reaction-diffusion system to model the dynamics of a population of coupled SMCs through passive electrical coupling of adjacent cells. The reaction term in the model is based on our previous study on an isolated SMC \cite{hammed}. The model equations are 
\begin{align}
\label{eq:dimless1}
\frac{\partial V}{\partial \tau}&=D\frac{\partial^{2} V}{\partial X^{2}}-\bar{g}_{L}(V-\bar{v}_{L})-\bar{g}_{K}N(V-\bar{v}_{K})-\bar{g}_{\textrm{Ca}}M_{\infty}(V)(V-1), \\
\frac{\partial N}{\partial \tau}&=\lambda(V)(N_{\infty}(V)-N),
\label{eq:dimless2}
\end{align}
where $V(X,\tau)$ is the membrane potential and  $N(X,\tau)$ is the fraction of open $\textrm{K}^{+}$ channels. The system parameter $D \ge 0$ is the diffusion coefficient, $\bar{g}_{\rm L}$, $\bar{g}_{\rm K}$, and $\bar{g}_{\rm Ca}$ are conductances per unit area for the leak, potassium, and calcium currents respectively, while $\bar{v}_{\rm L}$ and $\bar{v}_{\rm K}$ are the corresponding Nernst reversal potentials (equilibrium potentials). The fraction of open calcium [potassium] channels at steady state $M_{\infty}$ [$N_\infty$] and the time scale for the opening of the potassium channel, $\lambda(V)$ are:
\begin{align}
\label{MINF}
M_{\infty}(V)&=\frac{1}2{}\left(1+\tanh\left(\frac{V-\bar{v}_{1}}{\bar{v}_{2}}\right)\right),\\
N_{\infty}(V)&=\frac{1}{2}\left(1+\tanh\left(\frac{V-\bar{v}_{3}}{\bar{v}_{4}}\right)\right),\\
\lambda(V)&=\psi\cosh\left(\frac{V-\bar{v}_{3}}{2\bar{v}_{4}}\right),
\label{eq:lambda}
\end{align}
where $\bar{v}_1$ and $\bar{v}_3$ measure the potential at which potassium and calcium channels are half-opened, $\psi$ is a time constant, and $\bar{v}_2$ and $\bar{v}_4$ are additional parameters. as listed in \cite{hammed}: $\bar{v}_{1}=-0.2813$, $\bar{v}_{2}=0.3125$, $\bar{v}_{3}=-0.1380$, $\bar{v}_{4}=0.1812$, $\psi=0.1665$, $\bar{v}_{\rm L}=-0.875$, $\bar{v}_{\rm K}=-1.125$, $\bar{g}_{\rm L}=0.25$, $\bar{g}_{\rm K}=1.0$, and $\bar{g}_{\rm Ca}=0.4997$.

In this paper we consider a one-dimensional spatial domain $\Omega = [-L,L]$ for the values of $X$. At the boundaries, $X = \pm L$, we use no-flux boundary conditions, however we are primarily concerned with the spatiotemporal patterns that emerge away from the boundaries.
\section{The dynamics of a single cell}\label{sec:3}
\setcounter{equation}{0}
In this section we summarise the dynamics of \eqref{eq:dimless1}--\eqref{eq:dimless2} in the absence of diffusion, i.e.~$D=0$. We show how stable oscillations are created through either Type I or Type II excitability. This is important to the nature of the spatiotemporal dynamics described in Section~\ref{sec:5}. More details on the dynamics of a single cell can be found in \cite{hammed,mythesis}. 

\begin{figure}[htbp]
\centering
\begin{subfigure}{.4\textwidth}
  \centering
  \caption{}
  \includegraphics[width = \textwidth]{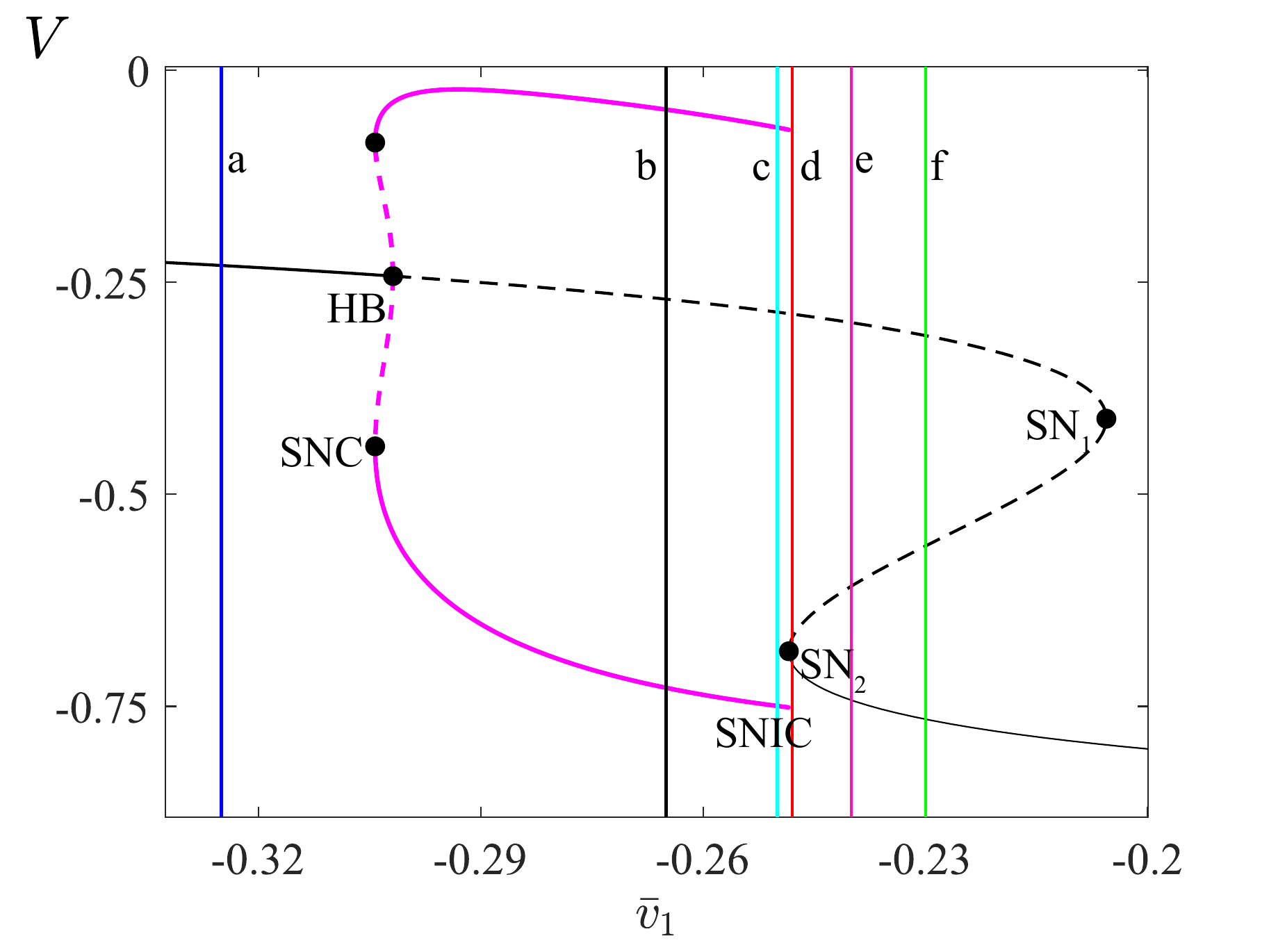}
  \label{fig:v1bif}
\end{subfigure}\\%
\begin{subfigure}{.4\linewidth}
    \centering
    \caption{}
    \includegraphics[width=\textwidth]{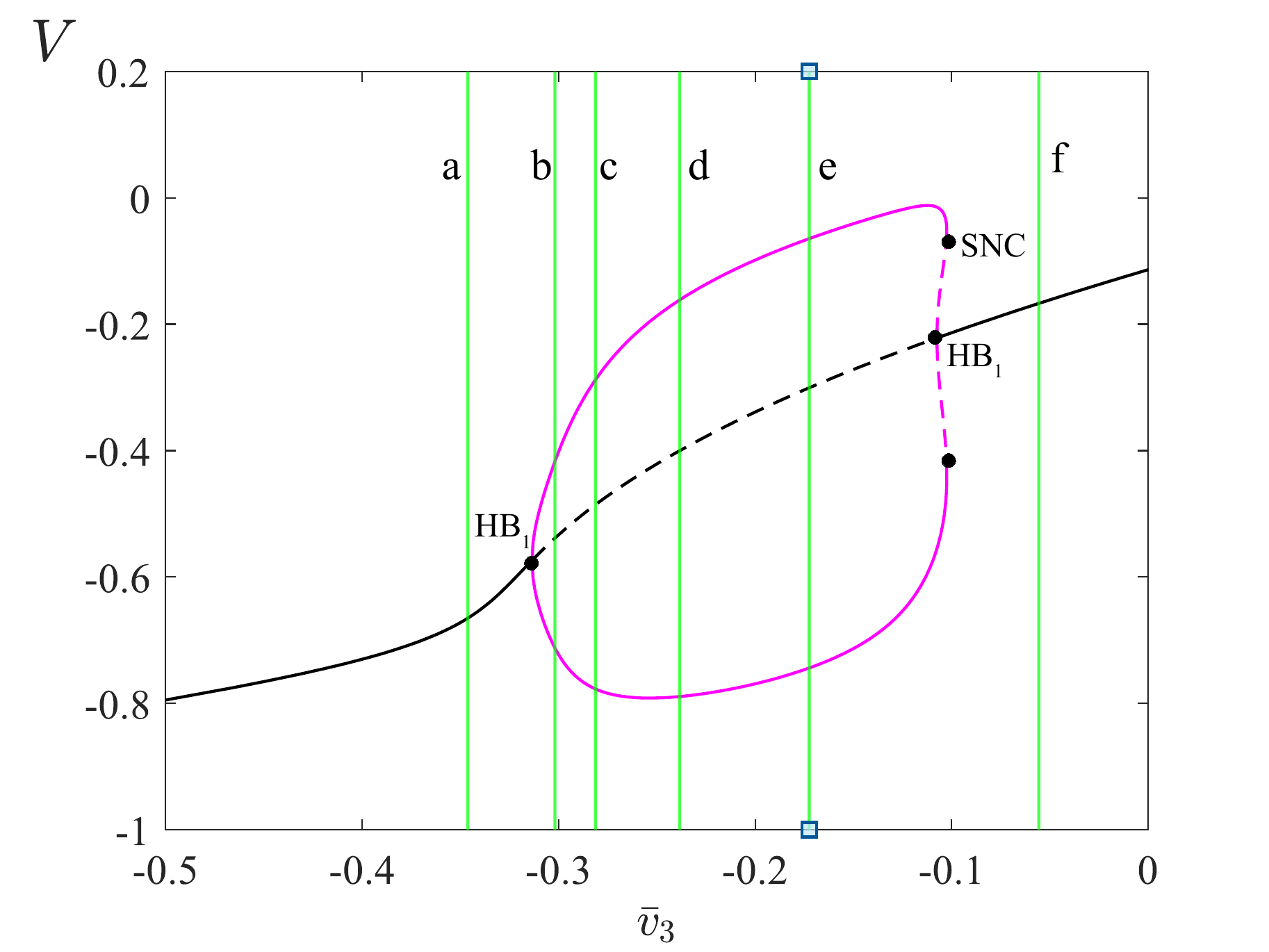}
     \label{fig:v3bif}
  \end{subfigure}%
\begin{subfigure}{.4\textwidth}
  \centering
  \caption{}
 \includegraphics[width = \textwidth]{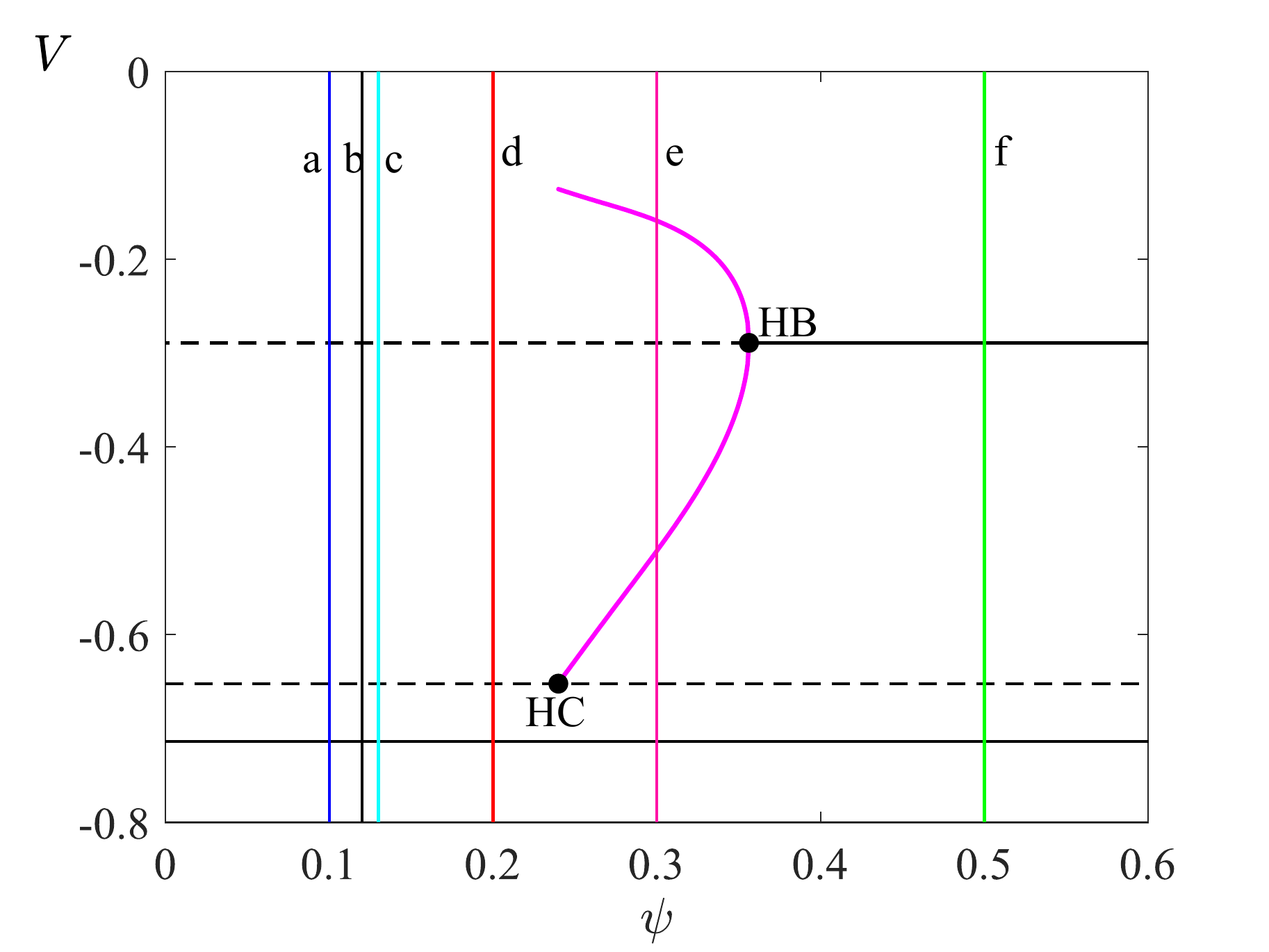}
  \label{fig:psi}
\end{subfigure}%
\caption{Bifurcation diagrams of \eqref{eq:dimless1}--\eqref{eq:dimless2} with $D=0$ using $\bar{v}_1$, $\bar{v}_3$, and $\psi$ as the bifurcation parameters. In each diagram all parameters (except the one being varied) are fixed at the values listed in Section~\ref{sec:2}. Black [magenta] curves correspond to equilibria [limit cycles]. Solid [dashed] curves correspond to stable [unstable] solutions. The vertical lines indicate the parameter values used in Figs.~\ref{fig:SPT_v1}, \ref{fig:SPT_v3}, and \ref{fig:SPT_psi}. HB: Hopf bifurcation; SN: saddle-node bifurcation (of an equilibrium); SNC: saddle-node bifurcation of a limit cycle; SNIC: saddle-node on an invariant circle bifurcation; HC: homoclinic bifurcation.}
\label{fig:onepar}
\end{figure}

Figs.~\ref{fig:onepar}(a)--(c) show bifurcation diagrams as $\bar{v}_{1}$, $\bar{v}_{3}$, and $\psi$ are varied from their values as listed in Section~\ref{sec:2}. These were computed numerically using {\sc auto} \cite{Doedel2012}. Spontaneous oscillations in \eqref{eq:dimless1}--\eqref{eq:dimless2} are triggered by a change in transmural pressure, therefore we use pressure-dependent parameters, $\bar{v}_{1}$ and $\bar{v}_{3}$, as bifurcation parameters. 

 Fig.~3.1a shows the result of varying $\bar{v}_1$.  The system has a unique equilibrium except between saddle-node bifurcations ${\rm SN}_1$ and ${\rm SN}_2$ where there are three equilibria: one stable (lower branch) and two unstable (middle and upper branch). As the value of $\bar{v}_1$ is increased from the smallest value shown in the diagram, the upper equilibrium branch loses stability in a subcritical Hopf bifurcation (HB). The unstable limit cycle produced here gains stability via a saddle-node bifurcation (SNC). Upon further increasing the value of $\bar{v}_{1}$ the stable limit cycle is destroyed at the saddle-node bifurcation ${\rm SN}_2$.  This is an example of a saddle-nodle on invariant circle bifurcation (SNIC) where the limit cycle is replaced by a heteroclinic connection between the two equilibria \cite{KuznetsovElementsTheory}. As a consequence, the period of the limit cycle approaches infinity as the bifurcation is approached. Here the system displays Type I excitability as stable oscillations are created in a SNIC bifurcation by appropriately decreasing the value of $\bar{v}_1$.
 
 Next we vary the value of $\bar{v}_3$. As shown in Fig.~\ref{fig:onepar}b, as we increase the value of $\bar{v}_3$ a unique equilibrium loses stability in a supercritical Hopf bifurcation $\rm{HB}_1$ then regains stability in a subcritical Hopf bifurcation $\rm{HB}_2$. The stable oscillations are created at $\rm{HB}_1$ with finite period. They subsequently lose stability at a saddle-node bifurcation and terminate at $\rm{HB}_2$. In this case the system displays Type II excitability since the periodic oscillations arises through a Hopf bifurcation.

 Finally Fig.~3.1c shows how the dynamics changes under variation to the value of $\psi$.  The system has three equilibria for all values of $\psi>0$. For relatively low and intermediate values of $\psi$, there exist one stable (lower branch) and two unstable (upper and middle branch) equilibria. By increasing $\psi$, a stable limit cycle emanates through a homoclinic bifurcation (HC) and upon further increase of $\psi$ terminates in a supercritical Hopf bifurcation (HB). As in Fig.~3.1b the excitability here is Type I.  Between the homoclinic and Hopf bifurcations the system is bistable as the limit cycle coexists with a stable equilibrium. As shown in \cite{hammed} for different parameter values the system has three coexisting stable solutions.
\setcounter{equation}{0}

\section{Linear stability analysis}\label{sec:4}
\setcounter{equation}{0}
 Alan Turing (\cite{Turing1952a}) hypothesised that spatially inhomogeneous patterns may arise in a reaction-diffusion system if a spatially homogeneous steady state is stable in the absence of diffusion and destabilised as result of diffusion. Such instability is referred to as diffusion-driven instability or Turing instability. The conditions required for the onset of Turing instability have been well studied \cite{Alonso,Shoji,Baurmann,Banerjee2012}. Here we perform a linear stability analysis of \eqref{eq:dimless1}--\eqref{eq:dimless2} around a spatially homogeneous steady state and show that the conditions for Turing instability are not satisfied for this system.

As shown in the previous section, in the absence of diffusion \eqref{eq:dimless1}--\eqref{eq:dimless2} typically has one or three equilibria and the stability of these depends on the values of the parameters.  Here let $(V^*,N^*)$ be a stable equilibrium of \eqref{eq:dimless1}--\eqref{eq:dimless2} with $D=0$ (i.e.~no diffusion) for some combination of parameter values. Then for \eqref{eq:dimless1}--\eqref{eq:dimless2} with $D > 0$, $(V^*,N^*)$ represents a spatially homogeneous state.    

Let $\big(V_{0}(X,\tau),N_{0}(X,\tau)\big)$ represent the perturbation of a solution to \eqref{eq:dimless1}--\eqref{eq:dimless2} from the steady state, i.e. 
\begin{equation}
\label{eq:perturb_space}
    \begin{pmatrix}
   V_{0}\\
    N_{0}
    \end{pmatrix}=\begin{pmatrix} V-
    V^{*}\\
    N-N^{*}
    \end{pmatrix}.
\end{equation}
By linearising \eqref{eq:dimless1}--\eqref{eq:dimless2} about $(V^*,N^*)$ we obtain the following leading-order approximation to the dynamics of the perturbation:
\begin{equation}
    \begin{pmatrix}
  V_{0}\\
    N_{0}
    \end{pmatrix}_{\tau}=\begin{pmatrix}
   D & 0\\
    0 & 0
    \end{pmatrix}\begin{pmatrix}
    V_{0}\\
    N_{0}
    \end{pmatrix}_{XX}+\begin{pmatrix}
   f_{V} & f_{N}\\
    g_{V} & g_{N}
    \end{pmatrix}\begin{pmatrix}
    V_{0}\\
    N_{0}
    \end{pmatrix}.
    \label{eq:linearised_plane}
\end{equation}
The second matrix in \eqref{eq:linearised_plane} is the Jacobian matrix of \eqref{eq:dimless1}--\eqref{eq:dimless2} evaluated at $(V^*,N^*)$.  By directly differentiating \eqref{eq:dimless2} with respect to $N$ we obtain
\begin{align}
\label{eq:gN}
g_{N}& =-\psi \cosh\left(\frac{V^*-\bar{v}_{3}}{2\bar{v}_{4}}\right). 
\end{align}
Formulas for the other three entries in the Jacobian matrix will not be needed. 

We now look for a solution to \eqref{eq:linearised_plane} of the form $(V_{0},N_{0})(X,\tau)=\boldsymbol{\beta} e^{(\lambda\tau+ikX)}$, where $\boldsymbol{\beta}$ is a constant vector, $\lambda$ is the growth rate of perturbation in time, and $k$ is the wave number. While there are many such solutions, we will show that all must have $\lambda<0$.  This implies that for any sufficiently small perturbation \eqref{eq:perturb_space}, the corresponding solution to \eqref{eq:dimless1}--\eqref{eq:dimless2} decays to $(V^*,N^*)$ as $t\to\infty$, hence the steady-state is not destabilised \cite{Murray2003}. By substituting the given form into \eqref{eq:linearised_plane},
\begin{equation}
   \begin{pmatrix}
   -k^2D+f_{V}-\lambda & f_{N}\\
   g_{V} & g_{N}-\lambda
    \end{pmatrix}\boldsymbol{\beta}=\begin{pmatrix} 0 \\ 0 \end{pmatrix}.
    \label{eq:homogeneous1}
\end{equation}
Equation \eqref{eq:homogeneous1} is homogeneous in $\boldsymbol{\beta}$, thus has a nontrivial solution only if the matrix in \eqref{eq:homogeneous1} is singular.

This implies
\begin{equation}
 \label{eq:tracedet}
    \lambda=\frac{T}{2}\pm\frac{\sqrt{T^{2}-4\Delta}}{2},
\end{equation}
where $T=-k^2D+f_{V}+g_{N}$ and $\Delta=-k^2Dg_{N}+f_{V}g_{N}-g_{V}f_{N}$ denote the trace and determinant of the matrix in \eqref{eq:homogeneous1} when $\lambda = 0$. By assumption $(V^*,N^*)$ is stable in the absence of diffusion, therefore
\begin{equation}
\label{eq:stabilitycriteria}
  f_{V}+g_{N}<0,  \qquad f_{V}g_{N}-f_{N}g_{V}>0.
\end{equation}
But from \eqref{eq:gN} we always have $g_{N}<0$ because $\psi>0$ for physical reasons. Therefore $T<0$ and $\Delta>0$, thus $\lambda<0$ for any $D>0$. Thus $(V^*,N^*)$ is not destabilised by the inclusion of diffusion and so the spatiotemporal patterns that we describe below are not due to Turing instability.

\section{Spatiotemporal dynamics of the full model}\label{sec:5}
\setcounter{equation}{0}
In this section we explore the effect of varying $\bar{v}_{1}$, $\bar{v}_{3}$ and $\psi$ on the spatiotemporal dynamics of the reaction-diffusion system \eqref{eq:dimless1}--\eqref{eq:dimless2}. Since the patterns are not due to Turing instability, as shown in Section~\ref{sec:4}, we will investigate spatiotemporal dynamics for a wide range of parameter values, in particular where the steady states may be stable or unstable. We show that a wide range of spatiotemporal patterns can occur, including travelling pulses, travelling fronts, and spatiotemporal chaos.

The system \eqref{eq:dimless1}--\eqref{eq:dimless2} was solved numerically by using the method of lines. We used a second-order central finite difference approximation to the spatial derivative using $1000$ $X$-values per unit interval, and a standard numerical scheme for the time derivative ({\sc ode15s} in {\sc matlab}) \cite{Schiesser2009AModels,Hiptmair2010}. All numerical simulations use no-flux boundary conditions for $X\in[-L,L]$ and initial conditions
\begin{align}
\label{eq:init}
V(0,X)=V^*+G(X) \hspace{2mm} \text{and} \hspace{2mm} N(0,X)=N^*,
\end{align}
where $(V^*,N^*)$ is a homogeneous steady state of \eqref{eq:dimless1}--\eqref{eq:dimless2}.  Different functions $G(X)$ (specified below) provide different perturbations from the steady state. A linear coordinate change can be applied to \eqref{eq:dimless1}--\eqref{eq:dimless2} to scale the value of $D>0$ to any positive number; in all simulations below we use $D=0.0001$.

\subsection{The effect of the parameters \texorpdfstring{$\bar{v}_{1}$, $\bar{v}_{3}$, and $\psi$}{TEXT} }\label{sec:SPTD}
Now we examine the spatiotemporal patterns exhibited by \eqref{eq:dimless1}--\eqref{eq:dimless2} for the values of $\bar{v}_{1}$, $\bar{v}_{3}$, and $\psi$ marked a-f in Figs.~\ref{fig:onepar}a to c. In this initial condition \eqref{eq:init} we use the Gaussian perturbation,
\begin{equation}
\label{eq:gauss}
    G(X)=A_0 \,{\rm exp} \left( \frac{-X^2}{2 \sigma^2} \right),
\end{equation}
with $A_{0}=0.3$ and $\sigma=0.1$.

\begin{figure}[htbp]
\centering
  \begin{subfigure}[b]{.3\linewidth}
    \centering
    \caption{}
    \includegraphics[width=.99\textwidth]{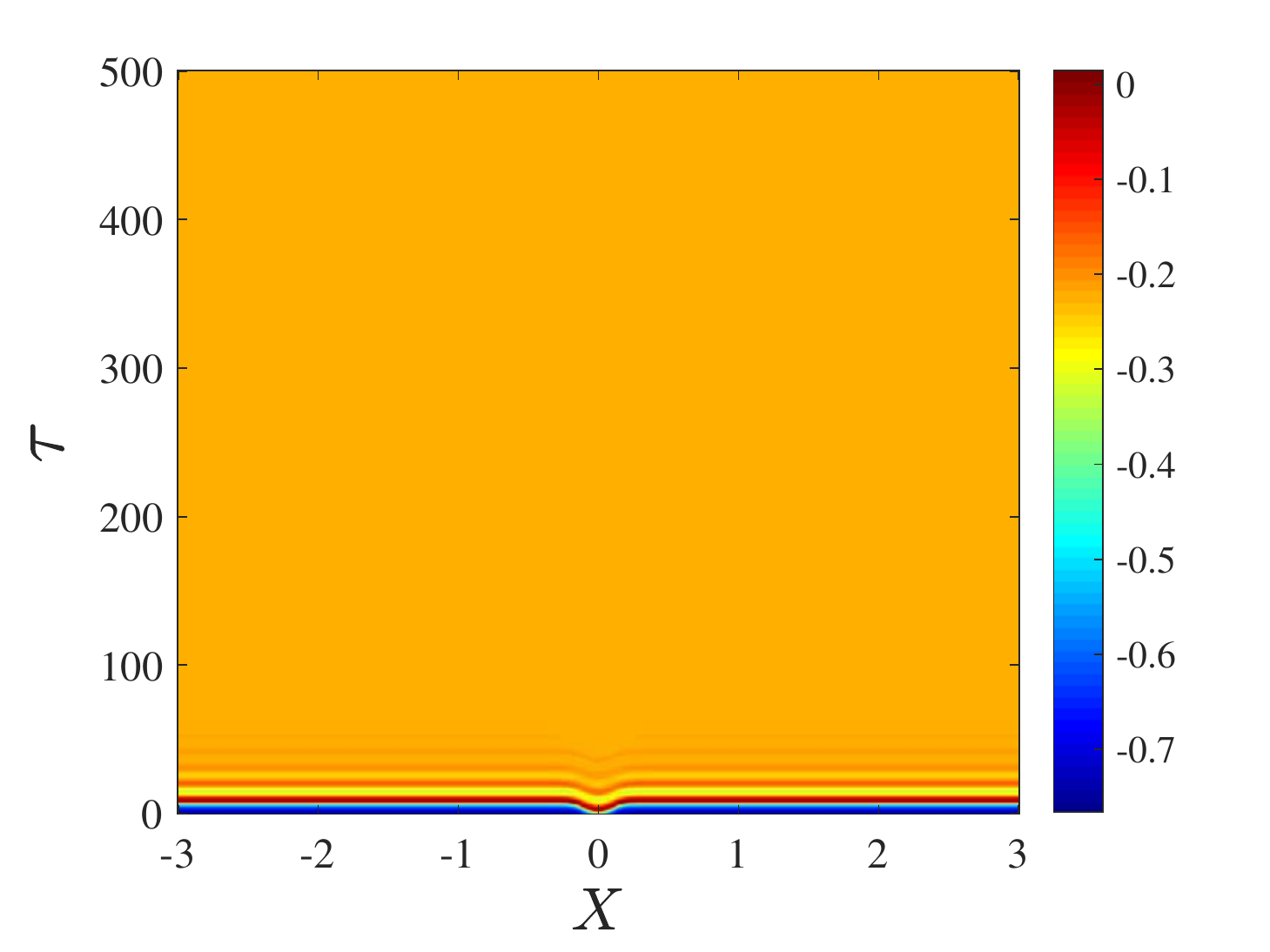}
     \label{fig:sptv10325}
  \end{subfigure}%
   \begin{subfigure}[b]{.3\linewidth}
    \centering
    \caption{}
    \includegraphics[width=.99\textwidth]{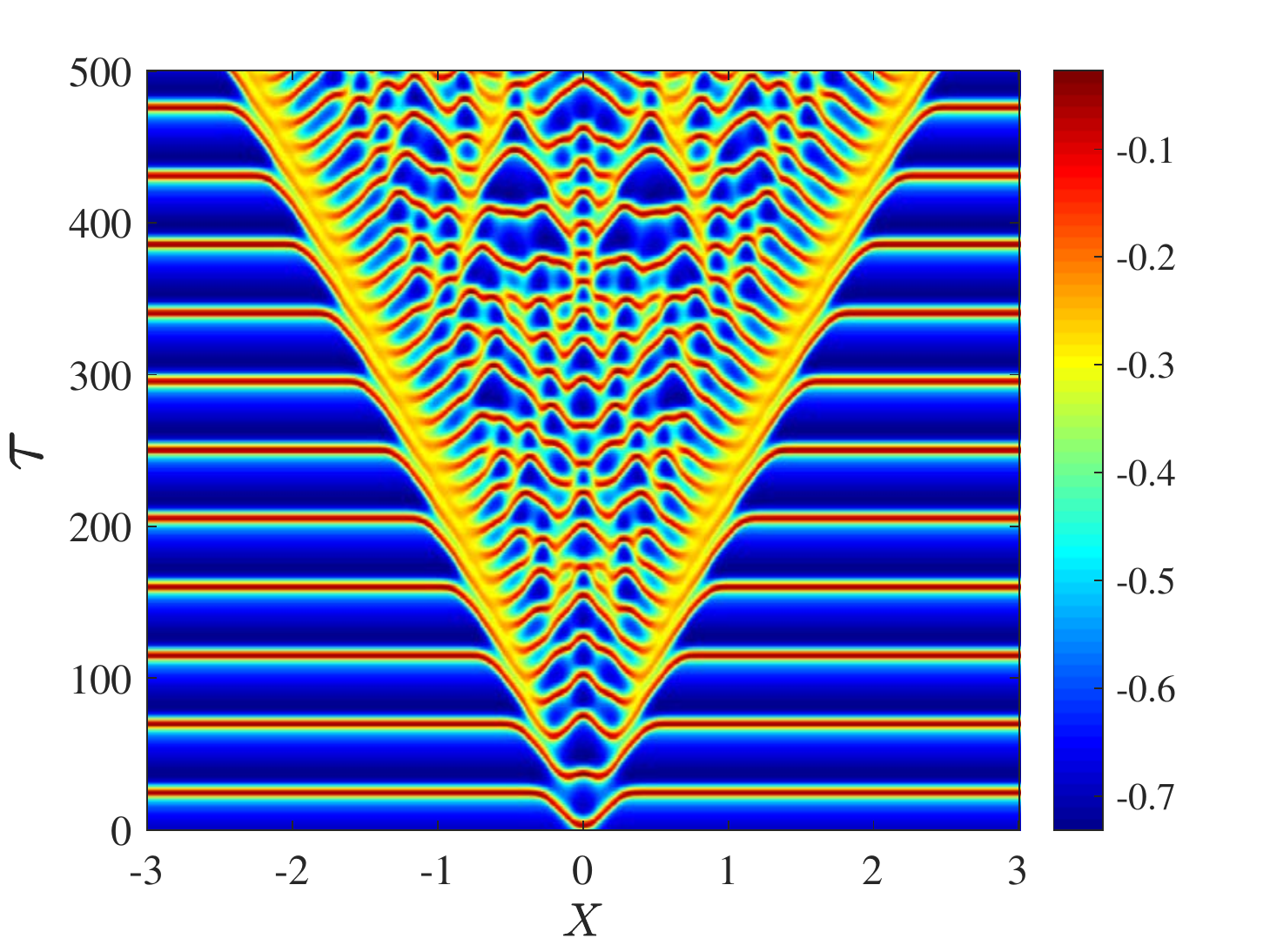}
     \label{fig:sptv10265}
  \end{subfigure}%
  \begin{subfigure}[b]{.3\linewidth}
    \centering
    \caption{}
    \includegraphics[width=.99\textwidth]{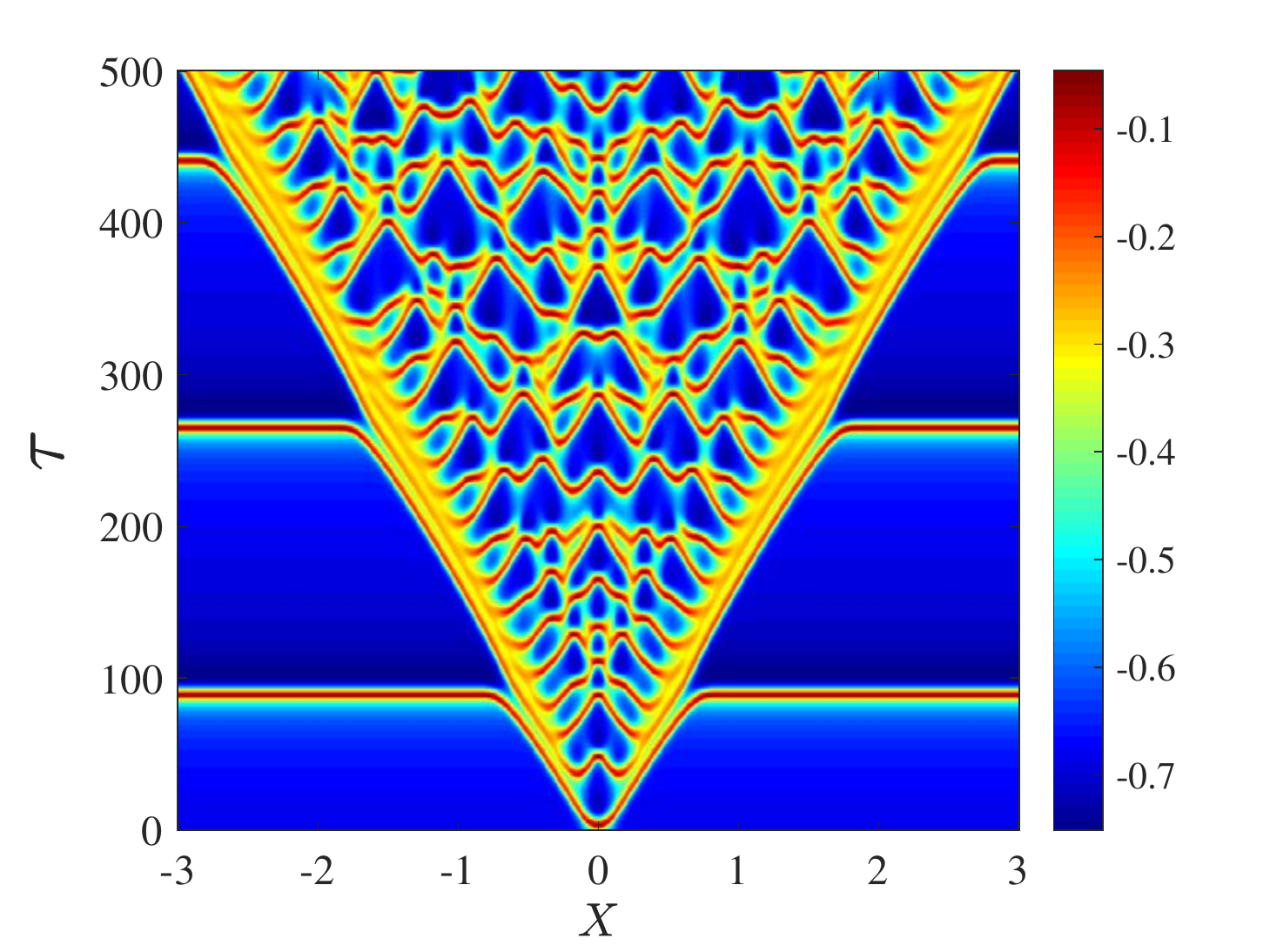}
     \label{fig:sptv1025}
  \end{subfigure}\\%
  \begin{subfigure}[b]{.3\linewidth}
    \centering
    \caption{}
    \includegraphics[width=.99\textwidth]{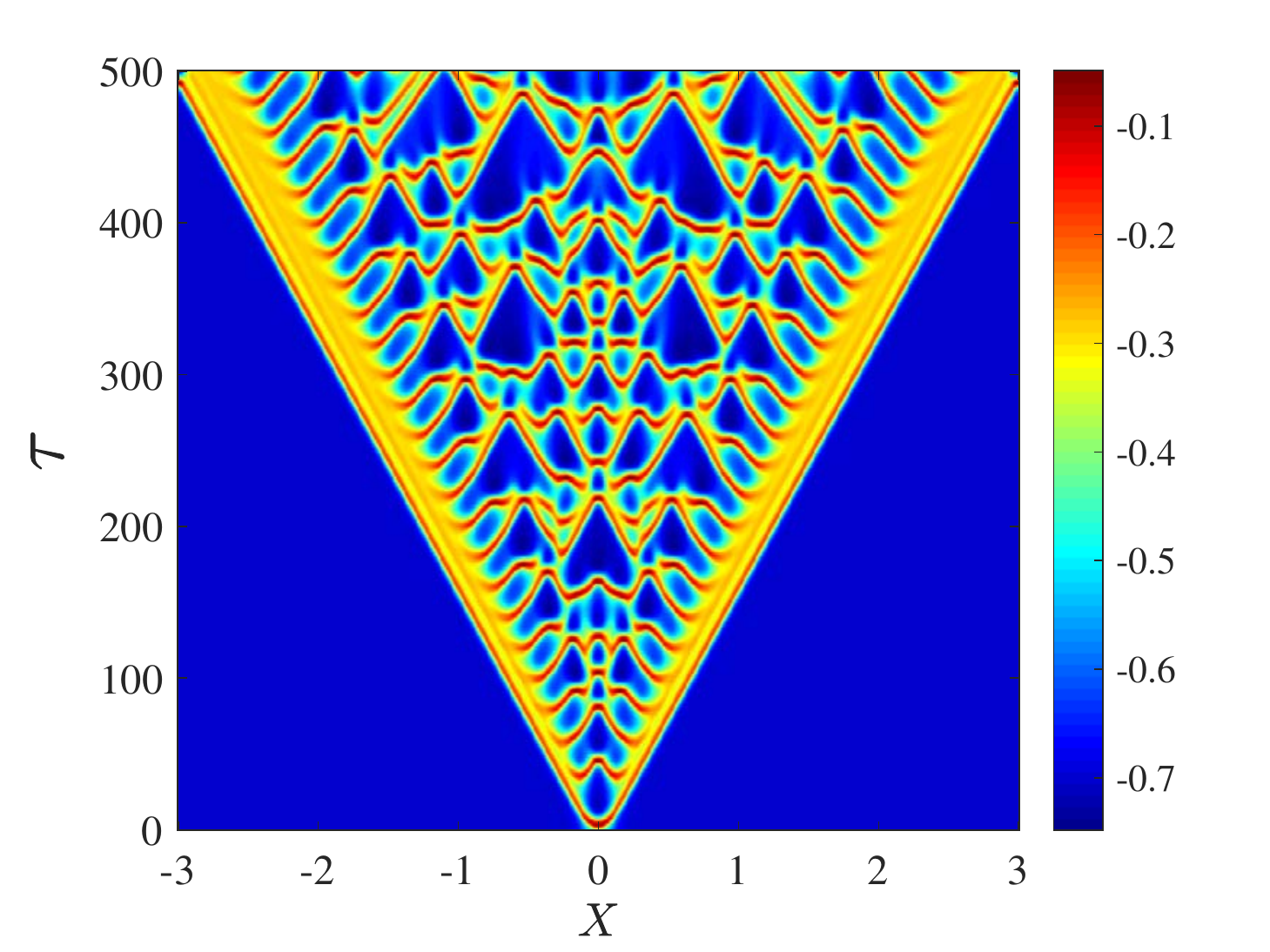}
     \label{fig:sptv10248}
  \end{subfigure}%
   \begin{subfigure}[b]{.3\linewidth}
    \centering
    \caption{}
    \includegraphics[width=.99\textwidth]{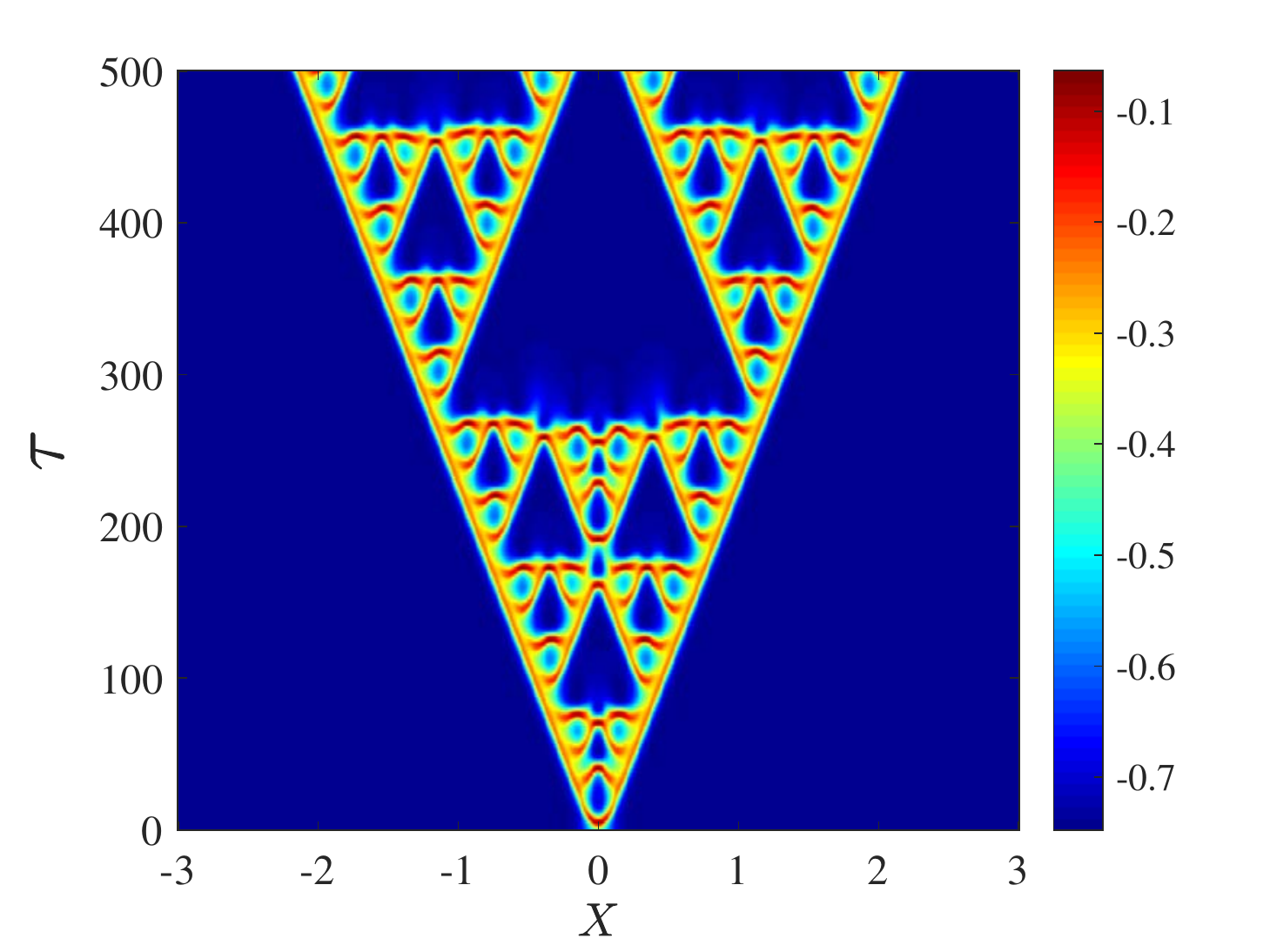}
     \label{fig:sptv1024}
  \end{subfigure}%
  \begin{subfigure}[b]{.3\linewidth}
    \centering
    \caption{}
    \includegraphics[width=.99\textwidth]{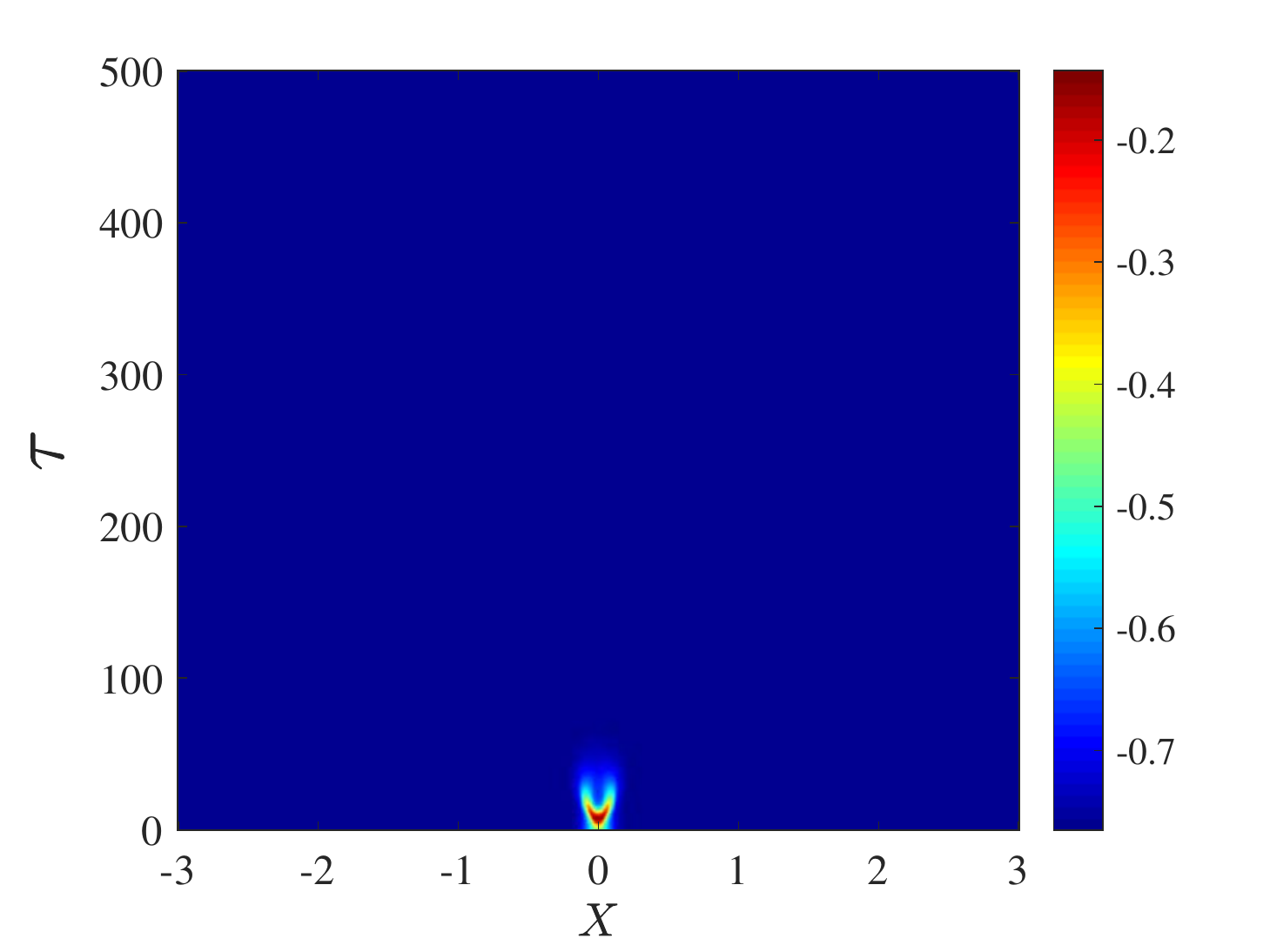}
    \label{fig:sptv1023}
  \end{subfigure}%
 \caption{Space-time plots of the membrane potential $V$ for the values of $\bar{v}_{1}$ marked in Fig.~\ref{fig:v1bif}. Specifically (a) $-0.325$; (b) $-0.265$; (c) $-0.25$; (d) $-0.248$; (e) $-0.240$; and (f) $-0.230$. The initial condition is \eqref{eq:init} with \eqref{eq:gauss}, using the upper equilibrium branch of Fig.~\ref{fig:v1bif} for the steady state $(V^*,N^*)$, and all other parameters are fixed as in Section~\ref{sec:2}.}
 \label{fig:SPT_v1}
\end{figure}

Fig.~\ref{fig:SPT_v1} shows the resulting spatiotemporal patterns for different values of $\bar{v}_{1}$. For low values of $\bar{v}_{1}$the system has a unique homogeneous steady state (the upper equilibrium branch in Fig.~3.1a).  This steady state is stable and the solution quickly converges to the steady state as in Fig.~5.1a.  Instead with $\bar{v}_1$ just to the right of the Hopf bifurcation, a complex spatiotemporal pattern emerges, Fig.~5.1b. The solution starts as a pulse at the centre of the domain due to the initial perturbation. Then the pulse splits into two propagating pulses that transition to time-periodic oscillations with inhomogeneous patterns at the back as they move across the domain. Outside the patterned region the solution is periodic corresponding to the limit cycle of the system with no diffusion. Similar behaviour is observed for values of $\bar{v}_{1}$ between the Hopf bifurcation and the SNIC bifurcation. For example in Fig.~\ref{fig:SPT_v1}c we have used $\bar{v}_1 = 0.25$.  This is very close to the SNIC bifurcation so now the oscillations outside the patterned region are of particularly high period.

\begin{figure}[htb]
\centering
 \begin{subfigure}[b]{.3\linewidth}
   \centering
   \caption{}
\includegraphics[width=.99\textwidth]{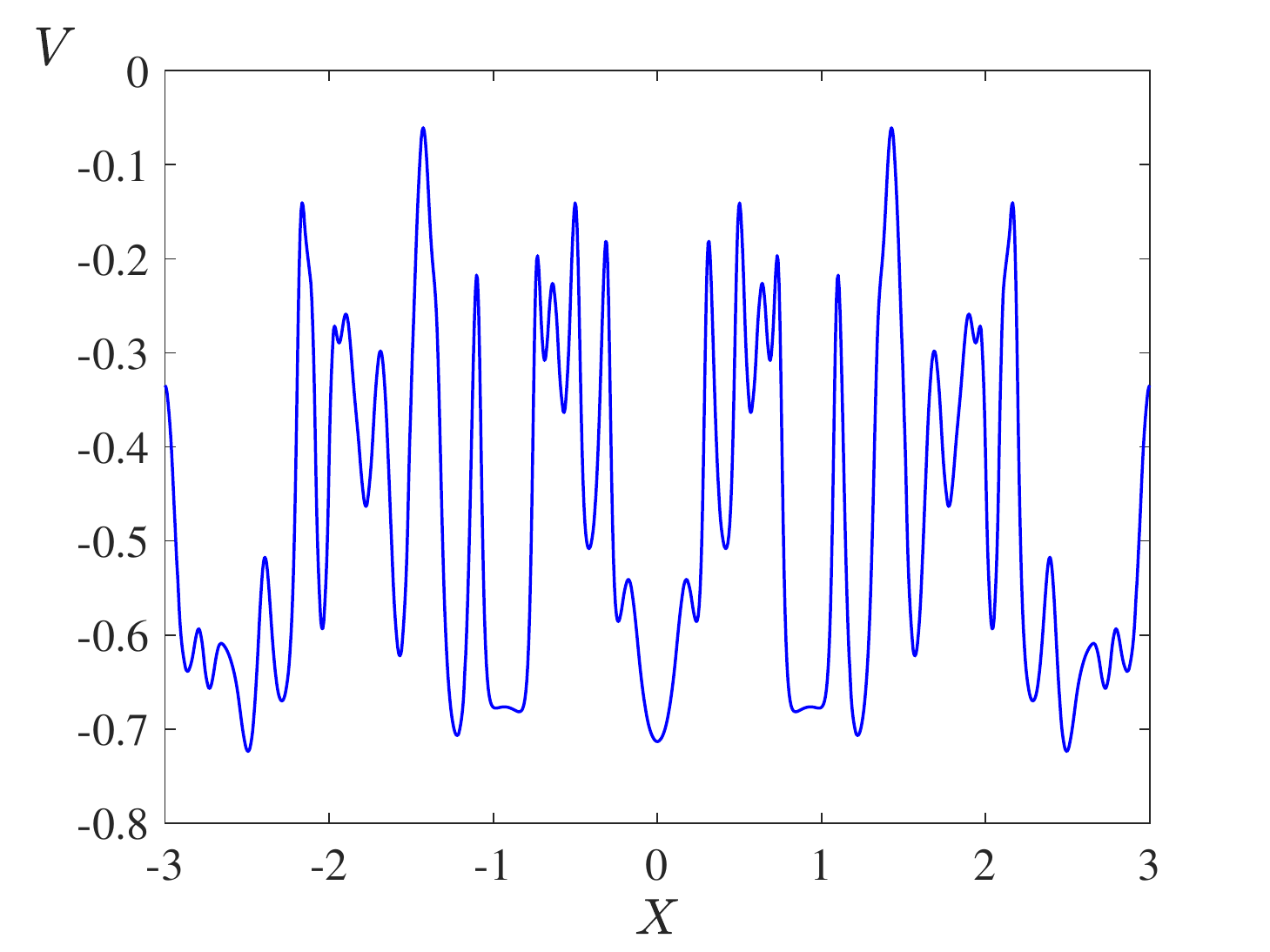}
     \label{fig:V_chaos}
  \end{subfigure}%
  \begin{subfigure}[b]{.3\linewidth}
   \centering
   \caption{}
\includegraphics[width=.99\textwidth]{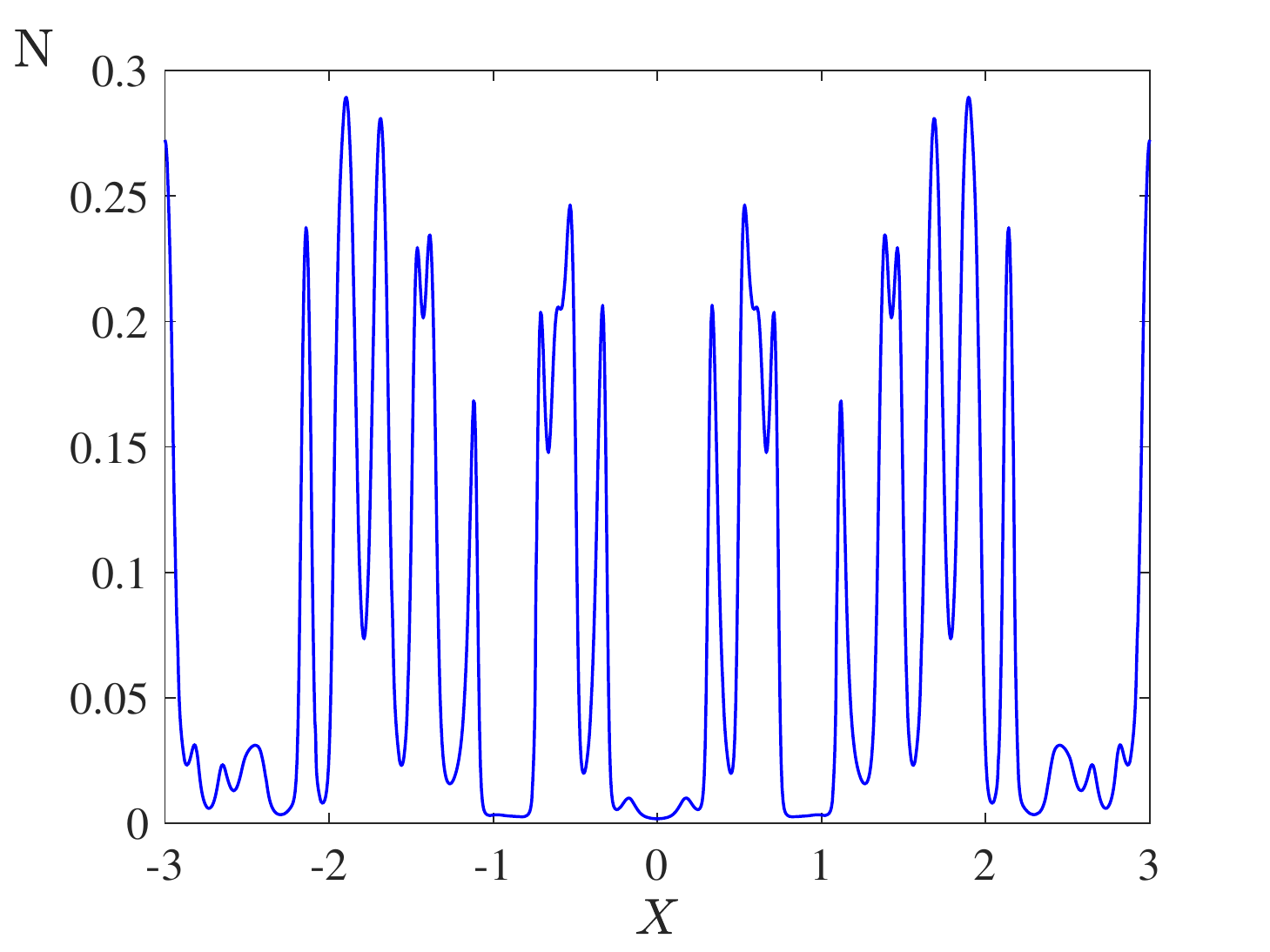}
     \label{fig:N_chaos}
  \end{subfigure}\\%
 \begin{subfigure}[b]{.3\linewidth}
\centering
\caption{}
\includegraphics[width=.99\textwidth]{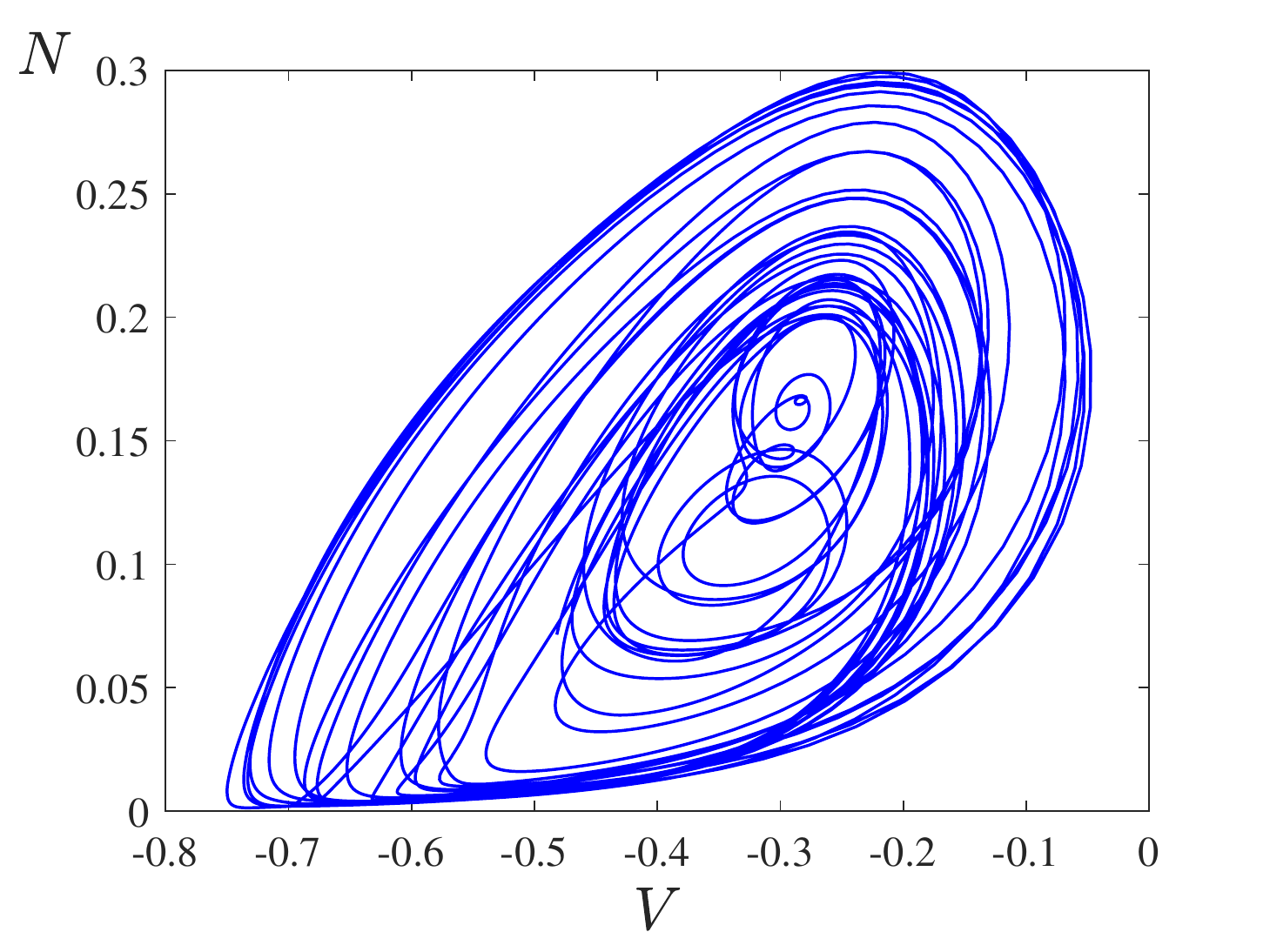}
\label{fig:chaos_pp}
 \end{subfigure}%
   \caption{The spatial distribution of the (a) membrane potential $V$; (b) fraction of open potassium channel $N$; (c) The temporal dynamics for $V$ and $N$ corresponding to patterns in (a) and (b) for time $\tau=1000$ at $\bar{v}_{1}=-0.25$.}
    \label{fig:chaos_pp1}
  \end{figure}

Beyond the SNIC bifurcation, as in Fig.~5.1d, we again observe complex spatiotemporal patterns but now oscillations do not occur outside the patterned region because the system with no diffusion no longer has a stable limit cycle. With a yet larger value of $\bar{v}_1$ the pattern forms a relatively ordered triangular structure bearing an interseting resemblance to the Sierpinski triangle. Numerical simulations performed over a longer time-scale suggest that this structure persists indefinitely. Fig.~\ref{fig:chaos_pp1} shows a typical profile of the solution at a large time.  However, by increasing the value of $\bar{v}_1$ further, as in Fig.~\ref{fig:SPT_v1}f, patterns are no longer observed.  Here the solution simply decays to the stable homogeneous steady state (the lower equilibrium branch of Fig.~\ref{fig:onepar}a).

\begin{figure}[htbp]
\centering
  \begin{subfigure}[b]{.3\linewidth}
    \centering
    \caption{}
    \includegraphics[width=.99\textwidth]{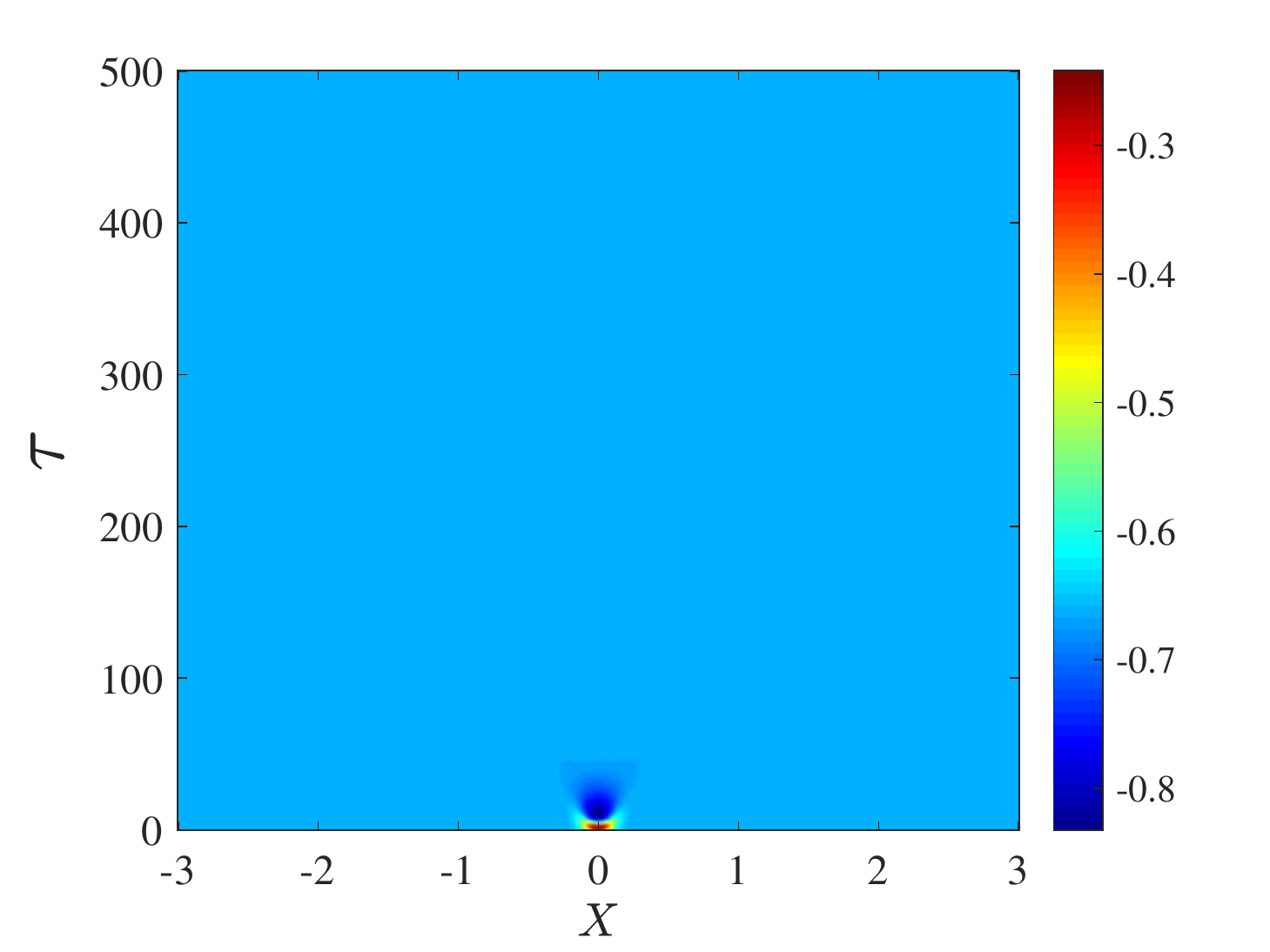}
    \label{fig:sptv303462}
  \end{subfigure}%
  \begin{subfigure}[b]{.3\linewidth}
    \centering
    \caption{}
    \includegraphics[width=.99\textwidth]{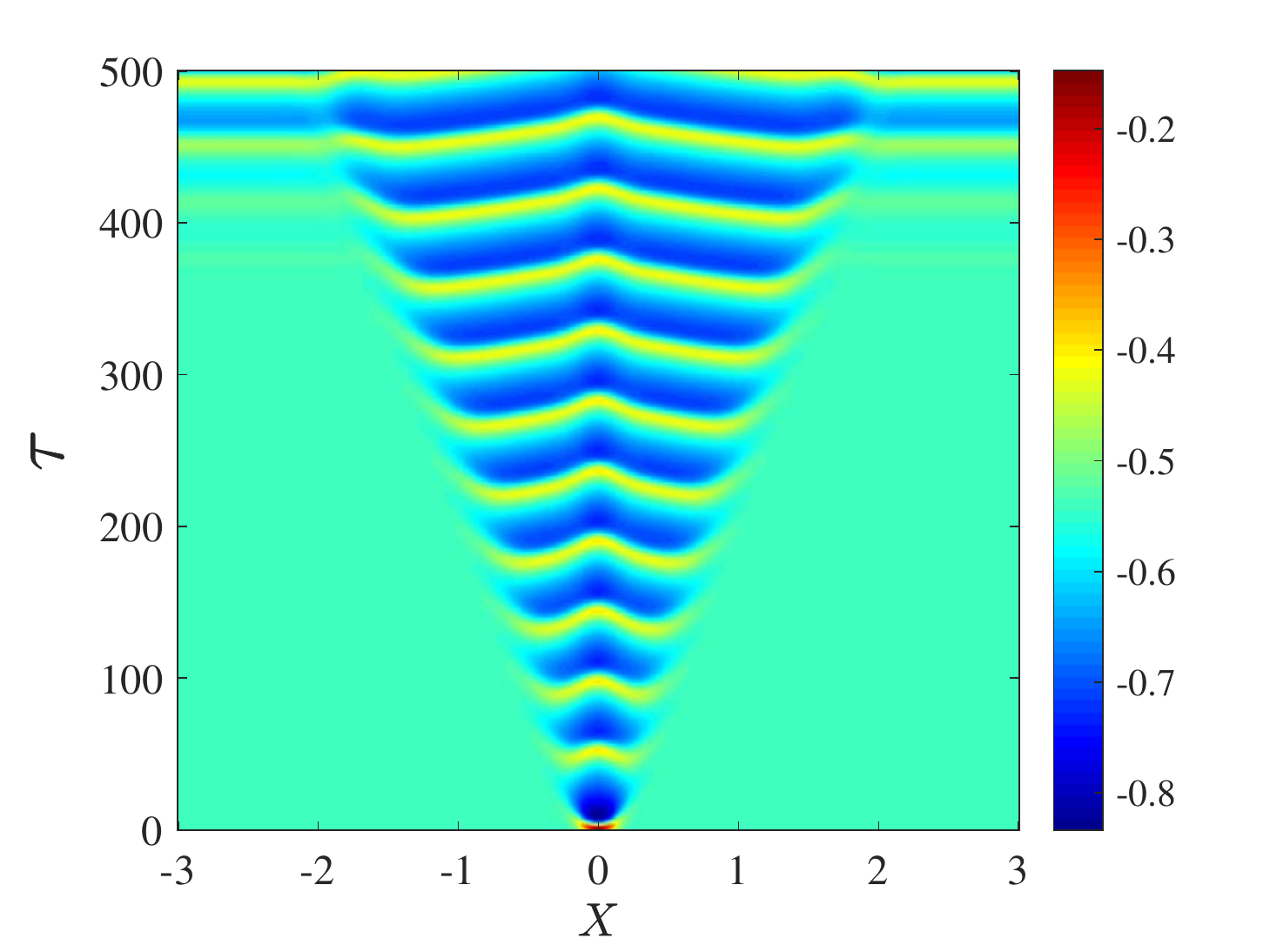}
     \label{fig:sptv303019}
  \end{subfigure}%
  \begin{subfigure}[b]{.3\linewidth}
    \centering
    \caption{}
    \includegraphics[width=.99\textwidth]{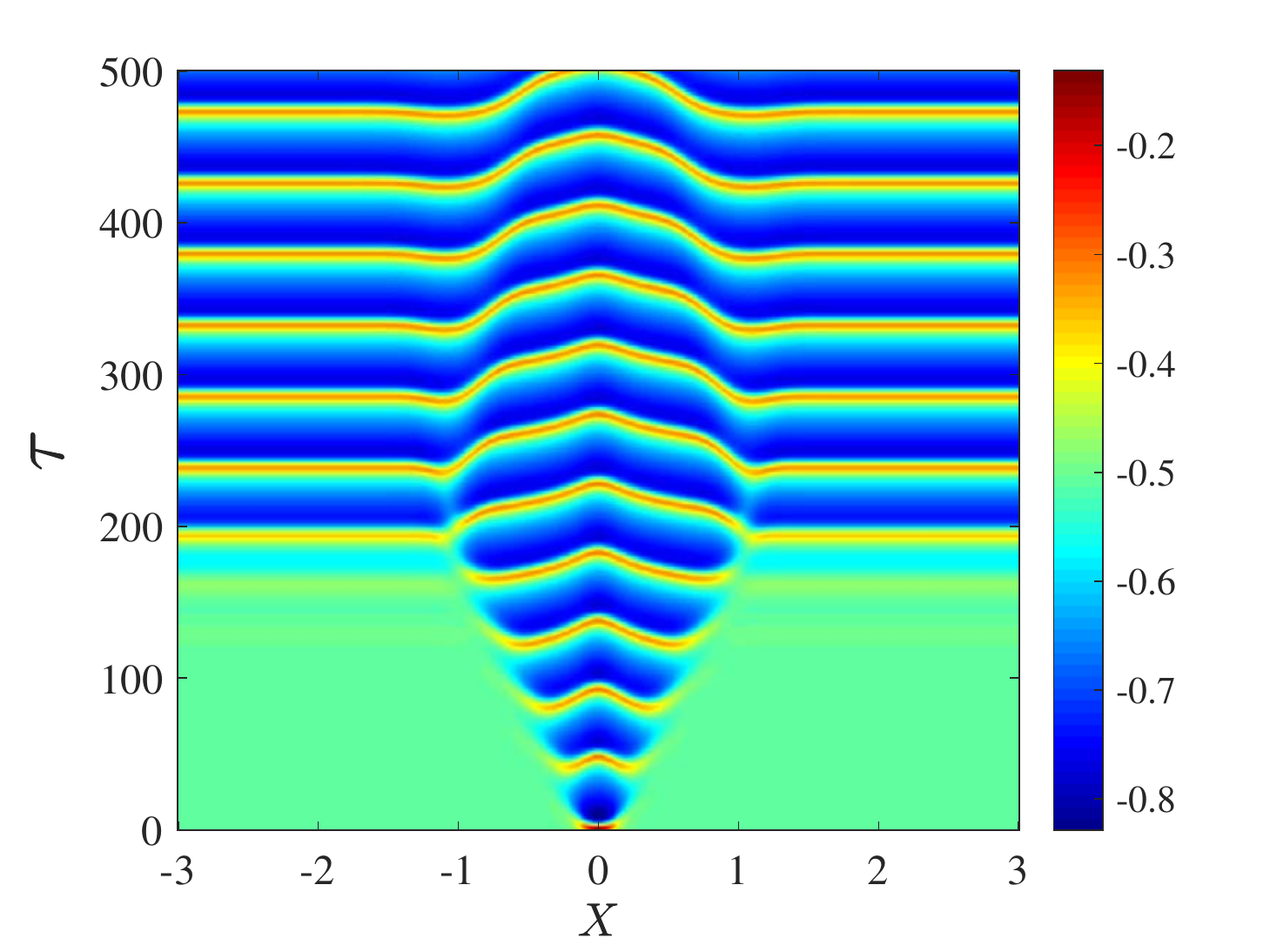}
     \label{fig:sptv302813}
  \end{subfigure}\\%
  \begin{subfigure}[b]{.3\linewidth}
    \centering
    \caption{}
    \includegraphics[width=.99\textwidth]{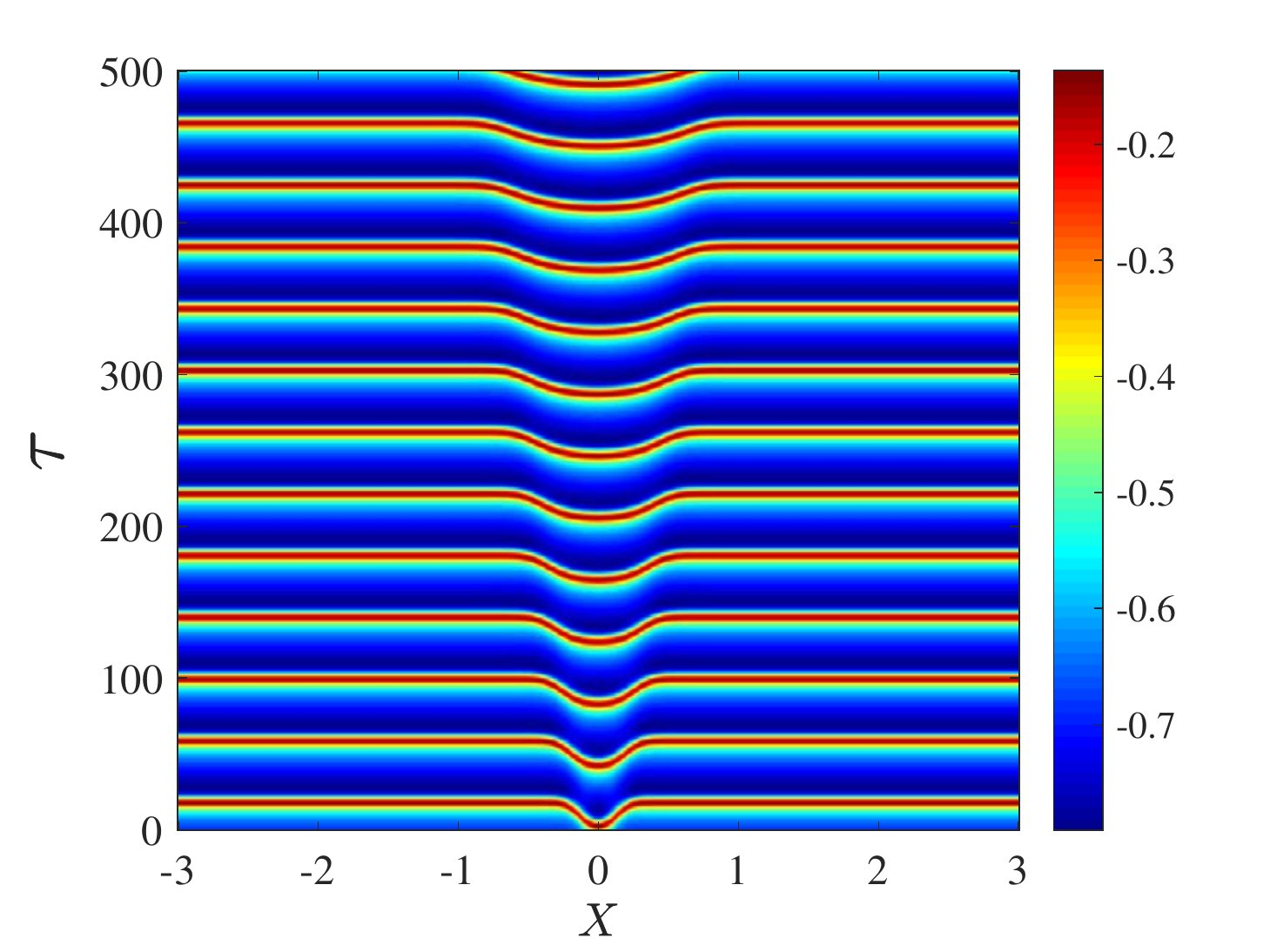}
     \label{fig:sptv302384}
  \end{subfigure}%
  \begin{subfigure}[b]{.3\linewidth}
    \centering
    \caption{}
    \includegraphics[width=.99\textwidth]{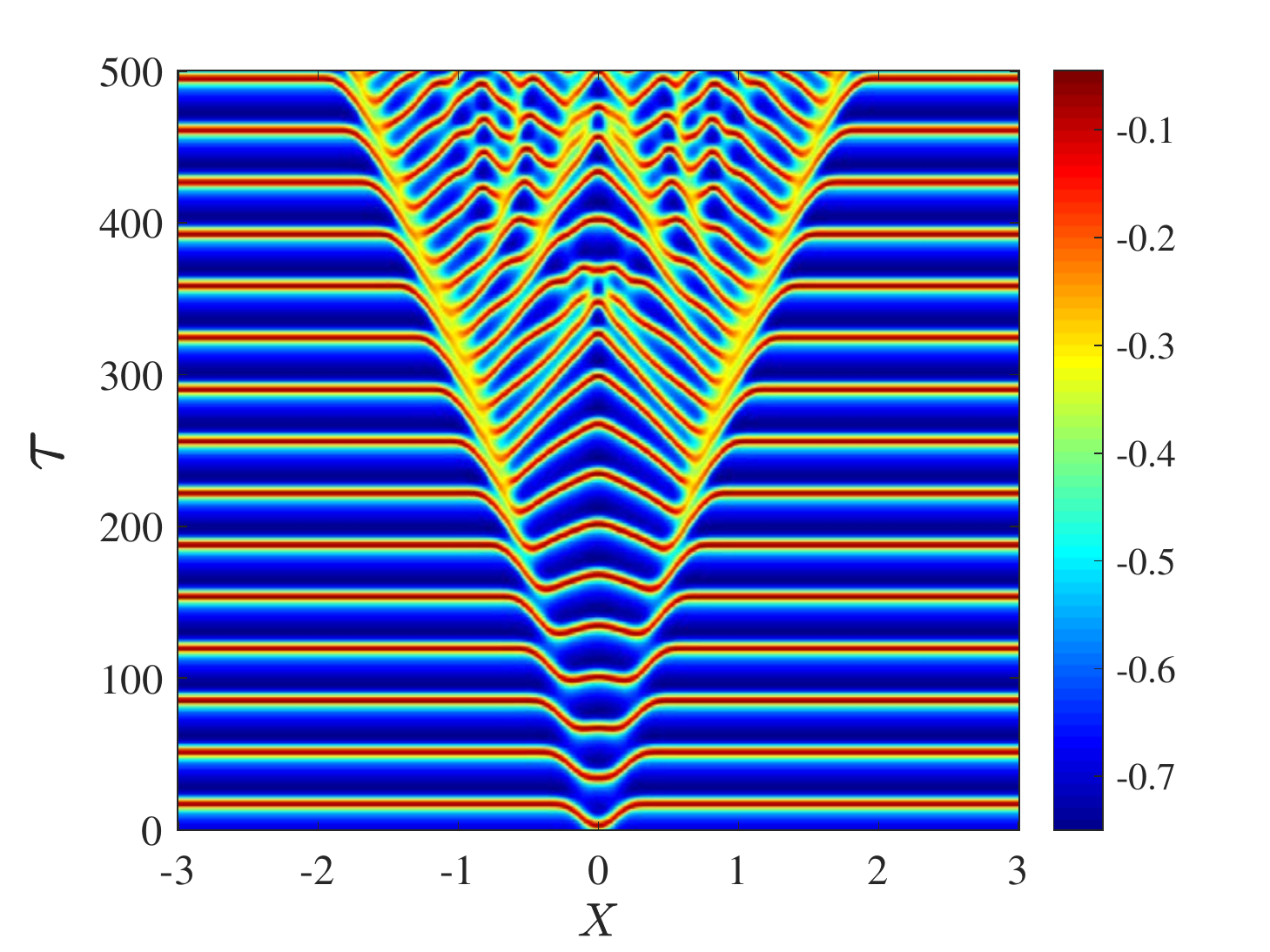}
     \label{fig:sptv301725}
  \end{subfigure}%
  \begin{subfigure}[b]{.3\linewidth}
    \centering
    \caption{}
    \includegraphics[width=0.99\textwidth]{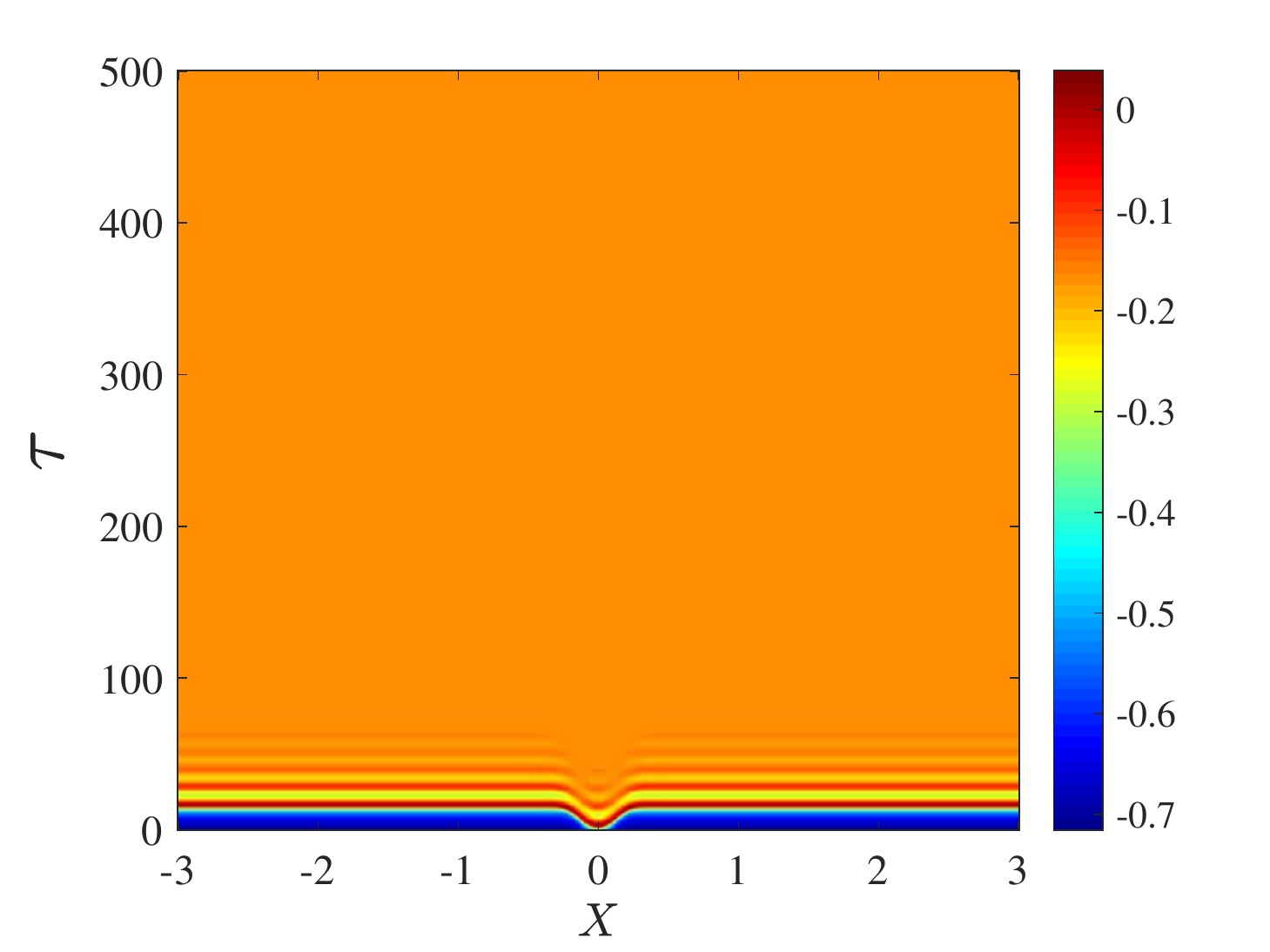}
     \label{fig:sptv300556}
  \end{subfigure}%
 \caption{Space-time plots of the membrane potential $V$ for the values of $\bar{v}_{3}$ marked in Fig.~\ref{fig:v3bif}. Specifically (a) $-0.3462$; (b) $-0.3019$; (c) $-0.2813$; (d) $-0.23842$; (e) $-0.1725$; and (f) $-0.05565$. The initial condition is \eqref{eq:init} with \eqref{eq:gauss}, using the stable and unstable equilibrium branch of Fig.~\ref{fig:v3bif} for the steady state $(V^*,N^*)$, and all other parameters are fixed as in Section~\ref{sec:2}.}
 \label{fig:SPT_v3}
\end{figure}

Now we study the spatiotemporal behaviour of the model by varying $\bar{v}_{3}$ and keeping all other parameters fixed as in Section~\ref{sec:2}. Recall that in this case the system in the absence of diffusion exhibits supercritical and subcritical Hopf bifurcations, see Fig.~\ref{fig:onepar}b. The results of numericals simulations are shown in Fig.~\ref{fig:SPT_v3}. For extremely low values of $\bar{v}_{3}$, the system returns quickly to the homogeneous steady state. Between the Hopf bifurcations, where the system in the absence of diffusion has a stable limit cycle, we observe mostly homogeneous oscillations corresponding to this limit cycle, see Fig.~\ref{fig:SPT_v3}b--e. In panels (b) and (c) away from $X=0$ where the perturbation is applied, it takes some time for the solution to settle to oscillatory behaviour because the initial condition is set very near the value unstable steady state. In panels (d) and (e) oscillations develop across the domain relatively quickly. In panel (e), which is just before the subcritical Hopf bifurcation, the initial stimulus creates a pulse of propagating action potentials. For values of $\bar{v}_3$ beyond the subcritical Hopf bifurcation and subsequent saddle-node bifurcation SNC (see Fig.~\ref{fig:onepar}b), periodic oscillations can be observed for a short time across the entire domain, then stabilise to the homogeneous steady state, as in Fig.~\ref{fig:SPT_v3}f.

\begin{figure}[htbp]
\centering
  \begin{subfigure}[b]{.3\linewidth}
    \centering
    \caption{}
    \includegraphics[width=0.99\textwidth]{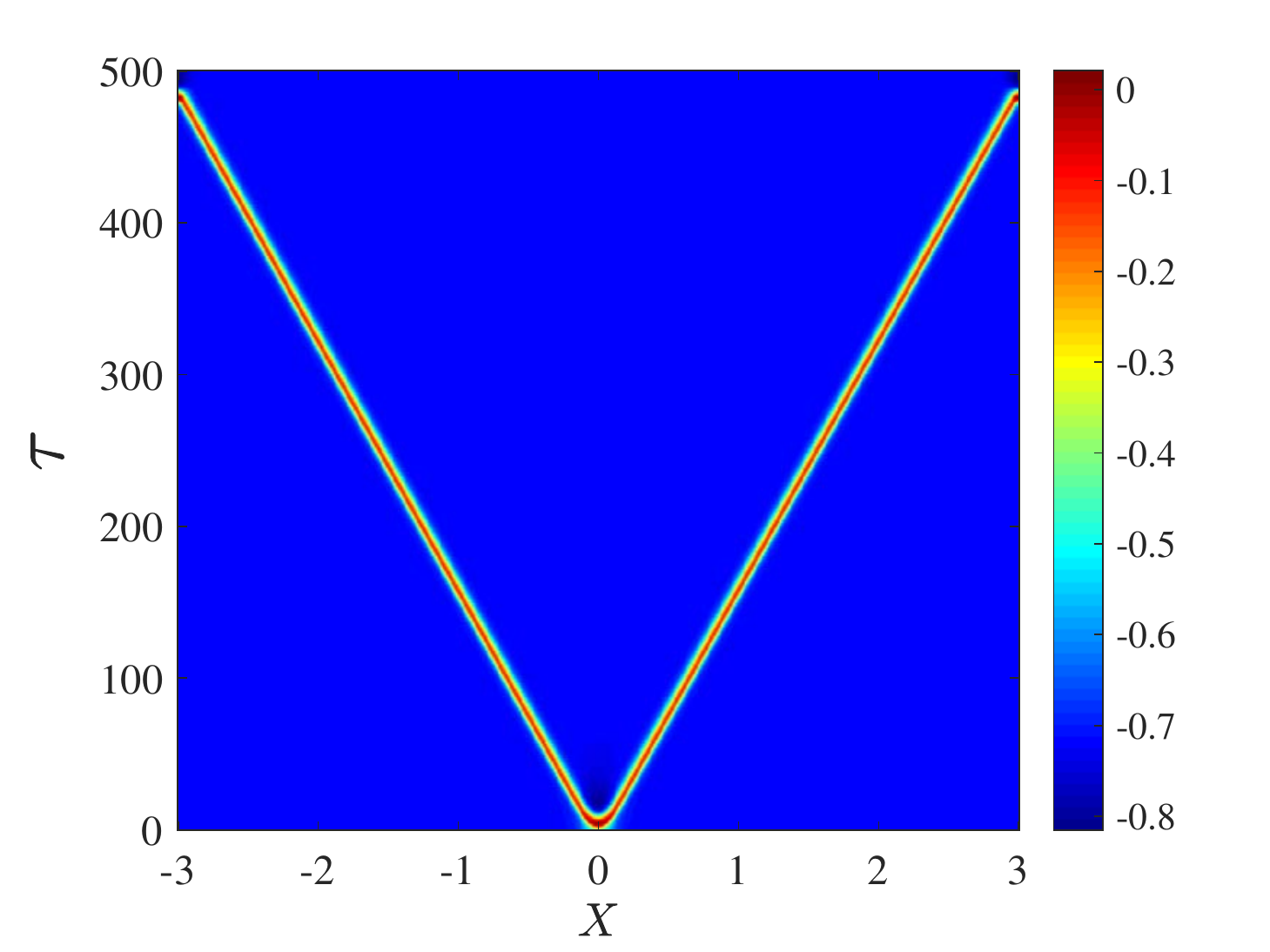}
    \label{fig:sptpsi01}
  \end{subfigure}%
  \begin{subfigure}[b]{.3\linewidth}
    \centering
    \caption{}
    \includegraphics[width=0.99\textwidth]{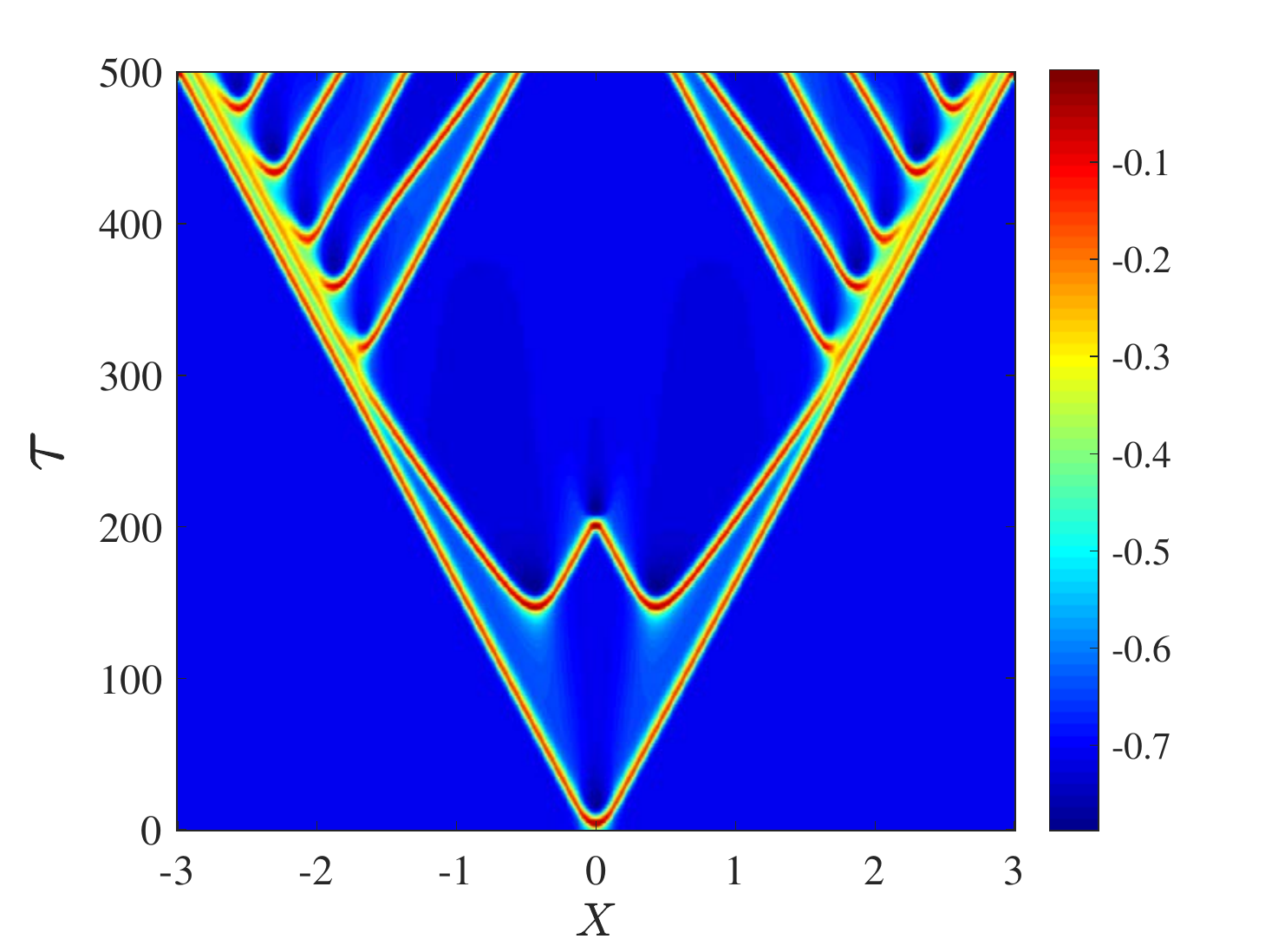}
     \label{fig:sptpsi012}
  \end{subfigure}%
  \begin{subfigure}[b]{.3\linewidth}
    \centering
    \caption{}
    \includegraphics[width=0.99\textwidth]{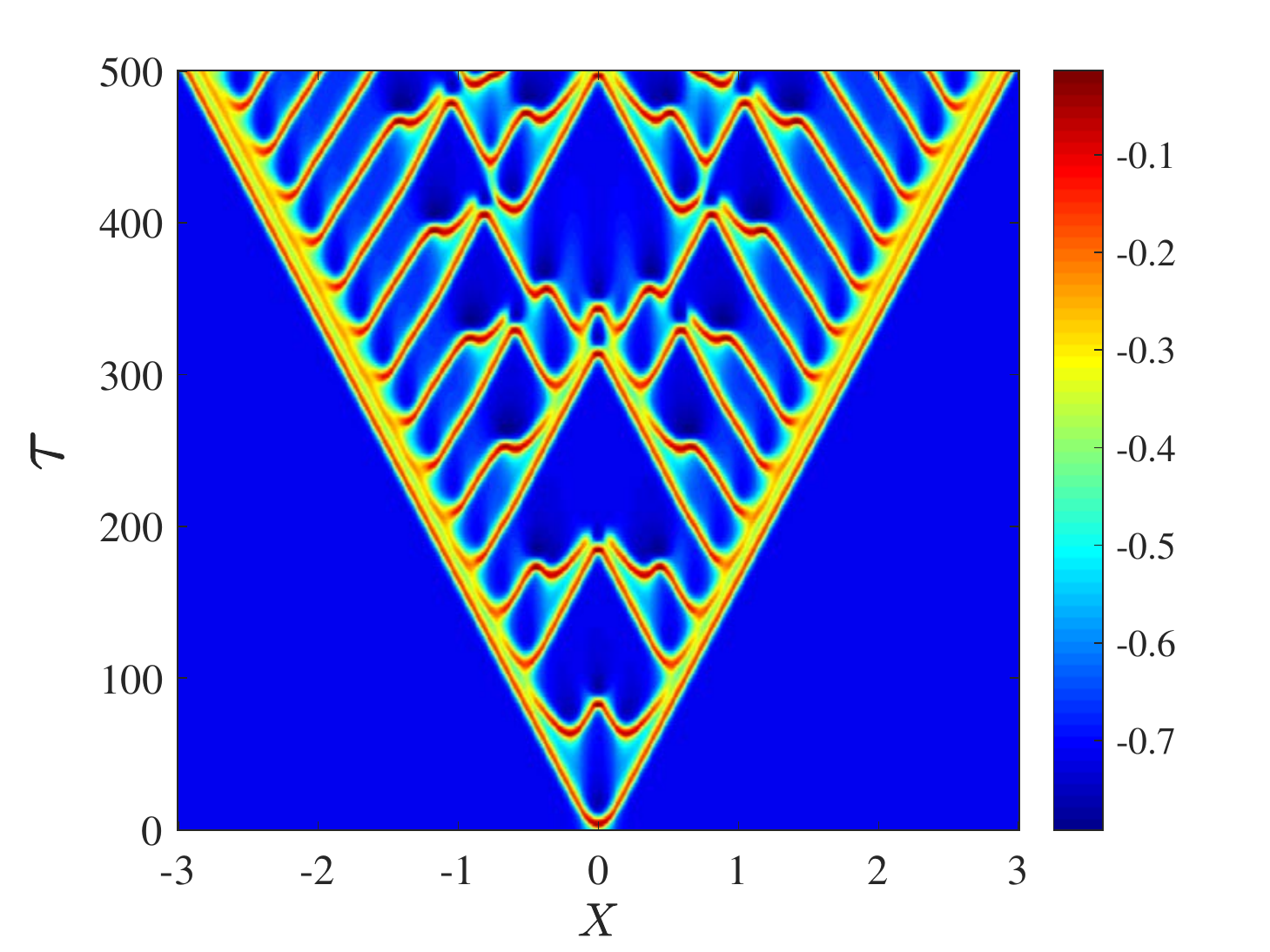}
     \label{fig:sptpsi013}
  \end{subfigure}\\%
  \begin{subfigure}[b]{.3\linewidth}
    \centering
    \caption{}
    \includegraphics[width=0.99\textwidth]{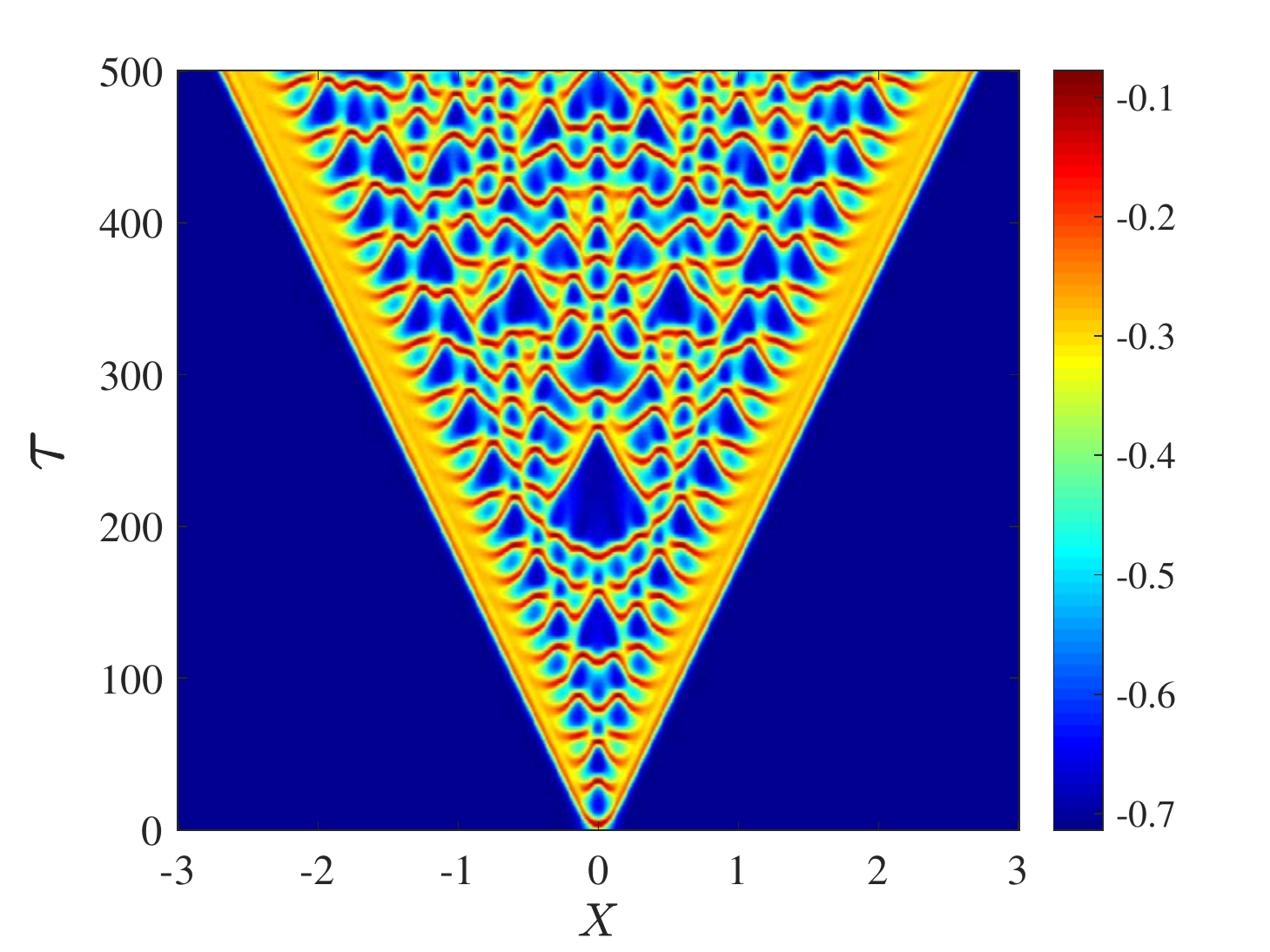}
     \label{fig:sptpsi02}
  \end{subfigure}%
  \begin{subfigure}[b]{.3\linewidth}
    \centering
    \caption{}
    \includegraphics[width=0.99\textwidth]{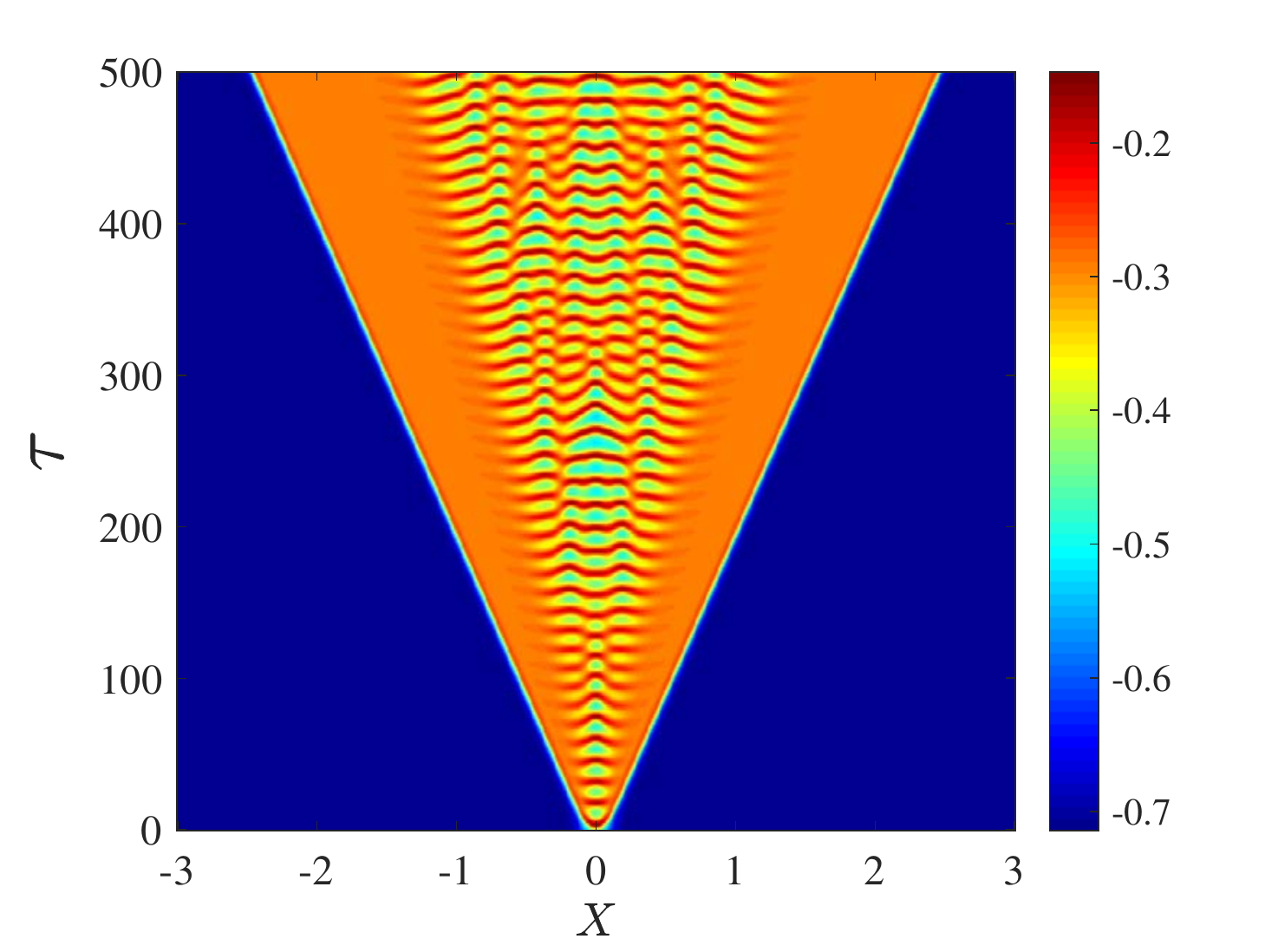}
     \label{fig:sptpsi03}
  \end{subfigure}%
  \begin{subfigure}[b]{.3\linewidth}
    \centering
    \caption{}
    \includegraphics[width=0.99\textwidth]{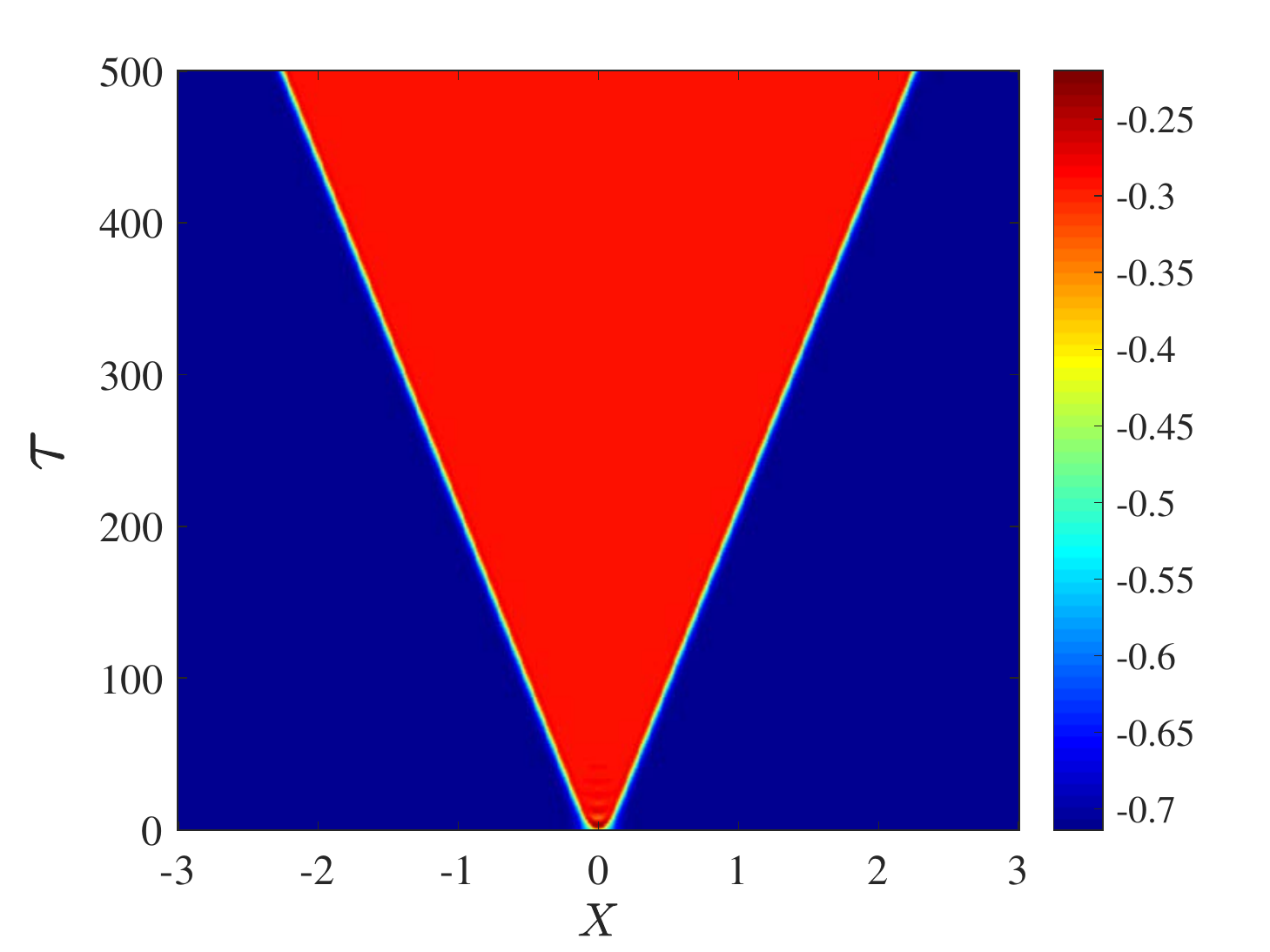}
     \label{fig:sptpsi05}
  \end{subfigure}%
 \caption{Space-time plots of the membrane potential $V$ for values of $\psi$ as marked in Fig.~\ref{fig:psi}. Specifically (a) $0.1$; (b) $0.12$; (c) $0.13$; (d) $0.2$; (e) $0.3$; and (f) $0.5$. The initial condition is \eqref{eq:init} with \eqref{eq:gauss}, using the upper equilibrium branch of Fig.~\ref{fig:v1bif} for the steady state $(V^*,N^*)$, and all other parameters are fixed as in Section~\ref{sec:2}. The solution transitions from propagating pulses travelling in opposite direction to complex spatiotemporal patterns to fronts travelling in opposite direction.}
 \label{fig:SPT_psi}
\end{figure}

Finally Fig.~\ref{fig:SPT_psi} the spatiotemporal patterns for the various values of $\psi$ marked in Fig.~\ref{fig:onepar}c. For extremely low values of $\psi$, the initial perturbation creates a pulse at the centre of the domain and as time progresses the pulse splits into two travelling pulses propagating in opposite directions at the same speed (Fig.~\ref{fig:SPT_psi}a). A slight increase in the value of $\psi$ leads to a destabilisation of the pulses that results in an initiation of secondary pulses travelling in the opposite direction to the primary pulses (Fig.~\ref{fig:SPT_psi}b). By increasing the value of $\psi$ further we are able to see within the $\tau = 500$ time-frame that the secondary pulses collide with one another and eventually irregular oscillations disseminate across the spatial domain (Fig.~\ref{fig:SPT_psi}c--d). Interestingly, as $\psi$ is varied past the homoclinic bifurcation, the unstable pulses transition to travelling fronts connecting a stable steady state to an unstable state with irregular oscillations at the back of the fronts (Fig.~\ref{fig:SPT_psi}e). As the value of $\psi$ is increased further, the upper equilibrium branch gains stability at the Hopf bifurcation so beyond this bifurcation the system has two stable steady states. In this case the fronts connect one stable steady state to the other (Fig.~\ref{fig:SPT_psi}f).

\begin{figure}[htbp]
\centering
  \begin{subfigure}[b]{.3\linewidth}
    \centering
    \includegraphics[width=.99\textwidth]{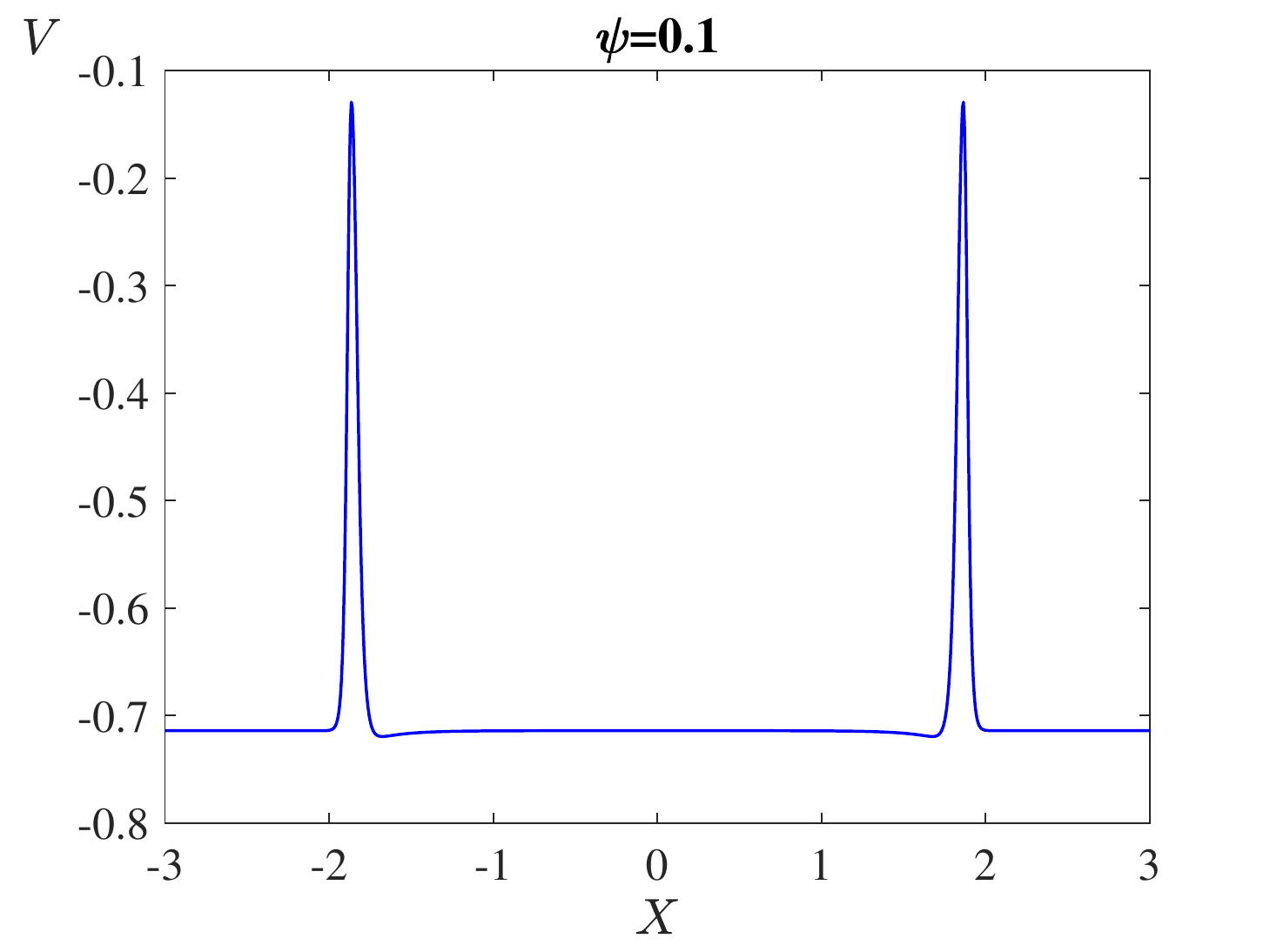}
   \label{fig:1a}
  \end{subfigure}%
  \begin{subfigure}[b]{.3\linewidth}
    \centering
    \includegraphics[width=.99\textwidth]{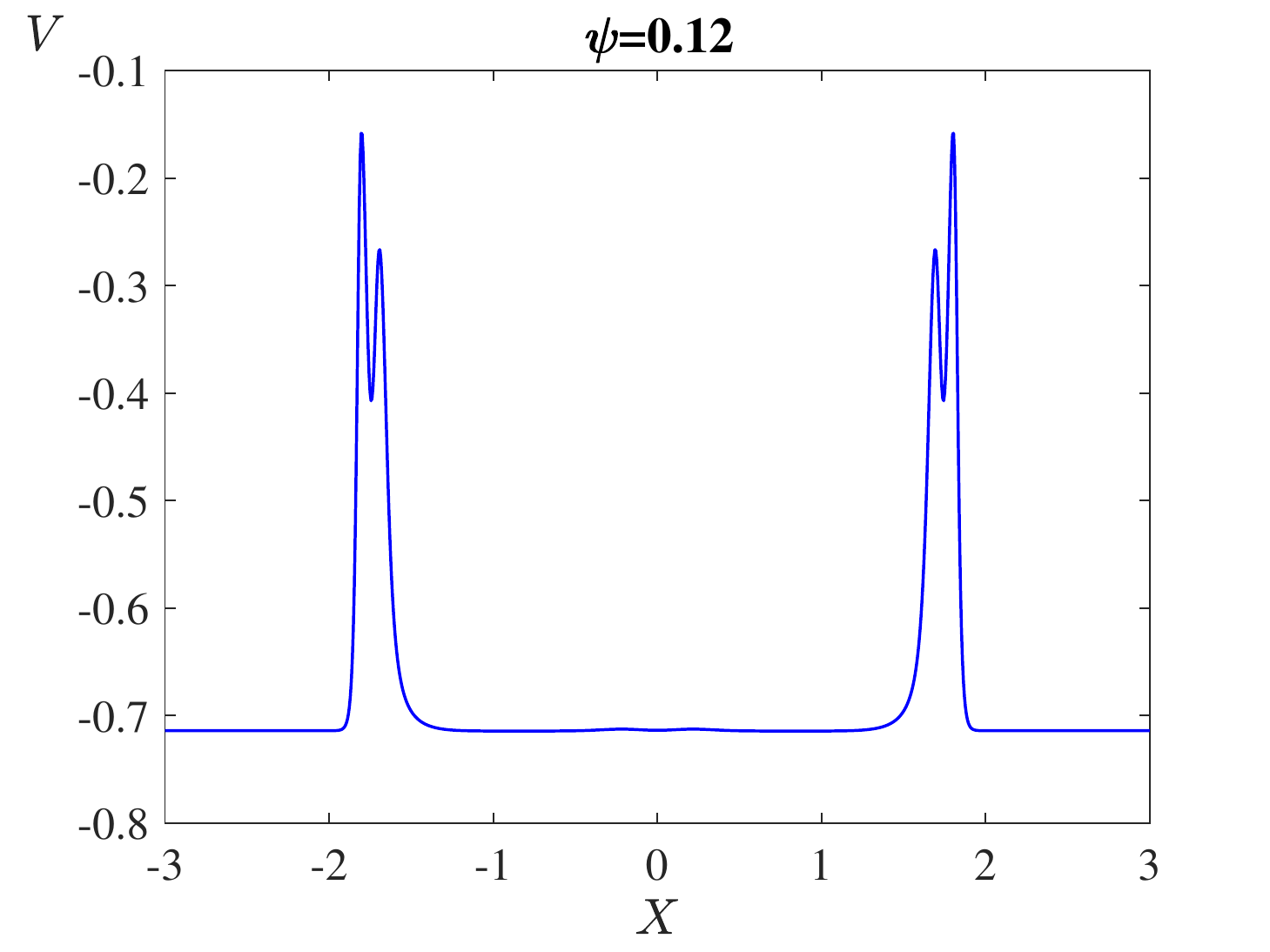}
    \label{fig:1b}
  \end{subfigure}%
  \begin{subfigure}[b]{.3\linewidth}
    \centering
    \includegraphics[width=.99\textwidth]{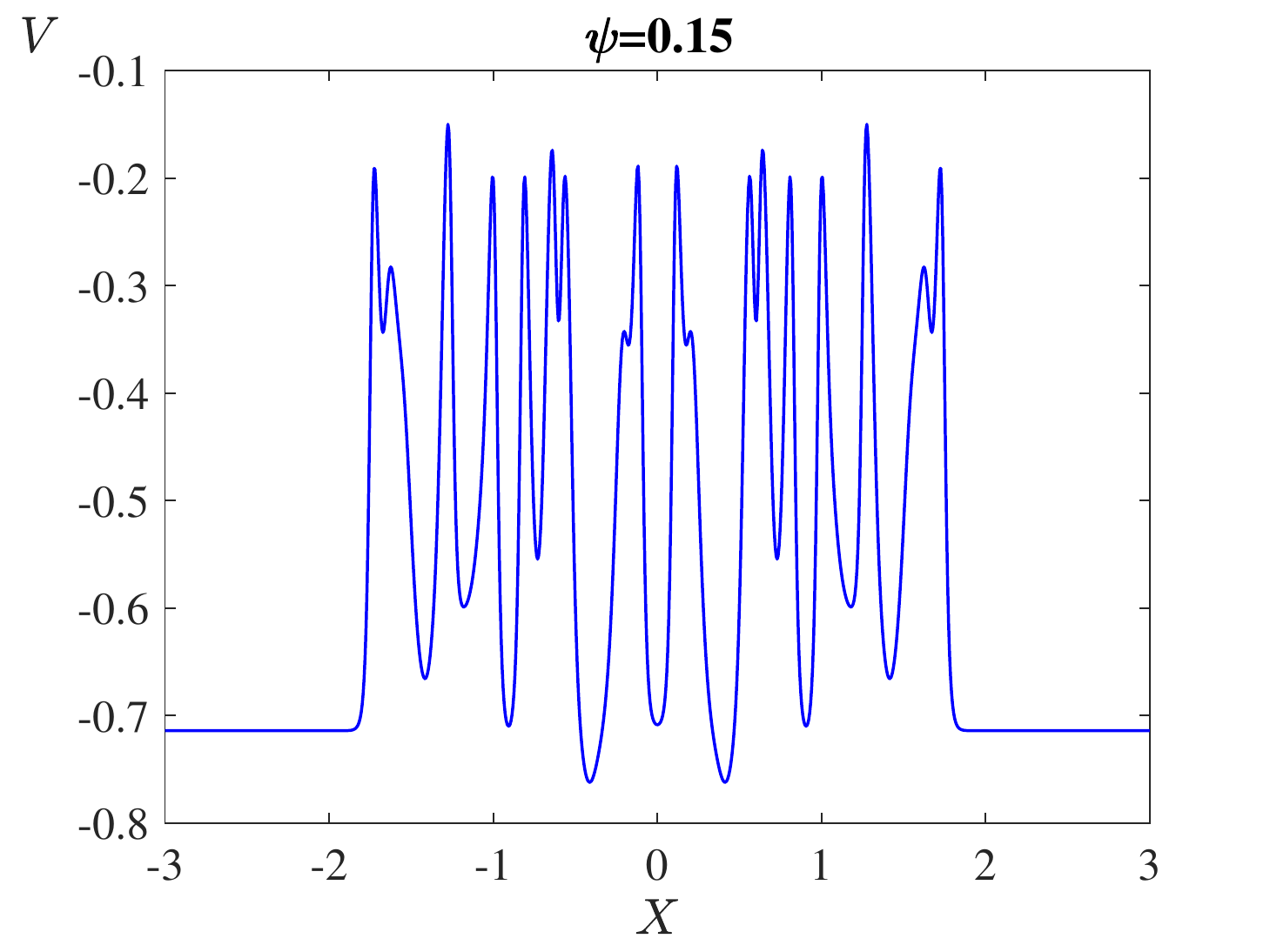}
    \label{fig:1c}
  \end{subfigure}\\%
  \begin{subfigure}[b]{.3\linewidth}
    \centering
    \includegraphics[width=.99\textwidth]{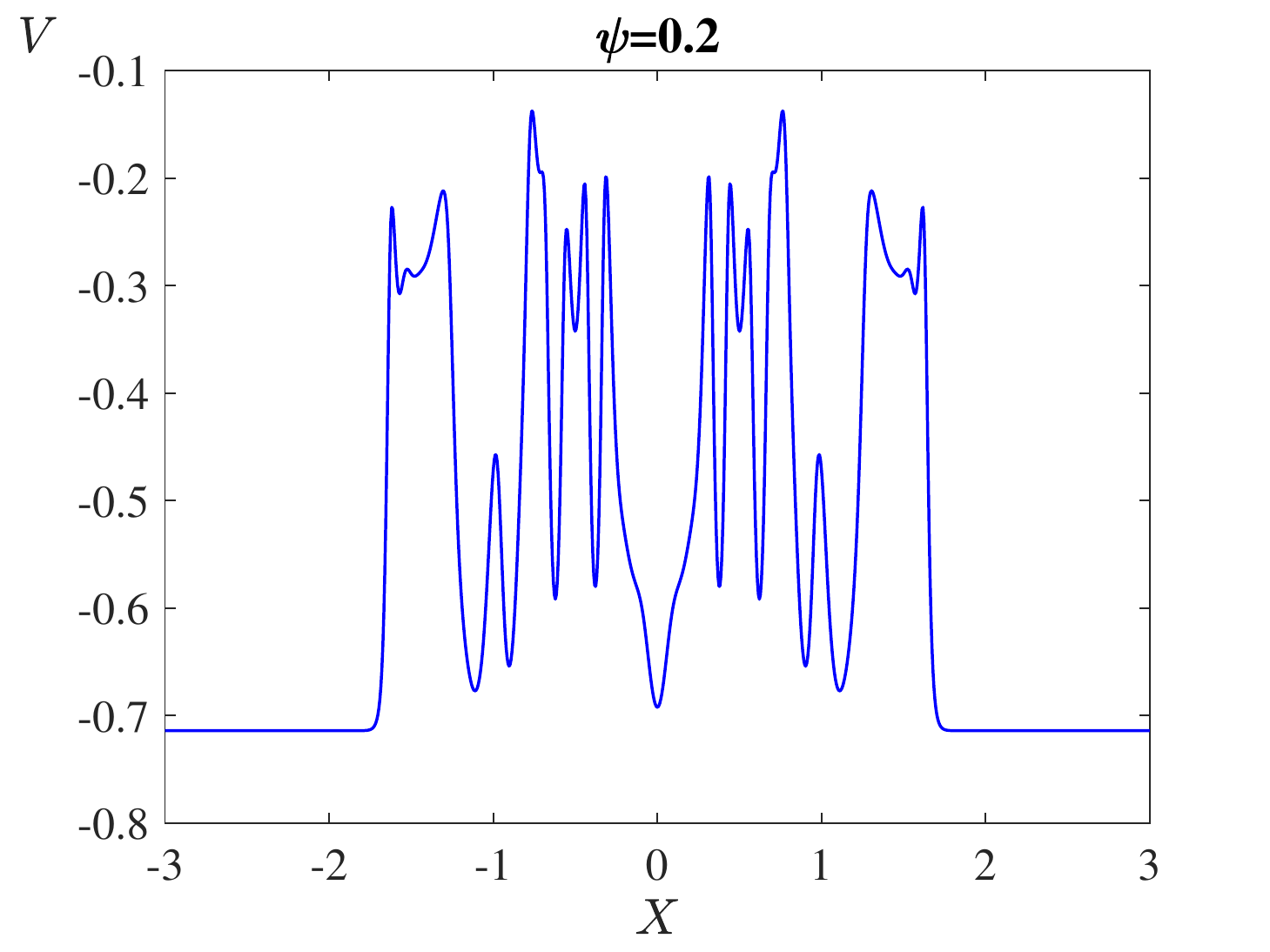}
  \end{subfigure}%
  \begin{subfigure}[b]{.3\linewidth}
    \centering
    \includegraphics[width=.99\textwidth]{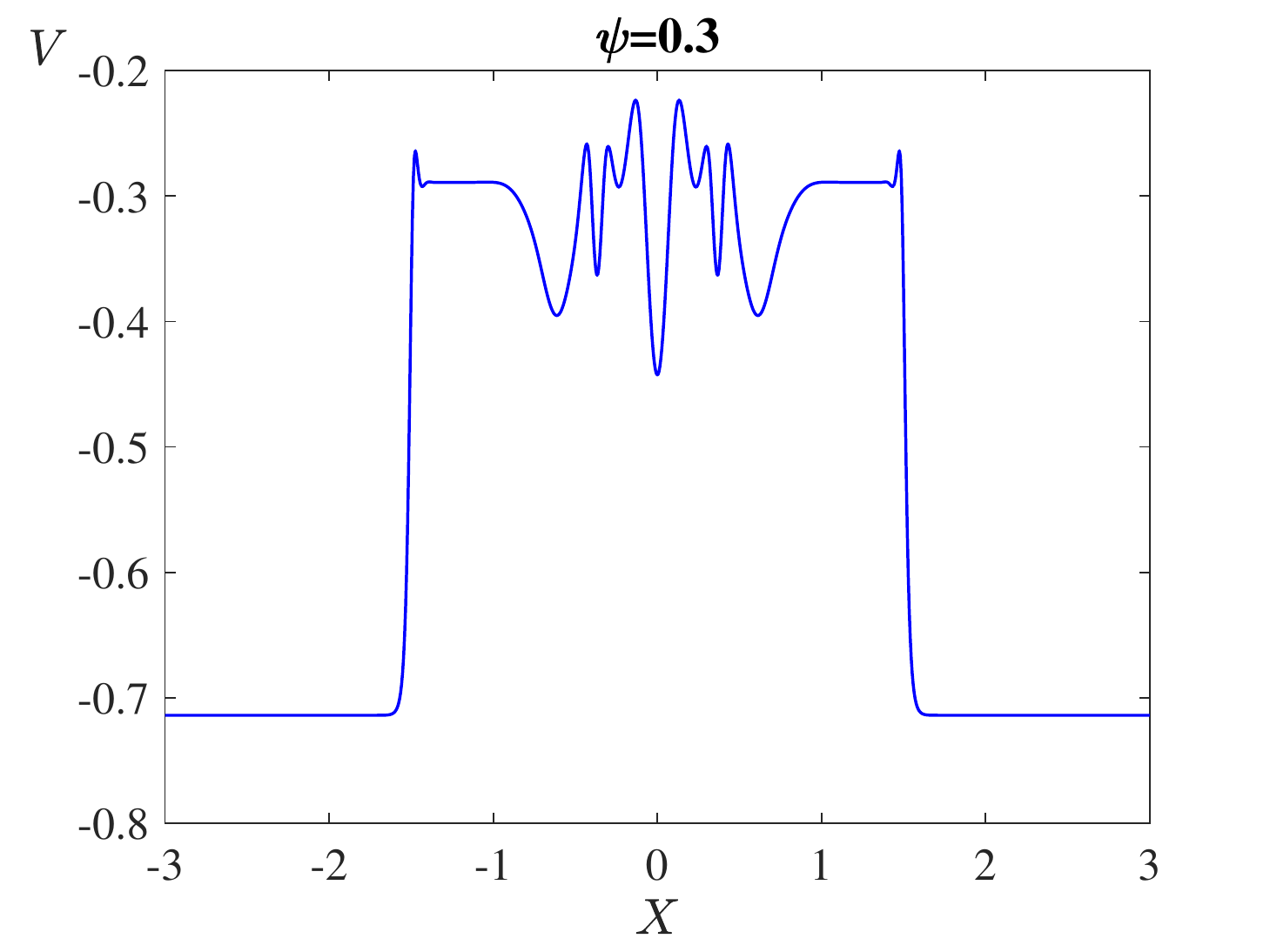}
  \end{subfigure}%
  \begin{subfigure}[b]{.3\linewidth}
    \centering
    \includegraphics[width=.99\textwidth]{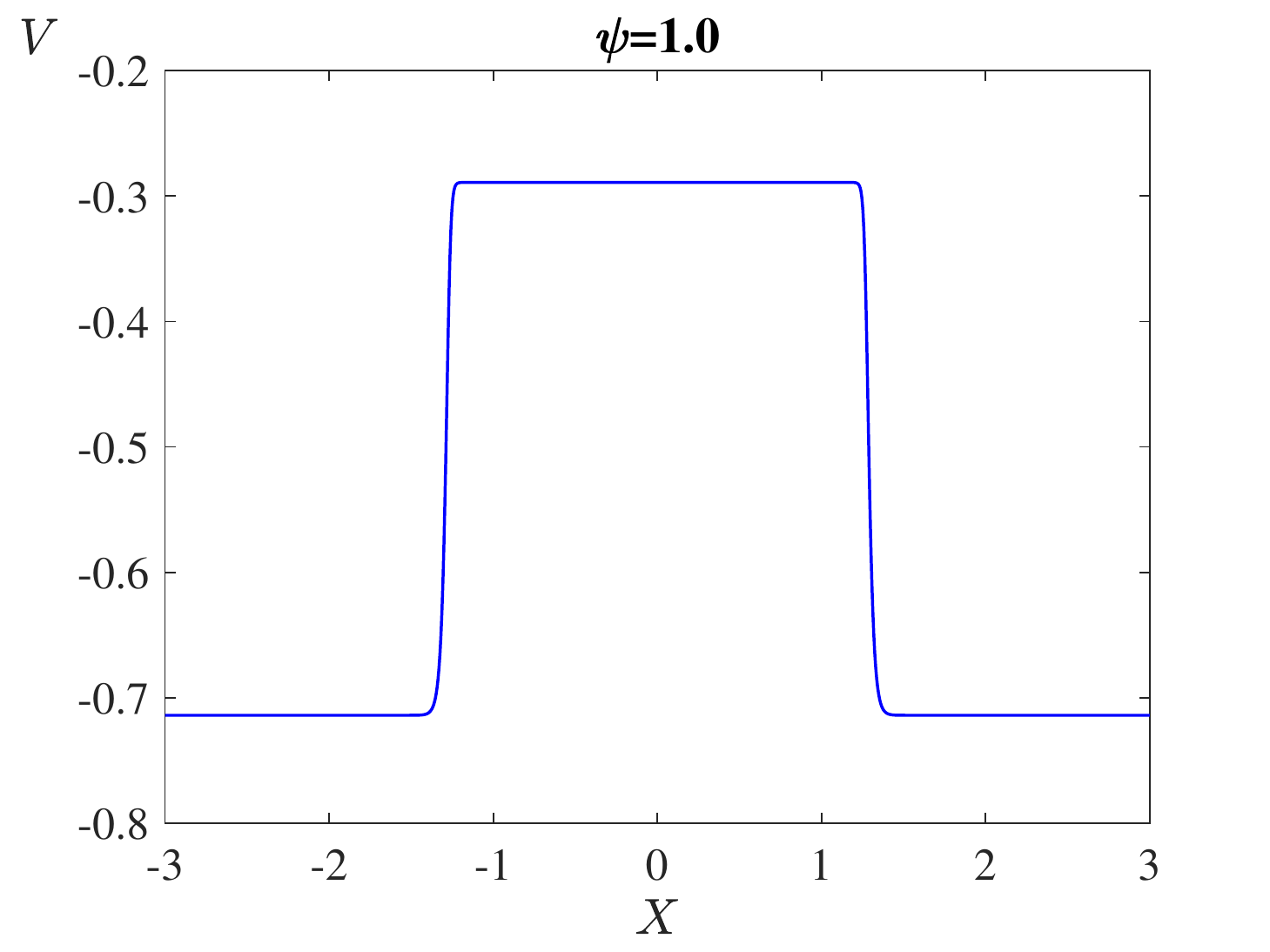}
  \end{subfigure}%
 \caption{Solution profiles at time $\tau=300$ showing the transitions from travelling pulses to spatiotemporal chaos and to fronts.}
 \label{fig:transition}
\end{figure}

Fig.~\ref{fig:transition} shows the solution at $\tau = 300$ for the six values of $\psi$ used in Fig.~\ref{fig:SPT_psi}. This shows how increasing the value of $\psi$ causes the two travelling pulses to transition into two travelling fronts via an intermediate phase of spatiotemporal chaos.

\subsection{Numerical simulations with alternate initial conditions}\label{sec:init_SP}
\begin{figure}[htbp]
\centering
  \begin{subfigure}[b]{.3\linewidth}
    \centering
    \caption{}
    \includegraphics[width=.99\textwidth]{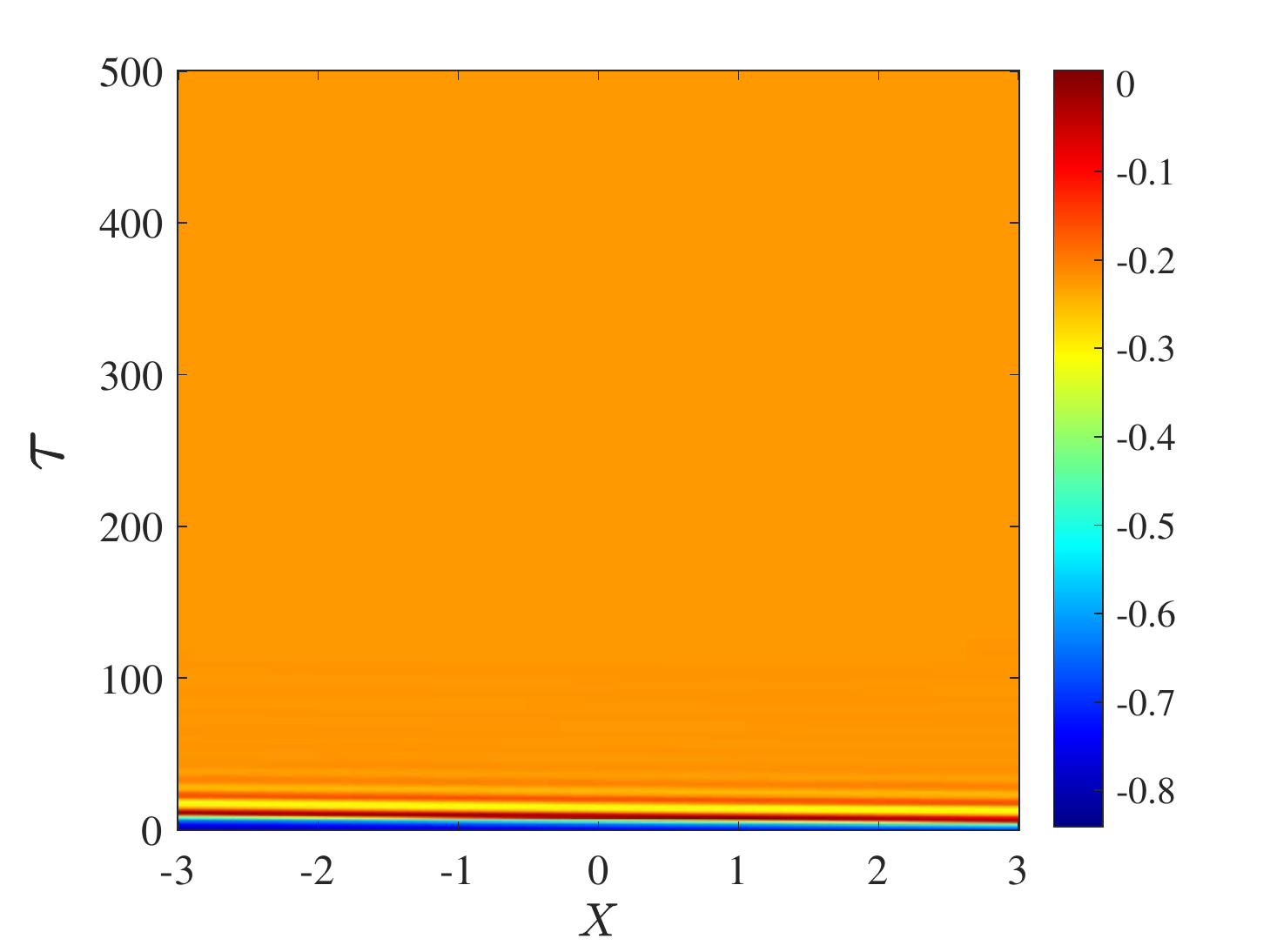}
    \label{fig:spt035v1}
  \end{subfigure}%
  \begin{subfigure}[b]{.3\linewidth}
    \centering
    \caption{}
    \includegraphics[width=.99\textwidth]{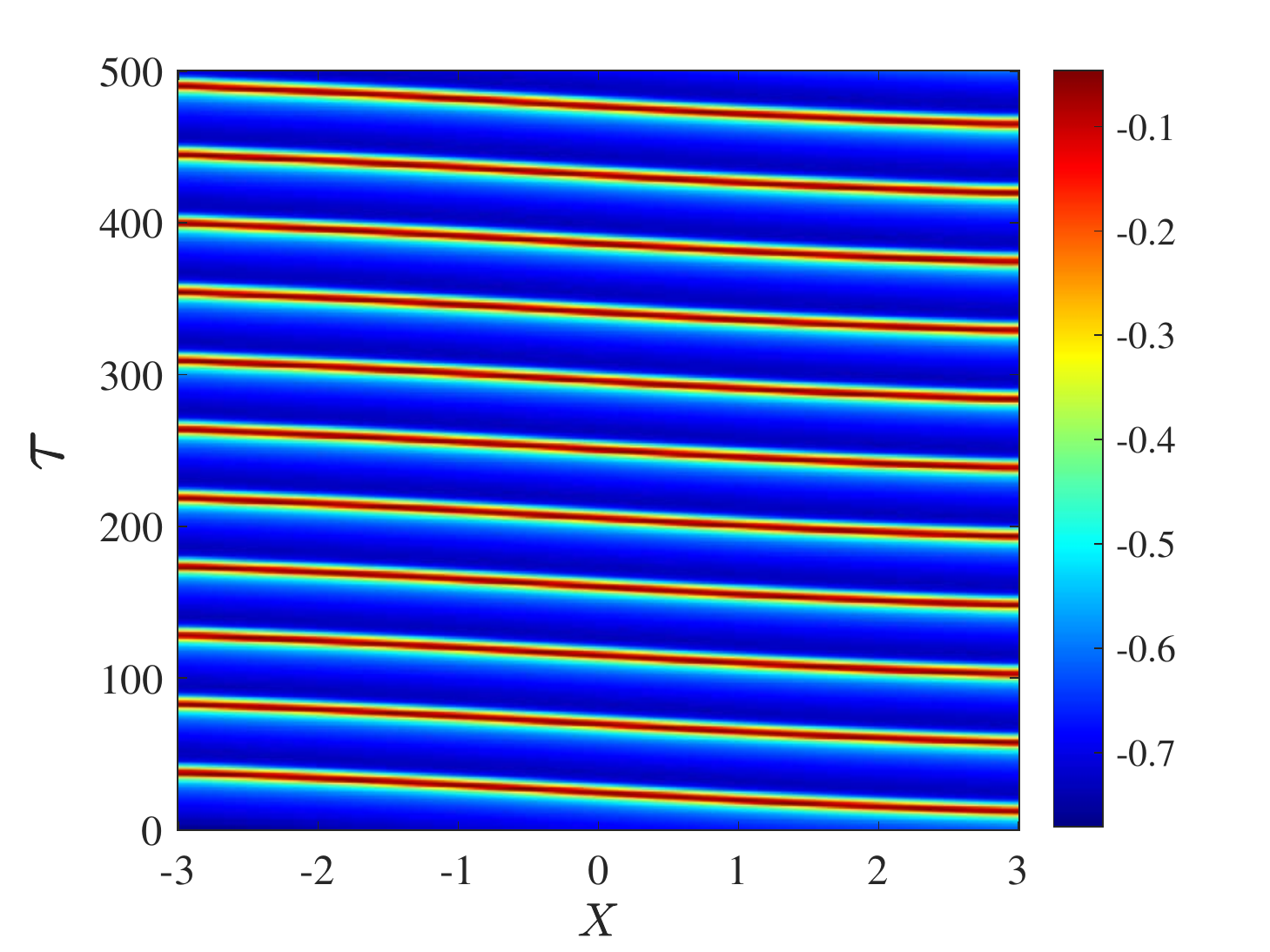}
     \label{fig:spt0265}
  \end{subfigure}%
  \begin{subfigure}[b]{.3\linewidth}
    \centering
    \caption{}
    \includegraphics[width=.99\textwidth]{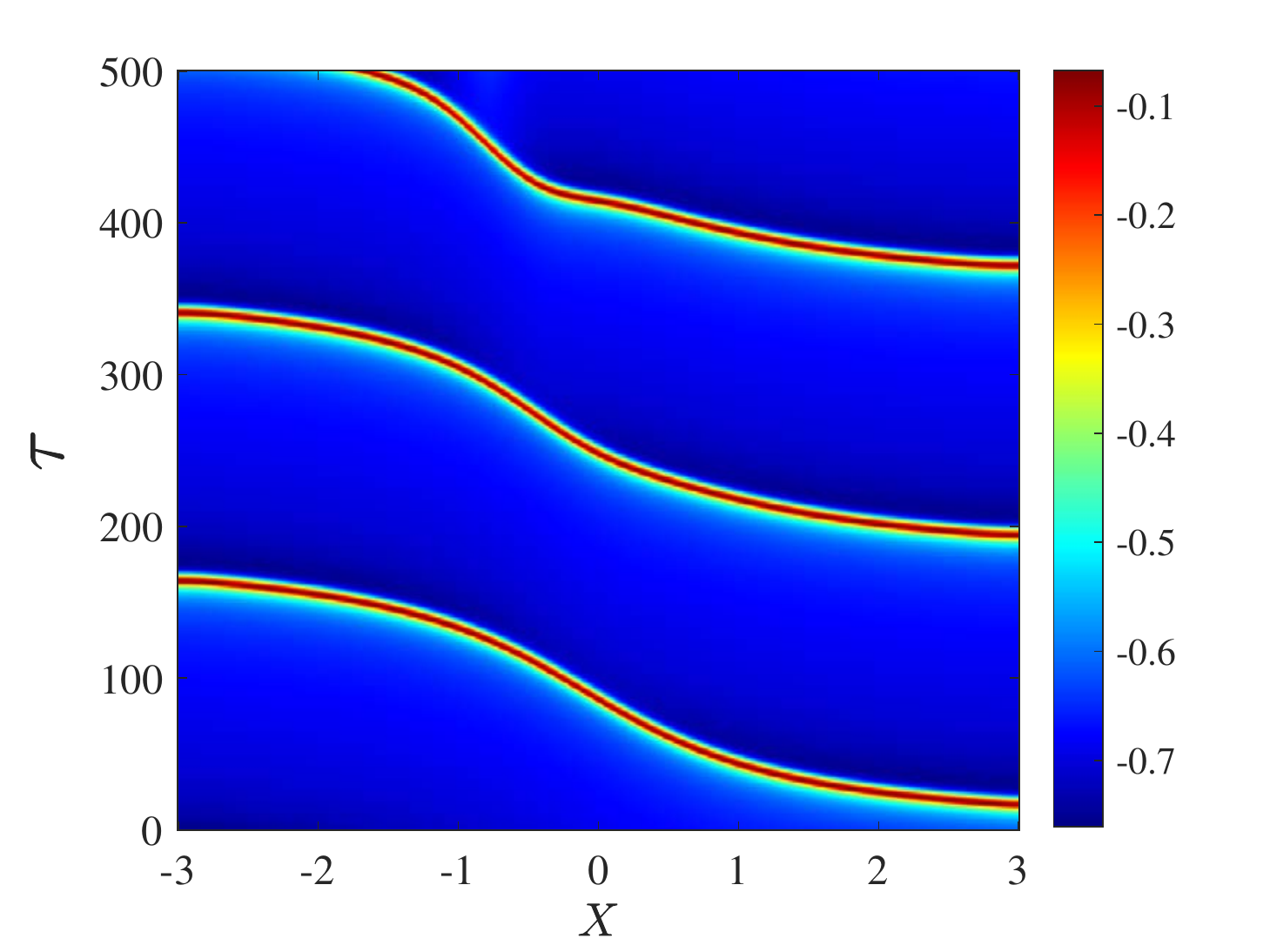}
     \label{fig:spt0255}
  \end{subfigure}\\%
  \begin{subfigure}[b]{.3\linewidth}
    \centering
    \caption{}
    \includegraphics[width=.99\textwidth]{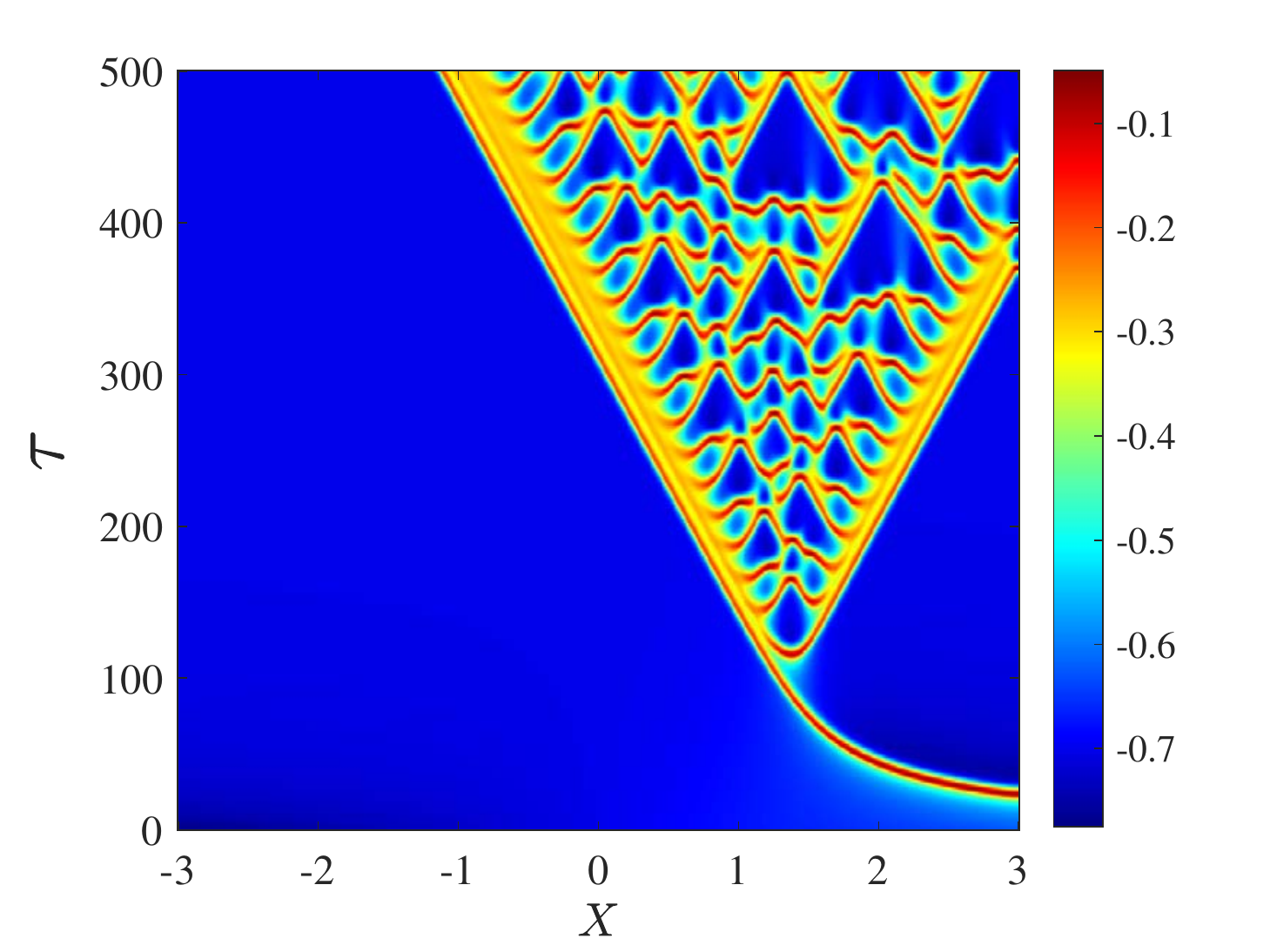}
     \label{fig:sptv0248}
  \end{subfigure}%
  \begin{subfigure}[b]{.3\linewidth}
    \centering
    \caption{}
    \includegraphics[width=.99\textwidth]{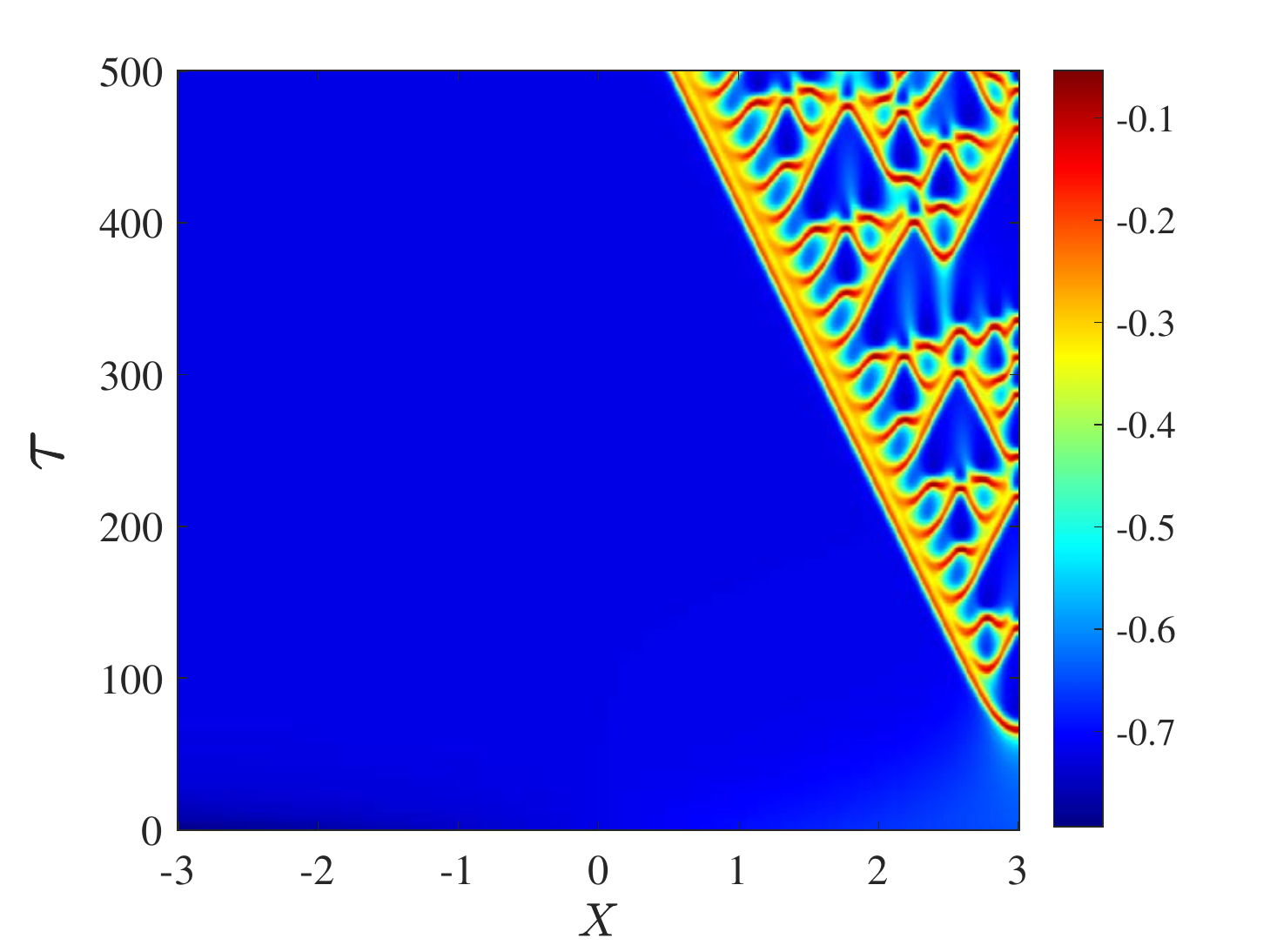}
     \label{fig:spt0246}
  \end{subfigure}%
  \begin{subfigure}[b]{.3\linewidth}
    \centering
    \caption{}
    \includegraphics[width=.99\textwidth]{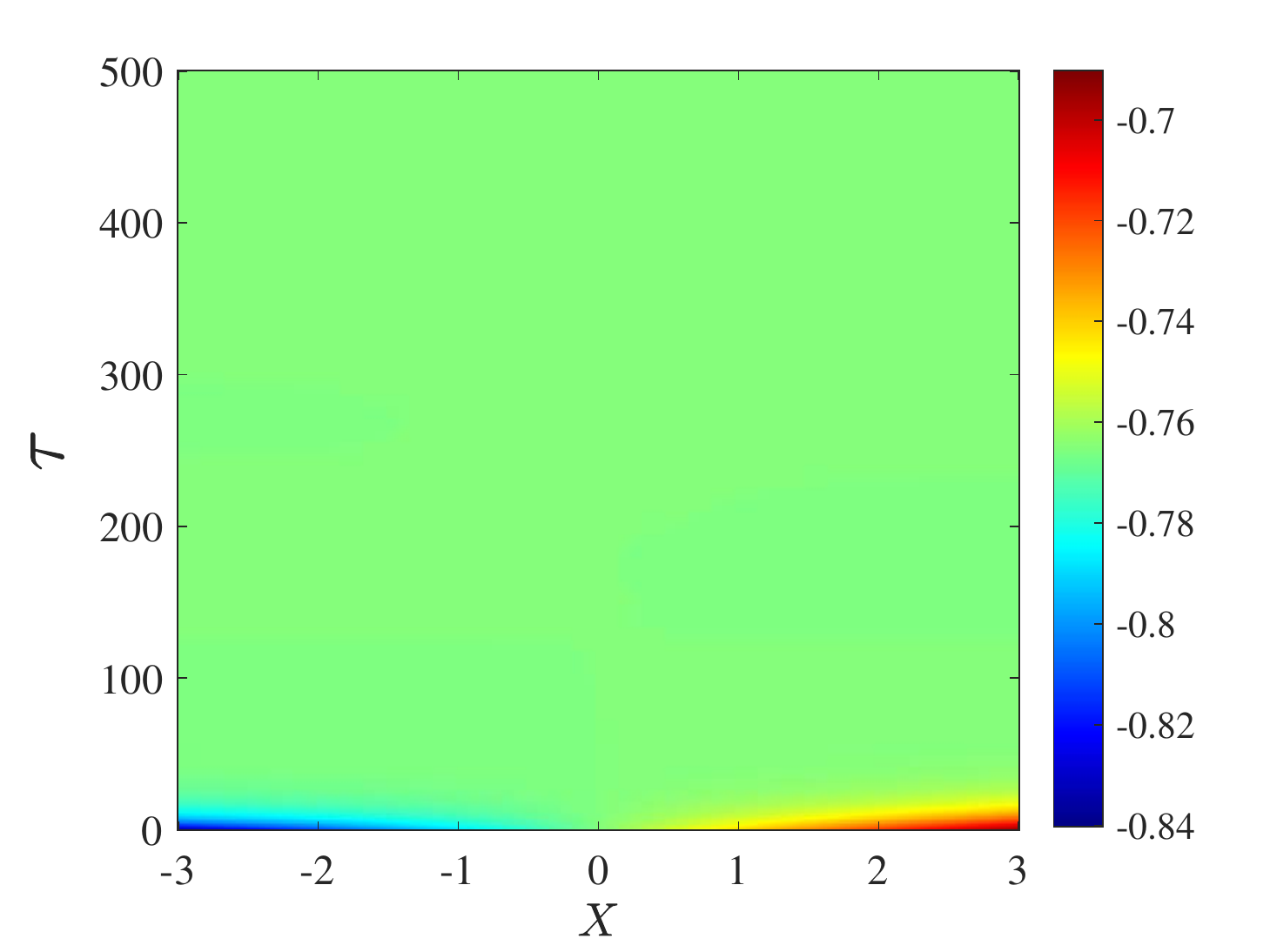}
     \label{fig:spt0245}
  \end{subfigure}%
 \caption{Space-time plots using the same parameter values as Fig.~\ref{fig:SPT_v1} but now with the perturbation function \eqref{eq:init2} in the initial condition \eqref{eq:init}.}
 \label{fig:SPT_v1_X}
\end{figure}

\begin{figure}[htb]
\centering
  \begin{subfigure}[b]{.35\linewidth}
    \centering
    \caption{}
    \includegraphics[width=.99\textwidth]{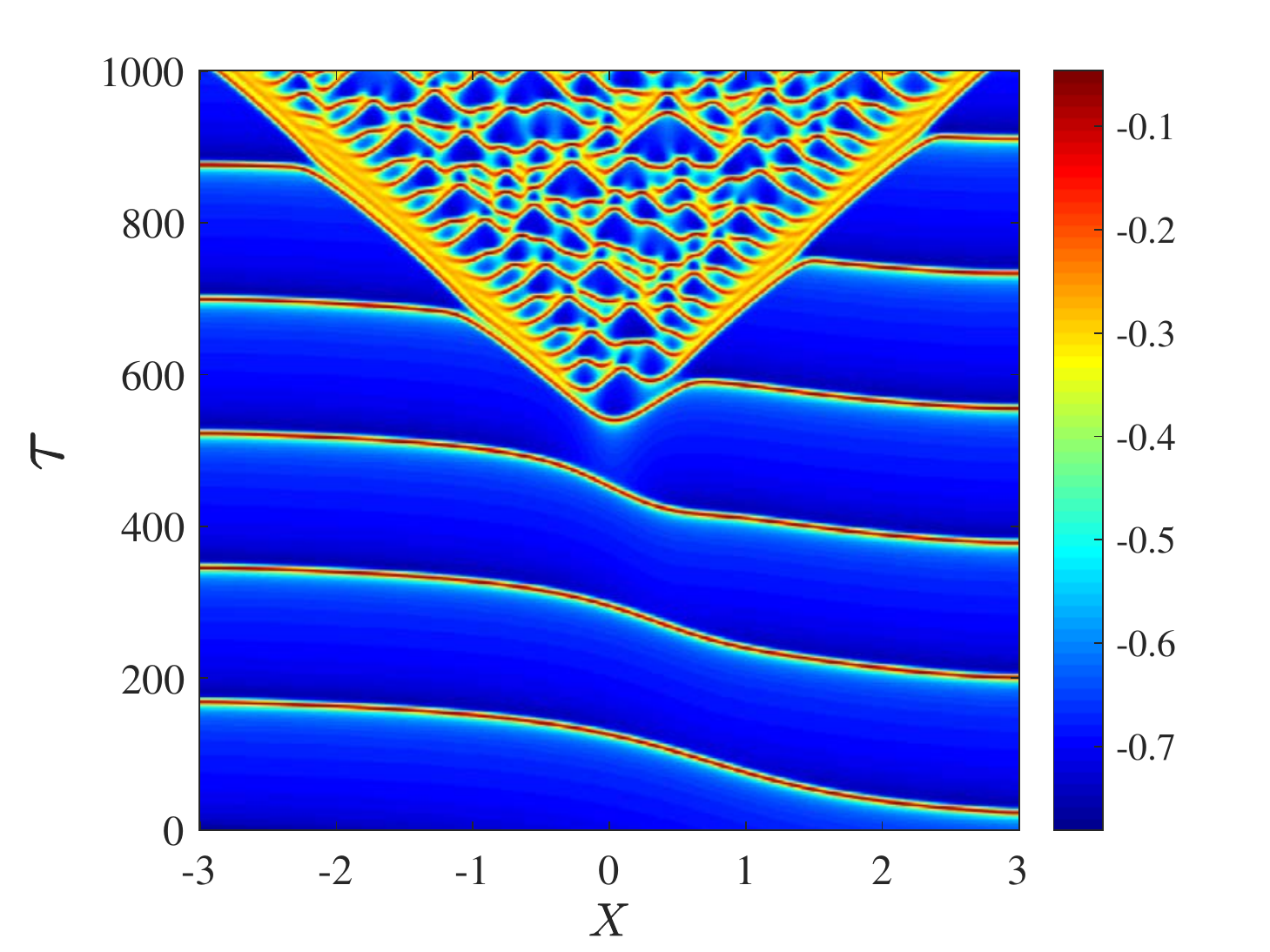}
    \label{fig:Asym-spt025v1}
  \end{subfigure}%
  \begin{subfigure}[b]{.35\linewidth}
    \centering
    \caption{}
    \includegraphics[width=.99\textwidth]{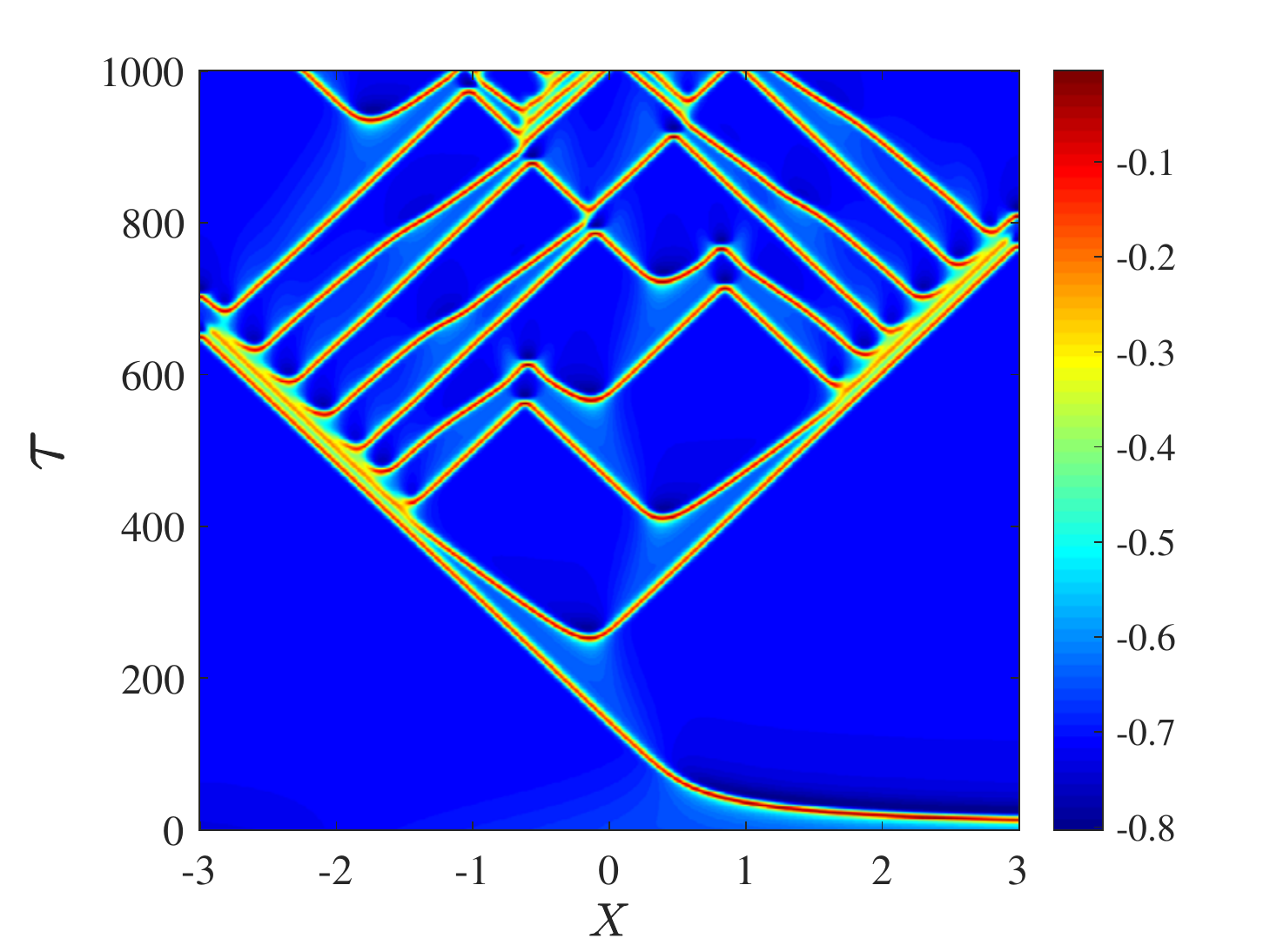}
     \label{fig:Asym_spt012psi}
  \end{subfigure}%
   \caption{Space-time plots of the membrane potential $V$ using a longer time-frame than other plots.  The parameter values and initial conditions in panel (a) are the same as Fig.~\ref{fig:spt0245}, and in panel (b) are same as Fig.~\ref{fig:sptpsi012}.}
    \label{fig:Asym_spt}
  \end{figure}

In this section we consider other perturbation functions $G(X)$ in the initial condition \eqref{eq:init} to investigate how the initial condition affects the patterns that develop. First we consider
\begin{equation}
\label{eq:init2}
    G(X)=\epsilon X,
\end{equation}
with $\epsilon=0.025$. Fig.~\ref{fig:SPT_v1_X} shows the resulting spatiotemporal patterns for different values of $\bar{v}_1$.  Specifically the six plots use the same parameter values as the corresponding plots in Fig.~\ref{fig:SPT_v1}. In panels (a) and (f) of Fig.~\ref{fig:SPT_v1_X} the solution simply settles to the stable equilibrium of the system in the absence of diffusion (as in Fig.~\ref{fig:SPT_v1}). In panels (b) and (c) the initial condition is insufficient to generate the spatiotemporal chaos that was observed in Fig.~\ref{fig:SPT_v1} within the $\tau=500$ timeframe. By simulating for a longer time we found that in (b) the solution appeared to converge to homogeneous oscillations matching the stable limit cycle of the system in the absence of diffusion,
while in (c) spatiotemporal chaos did arise shortly after $\tau = 500$, and this is shown in Fig.~\ref{fig:Asym-spt025v1}. Finally in panels (d) and (e) we do observe spatiotemporal chaos. The particular patterns that emerge appear to have the same features as those in Fig.~\ref{fig:SPT_v1} suggesting that for both initial conditions the solution is converging to the same attractor.

For other values of the parameters and other initial conditions we similarly observed that, broadly speaking, the dynamics of \eqref{eq:dimless1}--\eqref{eq:dimless2} settled to the same long-time behaviour as that described in Sect.~\ref{sec:SPTD}. For example using the parameter values of Fig.~\ref{fig:sptpsi012}, when the initial condition is changed from \eqref{eq:gauss} to \eqref{eq:init2} the result is Fig.~\ref{fig:Asym_spt012psi} which evidently exhibits a similar structure. We conclude that the profile of the initial perturbation does not seem to change the types of spatiotemporal patterns that are produced by the model.
  
  \section{Travelling wave analysis}\label{sec:6}
\begin{figure}[htbp]
\centering
   \begin{subfigure}[b]{.3\linewidth}
    \centering
     \caption{}
    \includegraphics[width=.99\textwidth]{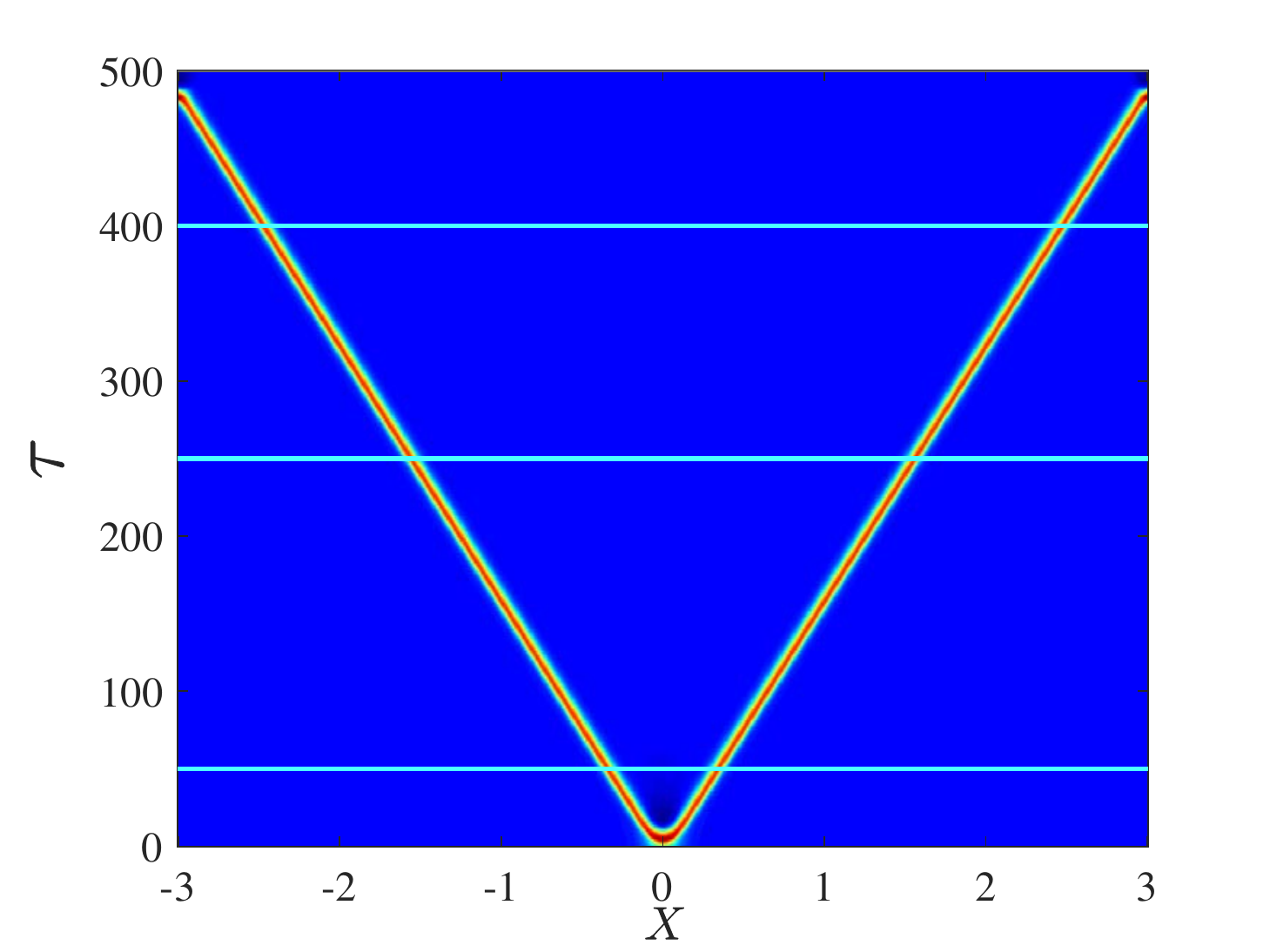}
  \end{subfigure}\\%
  \begin{subfigure}[b]{.3\linewidth}
    \centering
     \caption{}
    \includegraphics[width=.99\textwidth]{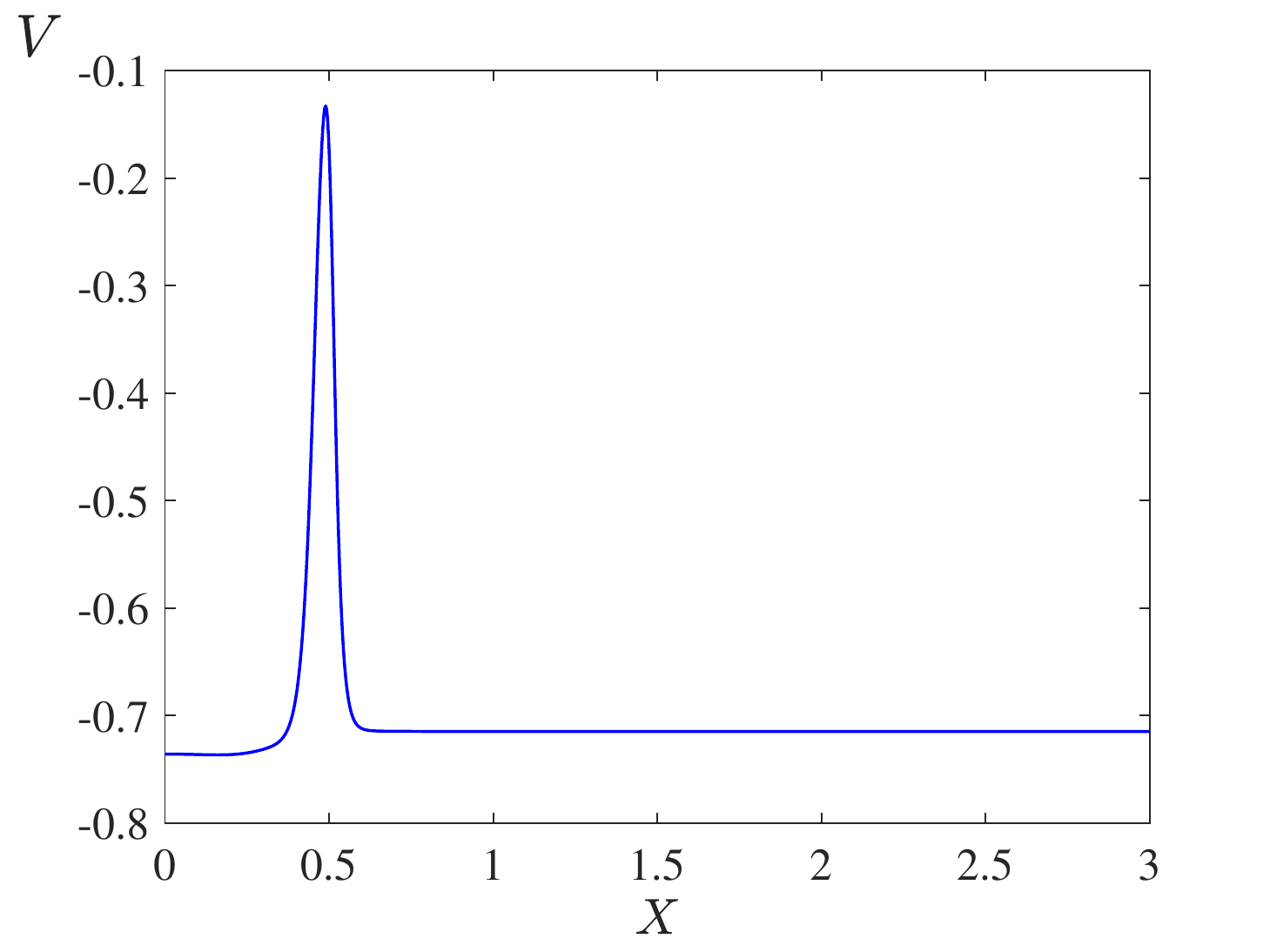}
  \end{subfigure}%
  \begin{subfigure}[b]{.3\linewidth}
    \centering
     \caption{}
    \includegraphics[width=.99\textwidth]{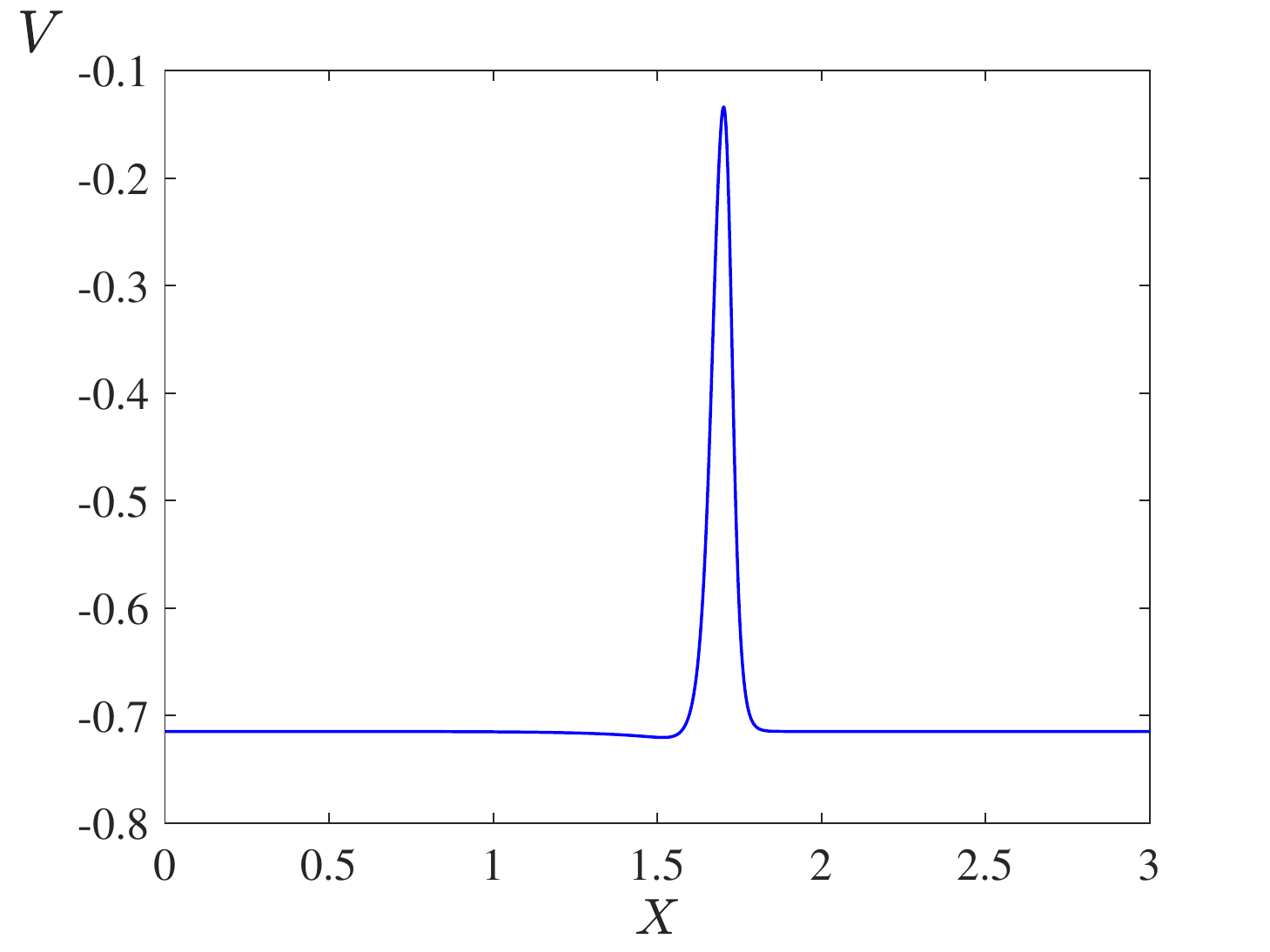}
  \end{subfigure}%
  \begin{subfigure}[b]{.3\linewidth}
    \centering
     \caption{}
    \includegraphics[width=.99\textwidth]{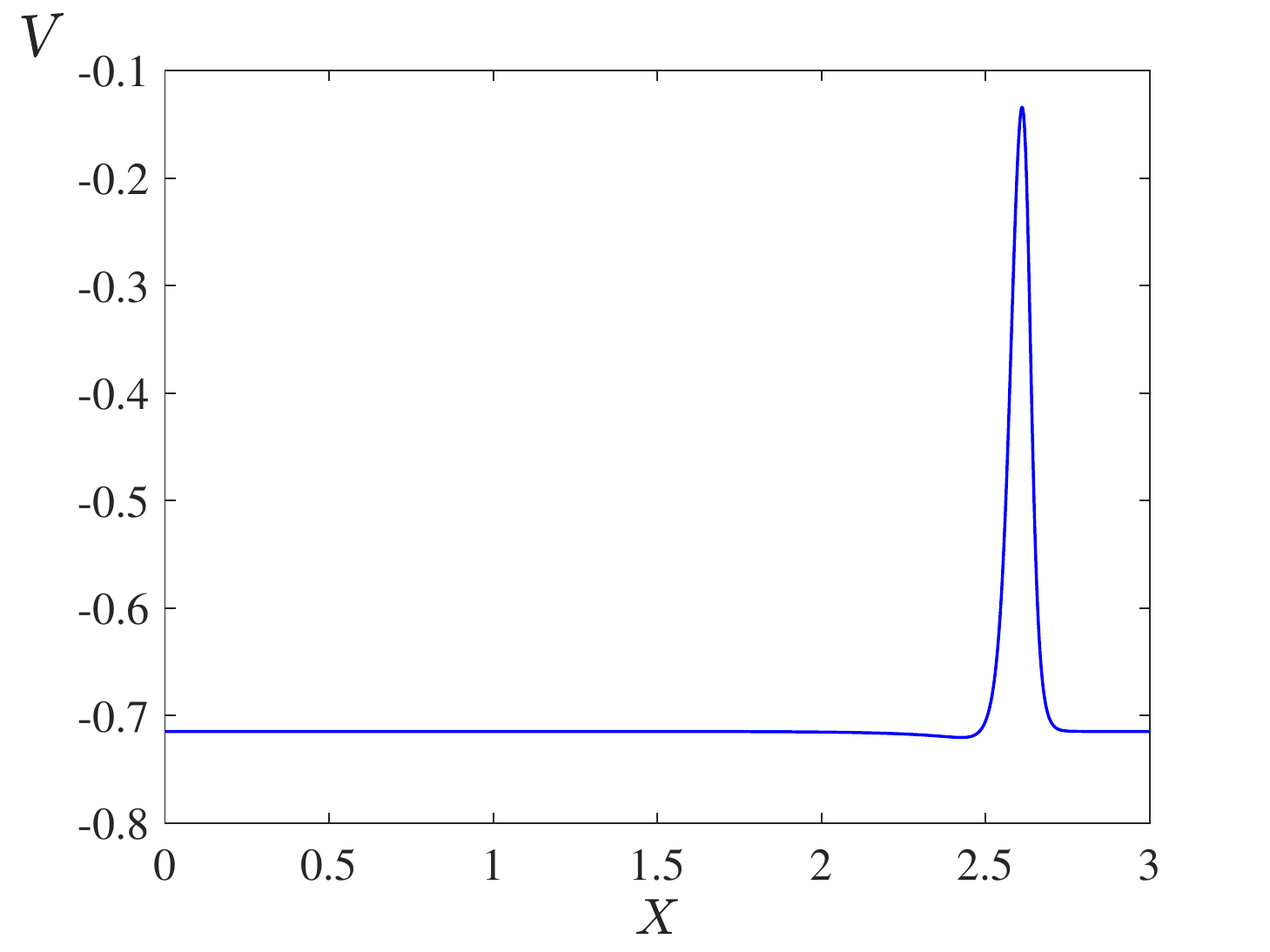}
  \end{subfigure}%
 \caption{A reproduction of Fig.~\ref{fig:SPT_psi}a and plots of the solution profile at the values of $\tau$ that are marked by horizontal lines.}
 \label{fig:propagating pulse}
\end{figure}

\begin{figure}[htbp]
\centering
   \begin{subfigure}[b]{.3\linewidth}
    \centering
     \caption{}
    \includegraphics[width=.99\textwidth]{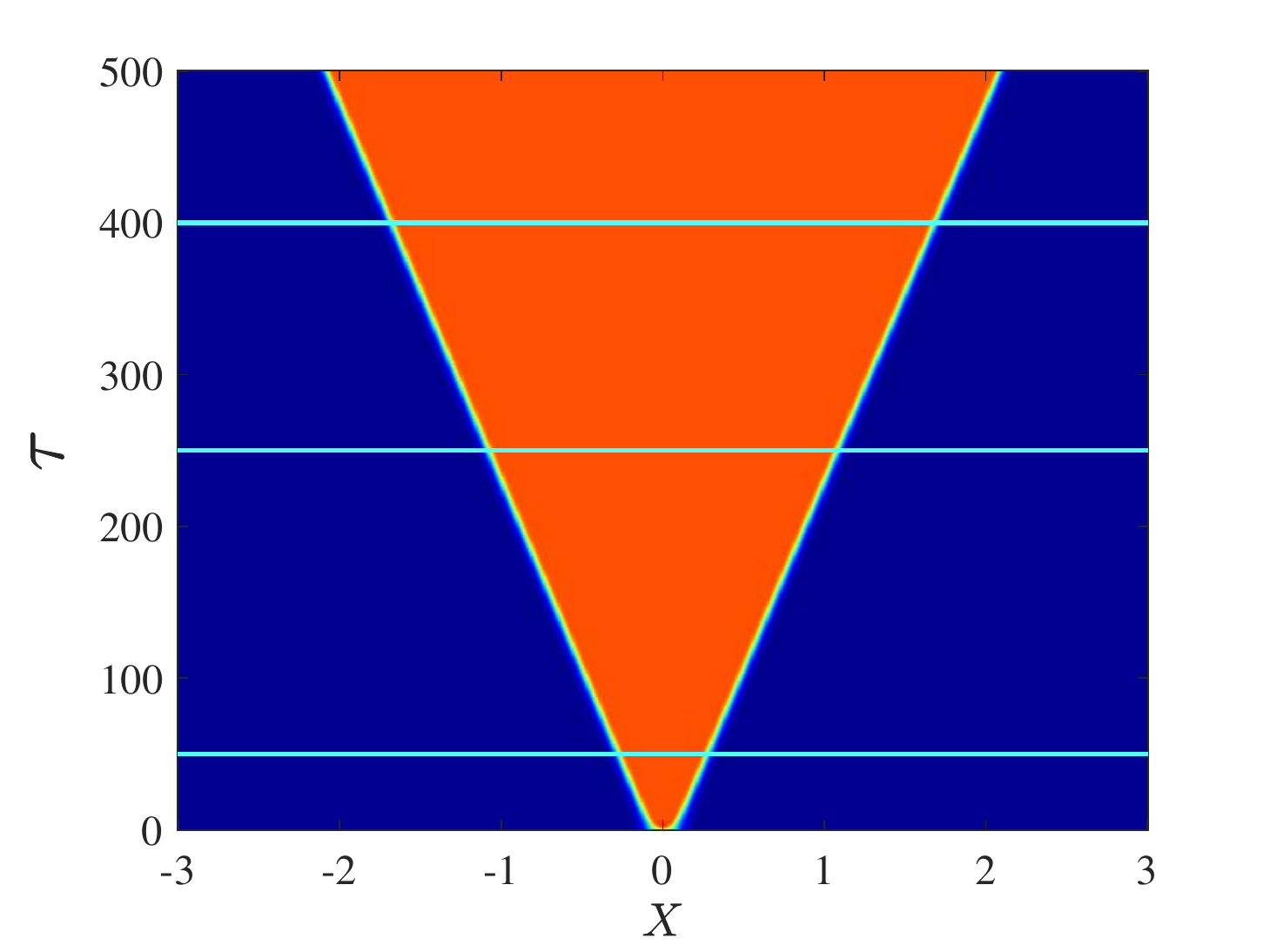}
  \end{subfigure}\\%
  \begin{subfigure}[b]{.3\linewidth}
    \centering
     \caption{}
    \includegraphics[width=.99\textwidth]{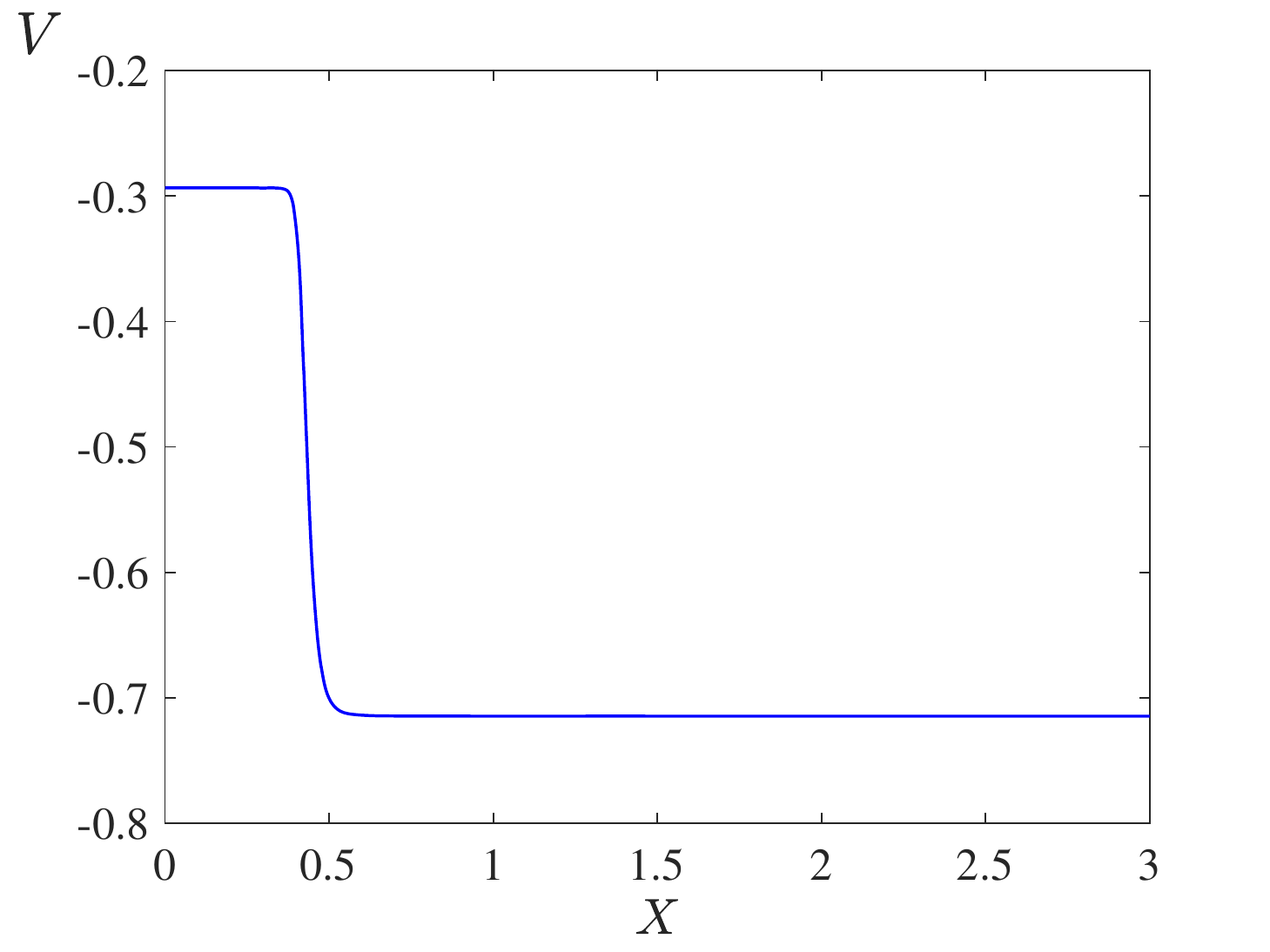}
  \end{subfigure}%
  \begin{subfigure}[b]{.3\linewidth}
    \centering
     \caption{}
    \includegraphics[width=.99\textwidth]{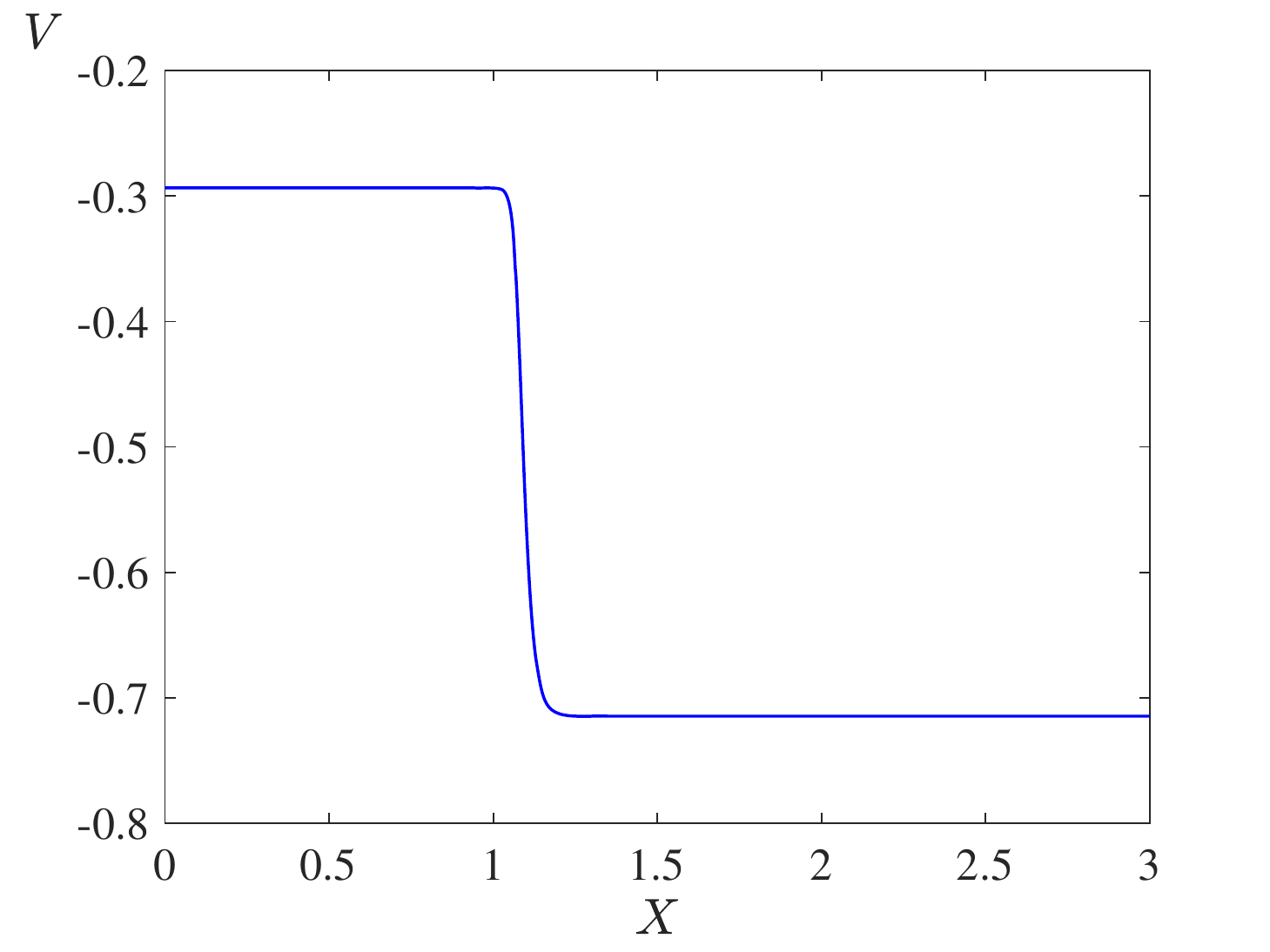}
  \end{subfigure}%
  \begin{subfigure}[b]{.3\linewidth}
    \centering
     \caption{}
    \includegraphics[width=.99\textwidth]{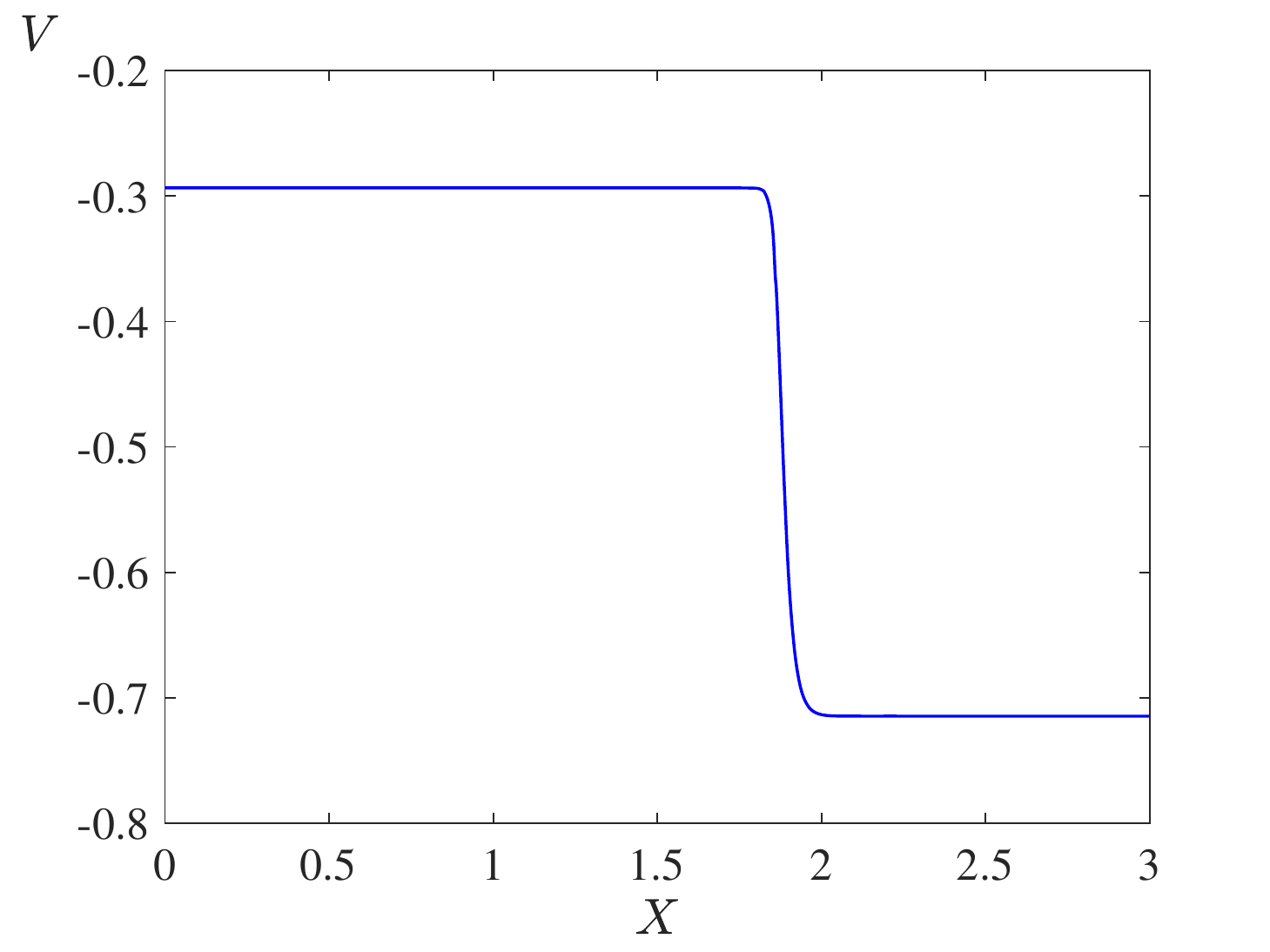}
  \end{subfigure}%
 \caption{A reproduction of Fig.~\ref{fig:SPT_psi}f and plots of the solution profile at the values of $\tau$ that are marked by horizontal lines.}
 \label{fig:propagating front}
\end{figure}

  \setcounter{equation}{0}
For the travelling waves analysis we will focus on the values of $\psi$ where the numerical simulations of \eqref{eq:dimless1}--\eqref{eq:dimless2} result in travelling pulses and fronts, respectively. For example, when $\psi=0.1$ two stable counter-propagating pulses are created, and they travel across the domain at approximate speed $c=0.006182$ (see Fig.~\ref{fig:sptpsi01}).  Fig.~\ref{fig:propagating pulse}a shows the pulses and Fig.~\ref{fig:propagating pulse}b--d are solution profiles at times $\tau=50$, $250$, $400$. Also, when $\psi=0.5$ two stable counter-propagating fronts are created, and they travel across the domain at speed $c=0.004155$ (see Fig.~\ref{fig:sptpsi05}). The fronts are shown in Fig.~\ref{fig:propagating front}a and Fig.~\ref{fig:propagating front}b--d are solution profiles at the same three times. The given wave speeds have been estimated directly from the numerical simulation results.

In the coming section, we introduce the travelling wave variable to transform \eqref{eq:dimless1}--\eqref{eq:dimless2} to a set of three ODEs and approximate the travelling wave solutions numerically. This allows us to find the homoclinic and heteroclinic trajectories that correspond to the travelling pulse and front solutions, respectively. 

\subsection{Existence of Travelling Waves}
To describe the travelling wave profile we consider travelling waves with unknown wave speed $c>0$. By introducing the travelling wave variable, $\zeta= X-c\tau$, the model \eqref{eq:dimless1}--\eqref{eq:dimless2} becomes
\begin{equation}
    \begin{pmatrix}V\\N
    \end{pmatrix}_{\tau}=D\begin{pmatrix}V\\0 \end{pmatrix}_{\zeta\zeta}+\begin{pmatrix}cV\\cN\end{pmatrix}_{\zeta}+\begin{pmatrix}f(V,N)\\g(V,N)\end{pmatrix},
    \label{eq:TWS_matrix}
\end{equation}
where
\begin{align*}
f(V,N)&=-\bar{g}_{L}(V-\bar{v}_{L})-\bar{g}_{K}N(V-\bar{v}_{K})-\bar{g}_{\text{Ca}}M_{\infty}(V)(V-\bar{v}_{\text{Ca}}),\\
g(V,N)&=\lambda_N(V)\big(N_{\infty}(V)-N\big).
\end{align*}
Travelling waves are stationary solutions to \eqref{eq:TWS_matrix} and satisfy
\begin{equation}
    D\begin{pmatrix}V\\0 \end{pmatrix}_{\zeta\zeta}+c\begin{pmatrix}V\\N\end{pmatrix}_{\zeta}+\begin{pmatrix}f(V,N)\\g(V,N)\end{pmatrix}=0.
    \label{eq:TWODE_matrix}
\end{equation}
We rewrite \eqref{eq:TWODE_matrix} as system of first order ODEs with $':=\frac{d}{d\zeta}$ by introducing a new variable $W=V'$ to obtain
\begin{equation}
\begin{aligned}
V'&=W,\\
W'&=-\frac{1}{D}(cW+f(V,N)),\\
N'&=-\frac{1}{c}g(V,N).
\end{aligned}
\label{eq:TWODE}
\end{equation}
The boundary conditions are
\begin{align}
\label{eq:boundaryvalues}
 \lim_{\zeta\to +\infty}(V,W,N)(\zeta)=(V_{+},0,N_{+}), \hspace{2mm} \lim_{\zeta\to -\infty}(V,W,N)(\zeta)=(V_{-},0,N_{-}),
\end{align} 
where $(V_{\pm},N_{\pm})$ are equilibria of \eqref{eq:dimless1}--\eqref{eq:dimless2}. For a pulse $(V_+,N_+)$ and $(V_-,N_-)$ are the same equilibrium; for a front they are different equilibria.

A number of mathematical methods have been established to show the existence of travelling wave solutions in reaction-diffusion systems. These involve singular perturbation theory \cite{Merkin1996,cornwell2018}, variational techniques \cite{Chen2015}, and factorisation \cite{Achouri2016}. We use the shooting method \cite{Ermentrout2002SimulatingStudents} to identify travelling waves and approximate their wave speed. This was achieved by numerically computing solutions to the travelling wave ODEs (17) for initial points perturbed from an equilibrium in a direction tangent to either its stable manifold or unstable manifold. In either case this direction is given by an eigenvector of the Jacobian matrix of \eqref{eq:TWODE} evaluated at the equilibrium, and a formula for this matrix is provided in Appendix A.  We adjusted the value of $c$ until the solution was approximately homoclinic (in the case of a pulse) or heteroclinic (in the case of a front).

\begin{figure}[htbp]
\centering
  \begin{subfigure}[b]{.5\linewidth}
    \centering
    \caption{}
    \includegraphics[width=0.99\textwidth]{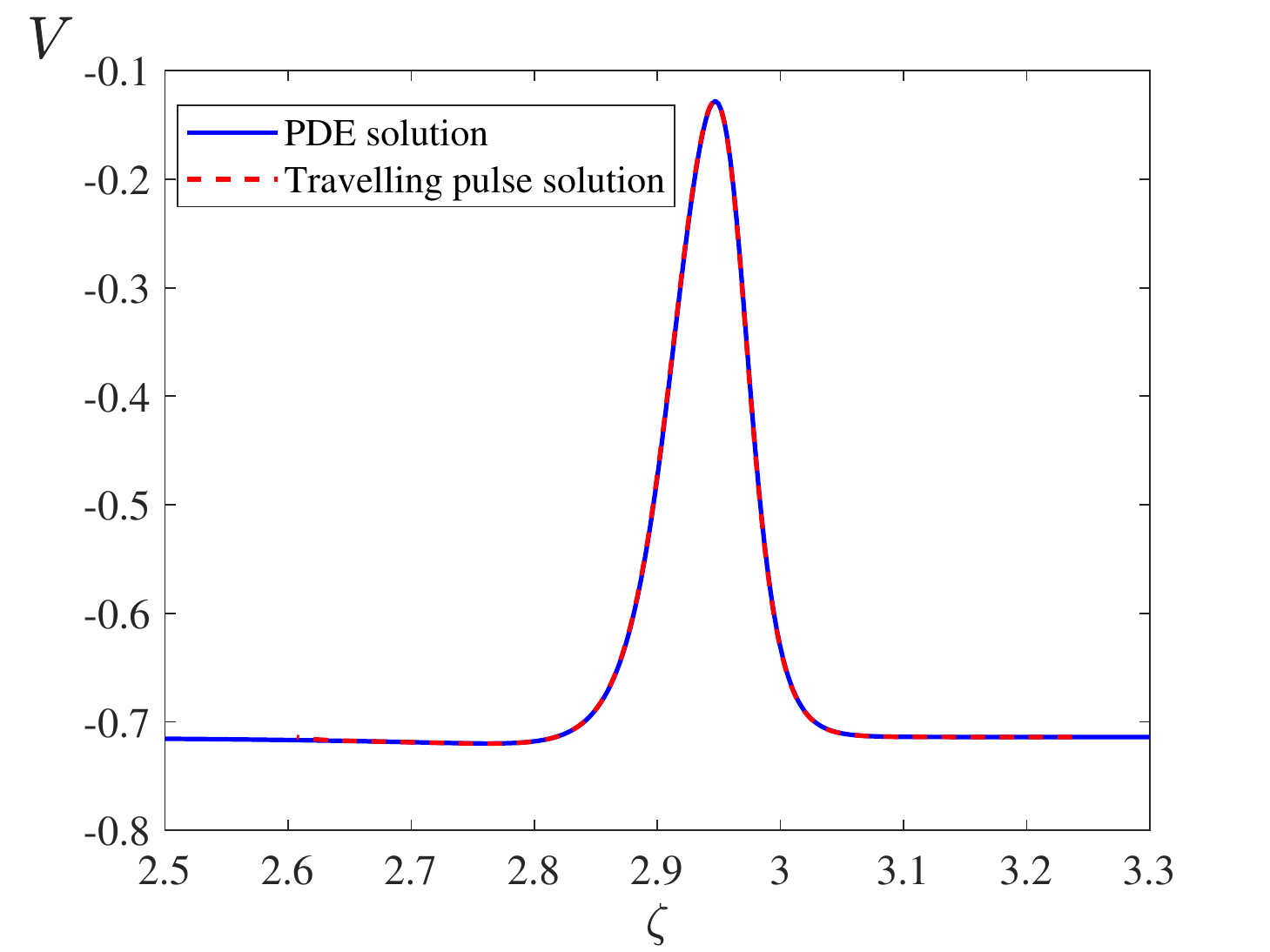}
    \label{fig:pulse}
  \end{subfigure}%
  \begin{subfigure}[b]{.5\linewidth}
    \centering
    \caption{}
    \includegraphics[width=0.99\textwidth]{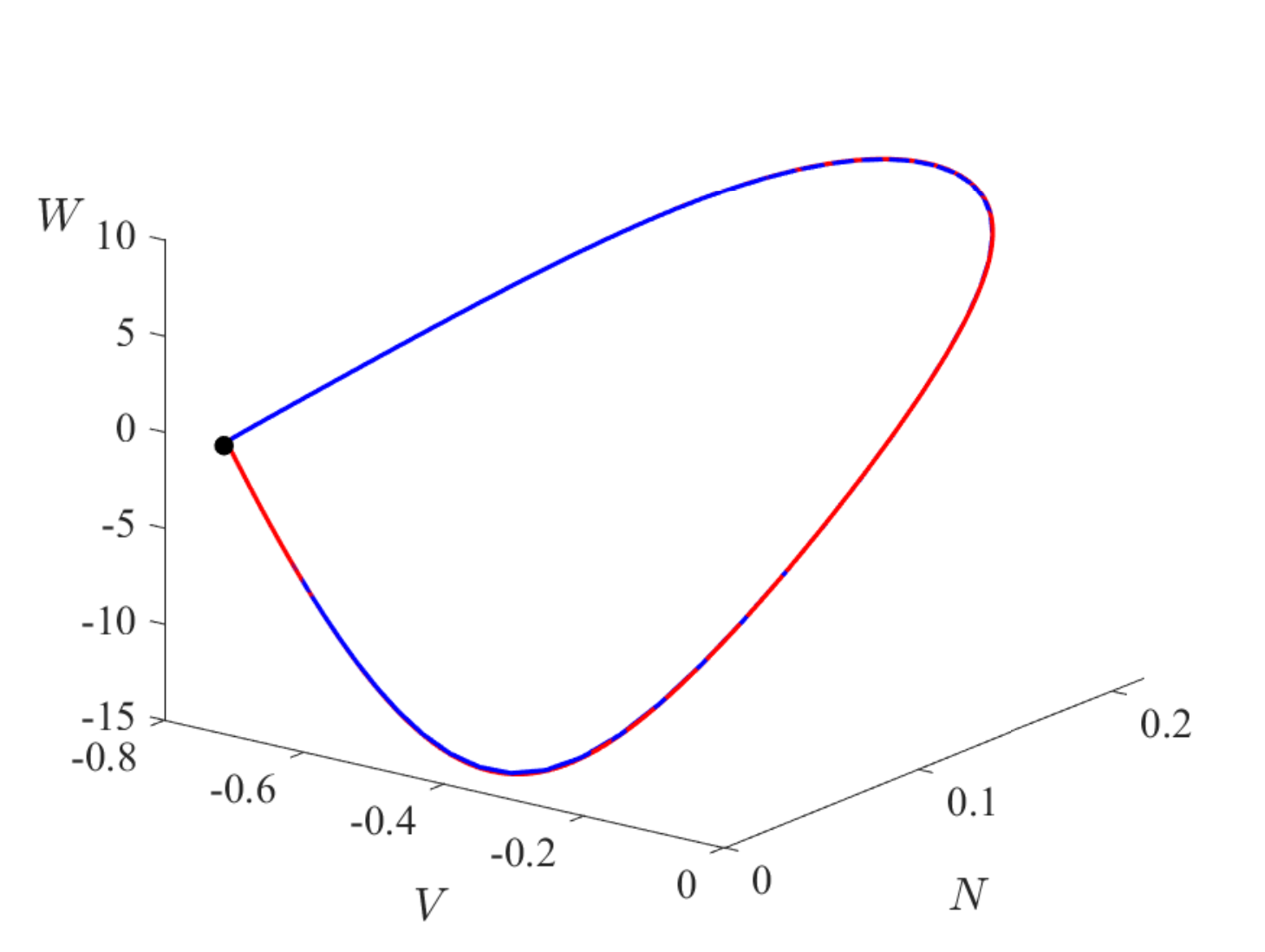}
     \label{fig:Hom_pulse}
  \end{subfigure}%
  \caption{(a) The solution profile of \eqref{eq:dimless1}--\eqref{eq:dimless2} and a solution to the travelling wave ODEs \eqref{eq:TWODE} with $c=0.006116$, using the same parameter values as Fig.~\ref{fig:SPT_psi}a (b) The same two solutions but plotted in the phase space of \eqref{eq:TWODE}.}
  \label{fig:pulse_hom}
  \end{figure}

We first consider the parameter values of Fig.~\ref{fig:sptpsi01} for which stable travelling pulses were observed.  The equilibrium associated with these pulses is the lower-most equilibrium branch of Fig.\ref{fig:onepar}c.  For the travelling wave ODEs \eqref{eq:TWODE} this equilibrium has a one-dimensional unstable manifold. By performing the shooting method we found that a solution approximating one branch of this manifold forms a homoclinic connection when $c = 0.006116$, approximately. This matches the speed of the pulses observed in Fig.~\ref{fig:sptpsi01}. A plot of the pulse profile for $V$ is shown in Fig.~\ref{fig:pulse} and its corresponding homoclinic trajectory in $(V,W,N)$ phase space is shown in Fig.~\ref{fig:Hom_pulse}. As expected the pulse profile extracted from our numerical solution to \eqref{eq:dimless1}--\eqref{eq:dimless2} matches the pulse solution obtained of the travelling wave ODEs \eqref{eq:TWODE}.
  
\begin{figure}[htbp]
\centering
  \begin{subfigure}[b]{.5\linewidth}
    \centering
    \caption{}
    \includegraphics[width=0.99\textwidth]{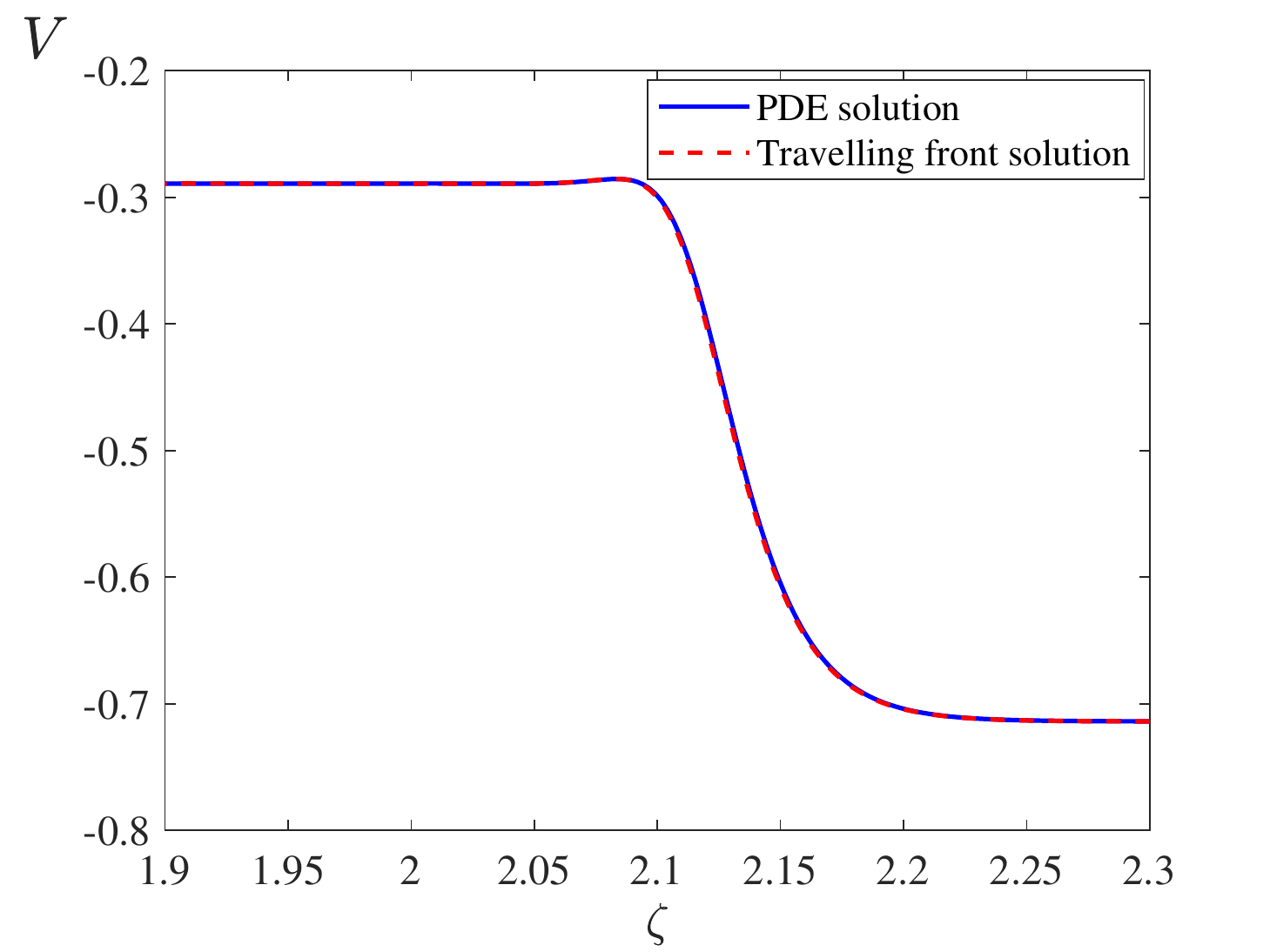}
    \label{fig:front_sol}
  \end{subfigure}%
  \begin{subfigure}[b]{.5\linewidth}
    \centering
    \caption{}
    \includegraphics[width=0.99\textwidth]{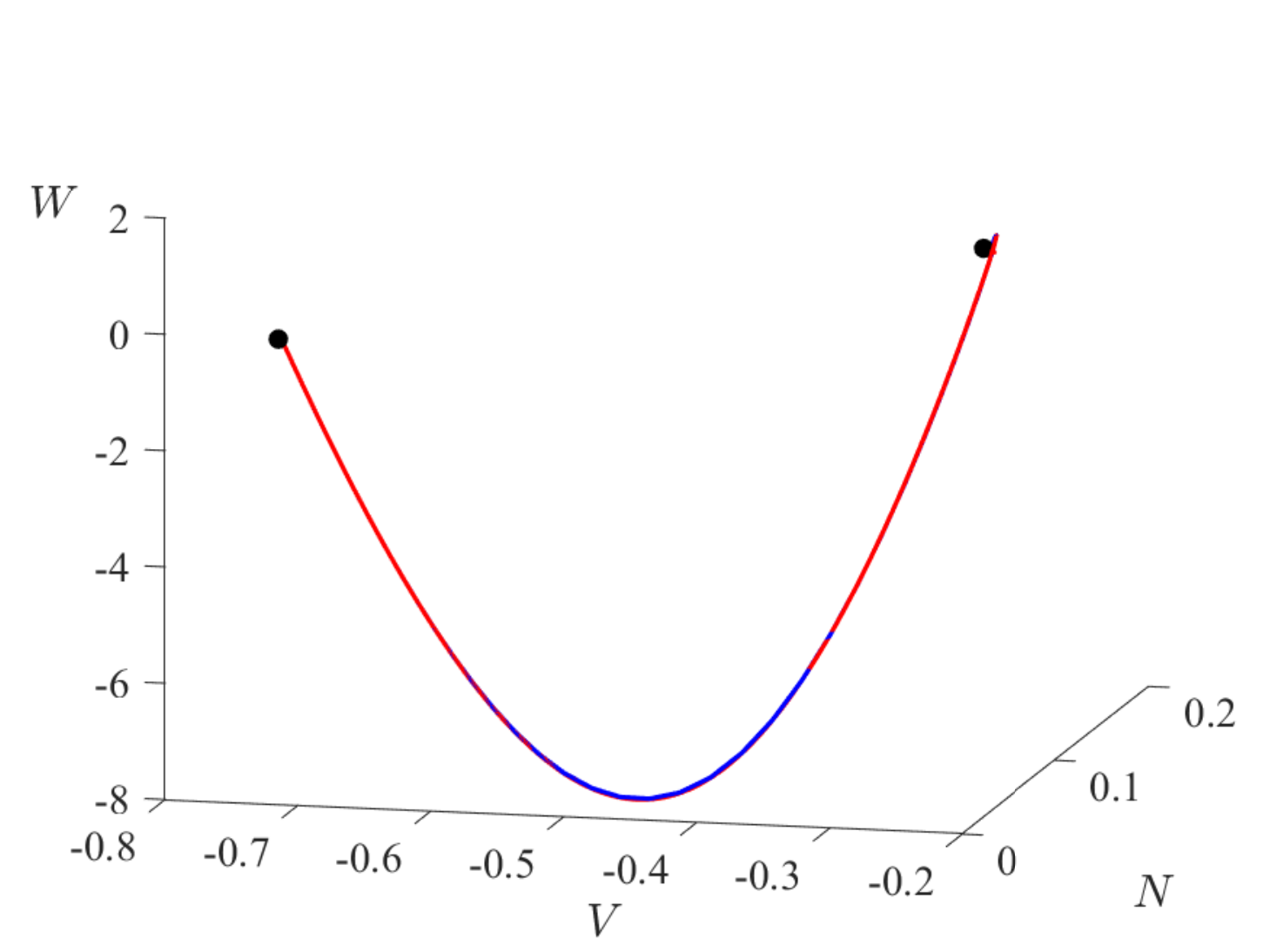}
     \label{fig:Het_front}
  \end{subfigure}%
  \caption{(a) The solution profile of \eqref{eq:dimless1}--\eqref{eq:dimless2} and a solution to the travelling wave ODEs \eqref{eq:TWODE} with $c=0.0043$ using the same parameter values as Fig.~\ref{fig:SPT_psi}f (b) The same two solutions but plotted in the phase space of \eqref{eq:TWODE}.}
  \label{fig:front}
  \end{figure}

  Now we consider the parameter values of Fig.~\ref{fig:SPT_psi}f for which our numerical solution produced two travelling fronts.  These connect the lower-most and upper-most equilibrium branches of Fig.\ref{fig:onepar}c.  As equilibria of \eqref{eq:TWODE} these have one-dimensional stable manifolds.  Consequently we solved \eqref{eq:TWODE} backwards in time from an initial point near the upper equilibrium and adjusted the value of $c$ until observing an approximately heteroclinic orbit.  This produced $c = 0.0043$, approximately, matching the wave speed observed in Fig.~\ref{fig:SPT_psi}f. The plot of the front profile for $V(\zeta)$ is shown in Fig.~\ref{fig:front_sol} and its corresponding heteroclinic trajectory in $(V,W,N)$ phase space is shown in Fig.~\ref{fig:Het_front}. The front obtained by the solution to \eqref{eq:dimless1}--\eqref{eq:dimless2} is also shown and seen to closely match the front profile of \eqref{eq:TWODE}.

\section{Conclusion}

In this work we have studied the collective dynamics of pacemaker SMCs through passive electrical coupling of adjacent cells. We presented a detailed analysis of spatiotemporal dynamics of the model in one-dimensional domain. The local dynamics of the reaction-diffusion system \eqref{eq:dimless1}--\eqref{eq:dimless2}
under variation of model parameters are analysed to study the electrical activity of an isolated SMC. We examined the dynamics of the membrane voltage $V$ via numerical bifurcation analysis with the potential at which the potassium and calcium channels are half-opened, $\bar{v}_{1}$ and $\bar{v}_{3}$, and the time constant, $\psi$, as bifurcation parameters. Upon varying the parameters the model exhibits different dynamical features including transitions from stable solutions to periodic firing. The model exhibited the two types of excitability for excitable cells depending on how parameters are varied. Variation of $\bar{v}_{1}$ results in Type I excitability and Type II excitability is observed as $\bar{v}_{3}$ is varied \cite{hammed}. 

With the help of linear stability analysis, we showed that the patterns are not due to Turing instability. In this present work, we have that the spatiotemporal patterns observed in the model are as a consequence of the nonlinear dynamics of the system in the absence of diffusion. That is the conditions for instability in the spatially extended model is the same as that of the non-spatial system.       

The main goal of this work is to examine spatiotemporal dynamics in a network of SMCs through propagation of action potential between adjacent cells. To investigate the evolution of the model dynamics as time increases we perturbed the centre of the domain with a Gaussian pulse. In general, the initial perturbation induced counter propagating pulses that travel across the domain as time progresses. The numerical simulations of the reaction-diffusion system showed a wide variety of spatiotemporal behaviours such as travelling pulses, travelling fronts and spatiotemporal chaos. Transitions between patterns are driven by the model parameters. For instance, the system transitions from stable travelling pulses to travelling fronts as the rate constant for the kinetics of the open ${\rm K}^{+}$ channel $\psi$ is varied. 

We also showed that the shape of the initial perturbation does not affect the type of spatiotemporal patterns exhibited by the model. The patterns observed when the initial perturbation is changed to a straight line are similar to when a Gaussian pulse is considered as the perturbation function. Therefore, we can conclude that the shape of initial perturbation does not appear to have much effect on different patterns that emerge as parameters are varied  provided an action potential is triggered by the initial perturbation.

We also investigated the travelling wave solutions to the reaction-diffusion system. We showed the existence of homoclinic trajectory connecting a steady state to itself (pulse solution) and heteroclinic trajectory connecting two different steady states (front solution) numerically with the shooting method, and it was shown that the estimated wave speed in the travelling wave ODE system is very close to the wave speed observed in the numerical simulations. It is worth exploring the stability of the travelling waves, this would provide a better understanding of how the spatiotemporal patterns transition from stable propagating waves to complex patterns. This analysis will be our focus in future work. In this work we have studied spatiotemporal dynamics in one spatial dimension, however the model can be extended to two-dimension setting. Hence, it would be interesting to further study two-dimensional patterns in the reaction-diffusion system.

\appendix
\section{The linearisation of the travelling wave ODE}\label{sec:A1}
Here we evaluate the Jacobian matrix of the travelling wave ODEs (17) at an arbitrary equilibrium $(V_{\pm},W_{\pm},N_{\pm})$. The Jacobian matrix is
\begin{align}
\label{eq:asym_pm}
J=\begin{pmatrix}0&1&0\\-\frac{1}{D}f_{V}&-\frac{c}{D}&-\frac{1}{D}f_{N}\\-\frac{1}{c}g_{V}&0&-\frac{1}{c}g_{N}\end{pmatrix},
\end{align} 
where 
\begin{align}
 f_{V}&=\Bigg[-\bar{g}_L-\bar{g}_{K}N_{\pm}-\frac{\bar{g}_{\text{Ca}}}{2\bar{v}_2}\left(1-\tanh^2\left(\frac{V_{\pm}-\bar{v}_1}{\bar{v}_2}\right)\right)(V_{\pm}-\bar{v}_{\text{Ca}}) \nonumber\\ & \hspace{6cm} -\frac{\bar{g}_\text{Ca}}{2}\left(1+\tanh\left(\frac{V_{\pm}-\bar{v}_{1}}{\bar{v}_{2}}\right)\right)\Bigg], \nonumber\\
f_{N}&=-\bar{g}_K(V_{\pm}-\bar{v}_K), \nonumber\\
g_{V} & =\frac{\psi}{2\bar{v}_4}\Bigg[ \Bigg\{\frac{1}{2}\left(1+\tanh\left(\frac{V_{\pm}-\bar{v}_{3}}{\bar{v}_{4}}\right)\right)-N_{\pm}\Bigg\}  \sinh\left(\frac{V_{\pm}-\bar{v}_{3}}{2\bar{v}_{4}}\right)\Bigg]\nonumber \\&\hspace{5cm} +\frac{\psi}{2\bar{v}_4} \Bigg[\cosh\left(\frac{V_{\pm}-\bar{v}_{3}}{2\bar{v}_{4}}\right)\left(1-\tanh^2\left(\frac{V_{\pm}-\bar{v}_3}{\bar{v}_4}\right)\right)\Bigg], \nonumber\\
g_{N}& =-\psi \cosh\left(\frac{V_{\pm}-\bar{v}_{3}}{2\bar{v}_{4}}\right). \nonumber
\end{align}
The eigenvalues of \eqref{eq:asym_pm} are the solutions to the characteristic equation
\begin{align}
\label{eq:cubicpoly}
    \lambda^{3}+P_{2}\lambda^{2}+P_{1}\lambda+P_{0}=0,
\end{align}
where
$$P_{2}=\frac{g_{N}D+c^{2}}{cD}, \hspace{0.3em}  P_{1}=\frac{g_{N}+f_{V}}{D}, \hspace{0.3em} \text{and} \hspace{0.3em}  P_{0}=\frac{f_{V}g_{N}-f_{N}g_{V}}{cD}.
$$

  



\end{document}